\documentclass{article}
\newcommand*\Sim{\includegraphics[scale=0.26]{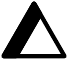}}
\newcommand*\etaSim{\includegraphics[scale=0.26]{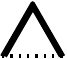}}
\newcommand{\rot}[1]{#1^{\mathrm{rot}}}
\usepackage{graphicx} %
\usepackage{biblatex}
\usepackage{tikz}
\usepackage{enumerate}
\usetikzlibrary{patterns} %

\usepackage{float}

\usepackage{microtype}
\usepackage{xurl}

\newcommand{\Rel}{\mbox{\bf Rel}}

\newcommand{\Hom}{\mbox{\rm Hom}}

\newcommand{\eSSet}{\mathrm{Set}^\varepsilon_\Delta}
\newcommand{\neSSet}{\mathrm{Set}^{\eta\varepsilon}_\Delta}
\newcommand{\SSet}{\mathrm{Set}_\Delta}
\newcommand{\eHorn}{\varepsilon-\mathrm{Horn}}

\DeclareRobustCommand{\ELambda}{\text{\reflectbox{$\bf\Lambda$}}}
\def\relto{\mapstochar\dashrightarrow}

\tikzset{LA/.style = {%
                      line width=#1, -{Straight Barb[length=3pt]}},
         LA/.default=1pt
        }

\usepackage{amsmath,amsthm,amssymb,cases,mathtools}

\usepackage{mathabx} %

\usepackage{tikz-cd}
\tikzset{%
    symbol/.style={%
        draw=none,
        every to/.append style={%
            edge node={node [sloped, allow upside down, auto=false]{$#1$}}}
    }
}
\tikzcdset{scale cd/.style={every label/.append style={scale=#1},
    cells={nodes={scale=#1}}}}

\title{Stability theorems for \(\varepsilon\)-simplicial sets}
\author{Dominik Lachman%
\thanks{Department of Algebra and Geometry,
Faculty of Science, Palack\'y University Olomouc,
17.~listopadu~12, 779~00~Olomouc, Czech Republic.}%
\kern0.35em\textsuperscript{,}%
\thanks{The author received support from the Czech Science Foundation grant 25-20013L.}}
\date{}

\theoremstyle{plain}
\newtheorem{theorem}{Theorem}[section]
\newtheorem{lemma}[theorem]{Lemma}
\newtheorem{corollary}[theorem]{Corollary}
\newtheorem{proposition}[theorem]{Proposition}

\theoremstyle{definition}
\newtheorem{defn}[theorem]{Definition}
\newtheorem{example}[theorem]{Example}

\theoremstyle{remark}

\begin{document}

\maketitle
\begin{abstract}
We establish stability results for internal homs of
$\varepsilon$-simplicial sets, a framework motivated by
nerves of Frobenius monoids in the category of relations.
Rotations are modelled through edge witnesses and
horn-filling conditions, without requiring unique fillers.
For $\varepsilon$-fibrant objects $X$ and $Y$, we prove
that the internal hom $\mathrm{Fun}(X,Y)$ satisfies
the $\varepsilon$-coherence conditions and admits fillers
for all $\varepsilon$-horns of dimension at least two.
Its subobject consisting of simplices whose vertices
correspond to $\eta$-edge-preserving maps is
$\varepsilon$-fibrant.
We also introduce $\eta\varepsilon$-simplicial sets
by adding explicit $\eta$-witnesses and prove that
the internal hom of two $\eta\varepsilon$-fibrant
objects is again $\eta\varepsilon$-fibrant.

Our proofs use adjunctions relating
$\varepsilon$-simplicial sets and paracyclic sets.
We characterise the $\varepsilon$-horn filling property
by requiring the corresponding counit to be a trivial
Kan fibration.
Consequently, for every object satisfying this property,
coherent choices of rotations of each simplex form
a contractible Kan complex.
\end{abstract}

\section{Introduction}

A theorem of Joyal states that, for any simplicial sets $X$ and $Y$,
if $Y$ is a quasicategory, then the internal hom $\mathrm{Fun}(X,Y)$
is also a quasicategory~\cite[Corollary~2.19]{Joyal2008}.
Here a quasicategory is a simplicial set admitting fillers for all
inner horns. In this article, we establish analogues of this stability
theorem for $\varepsilon$-fibrant $\varepsilon$-simplicial sets,
with particular attention to one-dimensional horn conditions.
The underlying notion of an $\varepsilon$-simplicial set was introduced
in~\cite{Lachman2026} to study Frobenius monoids in the category
$\Rel$ of sets and relations, encoding their rotational structure
through edge witnesses and lifting properties.

Frobenius monoids in $\Rel$ include structures arising in several
areas of mathematics. Groupoids~\cite{HCC2013},
effect algebras~\cite{PS,Lachman2026}, and
effect algebroids~\cite{Roumen2017,Lachman2026} give examples,
while each coherent configuration induces a relational Frobenius
monoid on its set of colours~\cite{JJL2025,Lachman2026}.
Interest in Frobenius structures is also motivated by their role
in categorical quantum mechanics~\cite{HV} and two-dimensional
topological quantum field theory~\cite{Kock2003}.

A simplicial approach to relational Frobenius monoids was initiated
by Mehta and Zhang~\cite{MZ}, who described them using simplicial
sets equipped with compatible rotational data.
Contreras, Keller, and Mehta~\cite{CKM} developed a related
description of Frobenius objects in the category of spans.

The coherent setting connects these constructions with the theory
of $2$-Segal spaces~\cite{DK2019}, in which polygonal decompositions
encode associativity. In particular, $2$-Segal sets describe
pseudomonoids in the bicategory of spans~\cite{Stern2021,MM2025Segal}.
Contreras, Mehta, and Stern~\cite{CMS} showed that Frobenius
structures on these pseudomonoids correspond to paracyclic
structures on the associated $2$-Segal sets.
The commutative Frobenius case is further developed by Marx
and Mehta~\cite{MM2025TQFT}, who relate it to coherent
span-valued two-dimensional TQFTs and give examples arising
from effect algebras.
A related simplicial generalisation of effect algebras and
effect algebroids is provided by the simplicial effects of
Okay, Castillo, and Stern~\cite{OCS2026}.

In the descriptions of Frobenius structures discussed above,
rotational data are explicitly specified, taking the form
of cyclic or paracyclic structures in the coherent setting.

The $\varepsilon$-simplicial approach requires only additional data
on edges. An $\varepsilon$-simplicial set is a simplicial set $X$
equipped with a set of witnesses and a map from this set to $X_1$.
Edges admitting a witness are called $\varepsilon$-edges or marked
edges; an edge may have several distinct witnesses.
A paracyclic structure, by comparison, specifies bijections
$t_n\colon X_n\to X_n$ in every dimension, subject to compatibility
identities with the face and degeneracy maps.
In the $\varepsilon$-simplicial setting, rotations are constructed
by filling $\varepsilon$-horns:
\begin{equation}\label{eq:rotByEhorn}
    \begin{tikzcd}
        \ELambda^n_k \arrow[r] \arrow[d,hook] & X \\
        \Sim^n \arrow[ur,dashed,"\exists"']
    \end{tikzcd}
    \qquad n\geq 1,\quad k\in\{0,n\}.
\end{equation}
Here $\Sim^n$ is the standard simplex with a single witness on the
edge from its first to its last vertex, and $\ELambda^n_k$ is the
corresponding outer horn with the induced marking.
For $n=1$, the horn consists of a single unmarked vertex.
In dimension two, a filler completes a triangle whose long edge
is marked; its remaining edge describes a rotation of the given edge.
Higher-dimensional fillers allow us to construct rotations of higher
simplices and relate their choices.

The $\varepsilon$-simplicial nerve of every relational Frobenius monoid
admits unique fillers for all $\varepsilon$-horns~\cite{Lachman2026}.
We study these horn-filling conditions without requiring uniqueness.
An $\varepsilon$-simplicial set is called $\varepsilon$-fibrant if
it admits fillers for all $\varepsilon$-horns and satisfies additional
$\varepsilon$-coherence conditions, also expressed as lifting
properties. These coherence conditions ensure that $\varepsilon$-edges
are closed under double rotation.

A concrete motivation for stability comes
from~\cite[Theorem~4.1]{Lachman2025Simplicial}: for any pair of effect
algebras $E$ and $F$, the $\varepsilon$-simplicial set
$\mathrm{Fun}(N(E),N(F))$ is isomorphic to $N(H)$ for some relational
Frobenius monoid $H$. Here $N$ is the nerve functor, and
$\mathrm{Fun}$ denotes the internal hom in the relevant presheaf
category. This result suggests investigating whether analogous
stability results hold for general $\varepsilon$-fibrant objects.

Our first stability theorem assumes $\varepsilon$-fibrancy of both
$X$ and $Y$. Under this assumption, $\mathrm{Fun}(X,Y)$ satisfies
the $\varepsilon$-coherence conditions and admits fillers for all
$\varepsilon$-horns of dimension at least two
(Theorem~\ref{thm:stability}).
The theorem does not settle the one-dimensional
$\varepsilon$-horn conditions for $\mathrm{Fun}(X,Y)$.
Nevertheless, we obtain a naturally defined $\varepsilon$-fibrant
subobject
\[
    \mathrm{Fun}^{\eta}(X,Y)\subseteq\mathrm{Fun}(X,Y).
\]
For an $\varepsilon$-fibrant object, the edges obtained by rotating
$\varepsilon$-edges are called $\eta$-edges
(Definition~\ref{def:etaEdge}).
The subobject $\mathrm{Fun}^{\eta}(X,Y)$ consists of those simplices
whose vertices correspond to simplicial maps preserving these edges;
its $\varepsilon$-fibrancy is proved in Corollary~\ref{cor:lastCor}.
If, in addition, $X$ has unique fillers for $\varepsilon$-horns in
dimensions one and two, then the entire internal hom
$\mathrm{Fun}(X,Y)$ is $\varepsilon$-fibrant
(Corollary~\ref{cor:stabilityUniqueFillers}).

To incorporate the preservation of $\eta$-edges into the ambient
category, we introduce $\eta\varepsilon$-simplicial sets, which carry
explicit $\eta$-witnesses in addition to $\varepsilon$-witnesses.
Extending the coherence conditions to this setting gives a notion
of $\eta\varepsilon$-fibrancy for which full stability under internal
homs holds:
\[
    X,Y\text{ are }\eta\varepsilon\text{-fibrant}
    \quad\Longrightarrow\quad
    \mathrm{Fun}(X,Y)\text{ is }\eta\varepsilon\text{-fibrant}
\]
(Theorem~\ref{thm:mainEta}).
This framework also yields the construction of
$\mathrm{Fun}^{\eta}(X,Y)$ in the $\varepsilon$-simplicial setting.

Allowing non-unique rotations has a concrete motivation in the
test-space approach to quantum logic. Test spaces and their
generalisations describe quantum logics in terms of events and
tests~\cite{FPR,Gudder1997EffectTestSpacesAndAlgebras}.
In this setting, the complement of an event is taken relative to
a test containing it and may depend on that test.
Passing from an algebraic effect test space to its effect algebra
identifies all tests with the top element, while the local complements
of an event represent the orthocomplement of its equivalence class.
The $\varepsilon$-simplicial description retains these distinctions:
tests are represented by marked edges, and local complements are
interpreted as rotations~\cite{Lachman2026}.

A second group of results describes coherent choices of rotations.
The relationship with paracyclic sets is expressed by a functor
\[
    P^*\colon
    \mathrm{Set}_{\Lambda_\infty}
    \longrightarrow\eSSet,
\]
which retains the underlying simplicial set and marks the rotations
of degenerate edges. The $\varepsilon$-horn filling conditions need
not hold for $P^*Z$ when $Z$ is an arbitrary paracyclic set.
We study the chain of adjunctions
\[
    P_!\dashv P^*\dashv P_*
\]
and prove that an $\varepsilon$-simplicial set $X$ admits fillers for
all $\varepsilon$-horns if and only if the counit
\[
    \epsilon'_X\colon P^*P_*X\longrightarrow X
\]
is a trivial Kan fibration (Proposition~\ref{prop:counitIsTriv}(i)).
Consequently, coherent choices of rotations of each simplex form
a contractible Kan complex (Corollary~\ref{cor:contractibleRotations}).
This conclusion requires only the $\varepsilon$-horn filling conditions.
For the stronger notion of $\varepsilon$-fibrancy, we show that $X$
is $\varepsilon$-fibrant if and only if a canonical enlargement of
the counit, $\widehat\epsilon_X\colon\widehat X\to X$, has the right
lifting property with respect to all monomorphisms of
$\varepsilon$-simplicial sets (Proposition~\ref{prop:counitIsTriv}).

The stability proofs use cellular decompositions obtained from these
adjunctions to reduce the problem to pushout-product calculations.
We study the weakly saturated classes generated by the
$\varepsilon$-horn and $\varepsilon$-coherence inclusions.
For example, if $B$ has no marked edges, the unit
$B\to P^*P_!B$ of the free paracyclic extension belongs to the
weakly saturated class generated by the $\varepsilon$-horn inclusions
(Proposition~\ref{prop:trivKanFibHelp}).
A key technical result shows that this class also contains the
pushout-products of $\varepsilon$-horn inclusions with boundary
inclusions of cycle objects (Proposition~\ref{prop:key}).

The article is organised as follows.
Section~2 reviews $\varepsilon$-simplicial sets and the tools for
working with lifting properties and weakly saturated classes.
Section~3 develops rotations of simplices and lifting problems,
and Section~4 introduces $\varepsilon$-coherence and
$\varepsilon$-fibrancy.
Section~5 establishes auxiliary results about cycles.
Section~6 studies the paracyclic adjunctions and the characterisations
of horn filling and fibrancy in terms of their counits.
Section~7 proves the stability theorem in the
$\varepsilon$-simplicial setting.
Finally, Section~8 introduces $\eta\varepsilon$-simplicial sets,
proves their stability under internal homs, and constructs the
$\varepsilon$-fibrant subobject $\mathrm{Fun}^{\eta}(X,Y)$.

\section{Preliminaries}\label{sec:preliminaries}

We assume basic familiarity with simplicial sets. We write $\Delta$
for the simplex category, whose objects are the finite linearly ordered
sets $[n]=\{0,\ldots,n\}$, $n\geq0$, and whose morphisms are
order-preserving maps. We denote the category of simplicial sets by
$\SSet=\mathrm{Set}^{\Delta^{\mathrm{op}}}$. For a simplicial set $X$,
we write $X_n=X([n])$ and denote its face and degeneracy maps by
$d_i$ and $s_i$, respectively.

The standard simplex $\Delta^n$ is the representable simplicial set
$\mathrm{Hom}_{\Delta}(-,[n])$. We identify an element $x\in X_n$
with the corresponding morphism $x\colon\Delta^n\to X$.
We write $\partial_i\Delta^n$ for the face opposite the vertex $i$,
$\partial\Delta^n$ for the boundary, and $\Lambda_k^n$ for the horn
obtained by taking the union of all faces except $\partial_k\Delta^n$.
For a finite linearly ordered set $I$, we use $\Delta^I$ for the
simplex with vertex set $I$, with the convention $\Delta^{\emptyset}
=\emptyset$. We write $y_{\mathcal C}$ for the Yoneda embedding of
a small category $\mathcal C$, omitting the subscript when the
category is clear.

\begin{defn}
Let $\Sim$ be the category obtained from $\Delta$ by freely adjoining
an object $[\varepsilon]$ and a morphism $[1]\to[\varepsilon]$.
An \emph{$\varepsilon$-simplicial set} is a presheaf on $\Sim$.
We denote their category by
\[
    \eSSet=\mathrm{Set}^{\Sim^{\mathrm{op}}}
\]
\end{defn}
Equivalently, an $\varepsilon$-simplicial set is a simplicial set $X$
equipped with a set $X_\varepsilon$ and a structure map
$X_\varepsilon\to X_1$. We use the same notation for an
$\varepsilon$-simplicial set and its underlying simplicial set.

The elements of $X_\varepsilon$ are called
\emph{$\varepsilon$-witnesses}, and the edges in the image of the
structure map are called \emph{$\varepsilon$-edges} or
\emph{marked edges}. The structure map need not be injective,
so an edge may have several distinct witnesses.
Morphisms are simplicial maps together with maps of witness sets
commuting with the structure maps. In particular, morphisms preserve marked edges.

We regard a simplicial set as an $\varepsilon$-simplicial set with
no $\varepsilon$-witnesses. This defines a fully faithful embedding
$\SSet\hookrightarrow\eSSet$. In particular, $\Delta^n$ denotes
the unmarked standard simplex also when considered in $\eSSet$;
under this convention, $\Delta^n=y_{\Sim}([n])$.

For $n\geq0$, let $\Sim^n$ be the $\varepsilon$-simplicial set whose
underlying simplicial set is $\Delta^n$ and whose set of
$\varepsilon$-witnesses consists of a single element lying over
the edge $\langle0,n\rangle$. In particular,
$\Sim^1=y_{\Sim}([\varepsilon])$, whereas $\Sim^0$ has one vertex
and one witness on its degenerate edge. The latter is the terminal
object of $\eSSet$ and is distinct from the unmarked object
$\Delta^0$. For a nonempty finite linearly ordered set $I$, we use
$\Sim^I$ for the corresponding marked simplex, with its unique
witness on the edge from the first to the last vertex.

Subobjects of an $\varepsilon$-simplicial set carry the restricted
structure maps. In particular, $\partial_i\Sim^n$ denotes the face
with the marking induced from $\Sim^n$: it contains the witness
precisely when it contains the edge $\langle0,n\rangle$.
This should be distinguished from a separately marked simplex
$\Sim^I$. In diagrams, a highlighted edge indicates a chosen
$\varepsilon$-witness on that edge.

For $n\geq1$ and $0\leq k\leq n$, let
$\ELambda_k^n\subseteq\Sim^n$ be the subobject whose underlying
simplicial set is $\Lambda_k^n$ and whose marking is induced from
$\Sim^n$. When $k\in\{0,n\}$, we call this an
\emph{$\varepsilon$-horn}. We write
\[
    \eHorn:=
    \bigl\{\ELambda_k^n\hookrightarrow\Sim^n
    \mid n\geq1, k\in\{0,n\}\bigr\}
\]
for the set of $\varepsilon$-horn inclusions.
Note that for $n=1$, each $\varepsilon$-horn is a single unmarked vertex.

We next recall the lifting and saturation constructions used throughout the paper.
Let $f\colon A\to B$, $g\colon X\to Y$ be morphisms of a category $\mathcal{C}$. We say that $f$ has the \emph{(unique) left lifting property} with respect to $g$, or equivalently, that $g$ has the \emph{(unique) right lifting property} with respect to $f$, if for every commutative square
\[
\begin{tikzcd}[ampersand replacement=\&]
A\arrow[r,"h_1"]\arrow[d,"f"]\&X\arrow[d,"g"]\\
B\arrow[r,"h_2"]\&Y
\end{tikzcd}
\]
there exists a (unique) morphism $l\colon B\to X$ such that $lf=h_1$ and $gl=h_2$. We denote this property by $f\boxslash g$ (or $f\boxslash_! g$).
In the case that for every $h_1$ there exists a morphism $l$ such that $lf=h_1$, we denote this property by $f\boxslash X$.

For a class $\mathcal I$ of morphisms, we write
\[
    \mathcal I^{\boxslash}
    :=\{p\mid i\boxslash p \text{ for every }i\in\mathcal I\},
    \qquad
    \prescript{\boxslash}{}{\mathcal I}
    :=\{i\mid i\boxslash p \text{ for every }p\in\mathcal I\}.
\]
Thus $\mathcal I\boxslash p$ means that
$p\in\mathcal I^{\boxslash}$. Similarly,
$\mathcal I\boxslash X$ means that every morphism in $\mathcal I$
has the left lifting property with respect to $X\to1$.

Let $\mathcal{I}$ be a class of morphisms in a category $\mathcal{C}$. We say that $\mathcal{I}$ is \emph{weakly saturated} if it is closed under pushouts, transfinite compositions, and retracts.
By $\overline{\mathcal{I}}$ we denote the closure of $\mathcal{I}$ under these three operations. It is well known (\cite{H2003}, Proposition 10.3.2) that for each morphism $h\in\mathcal{C}$, we have
\[
\mathcal{I}\boxslash h \implies \overline{\mathcal{I}}\boxslash h .
\]
Moreover, if $h$ admits the unique RLP with respect to $\mathcal{I}$, then $h$ admits the unique RLP with respect to $\overline{\mathcal{I}}$.

Similarly, we denote by $\widehat{\mathcal{I}}$ the closure of $\mathcal{I}$ under pushouts, transfinite compositions, and cancellation of a left factor, meaning that
\[
h_1h_2\in\widehat{\mathcal{I}}
\implies
h_2\in\widehat{\mathcal{I}}.
\]
We say that an object $X$ is $\mathcal{I}$-injective if $h\boxslash X$ for every $h\in\mathcal{I}$.

\begin{lemma}[\cite{AHS2007}, Lemma~2.14]\label{lem:fibHelp}
Let $\mathcal{C}$ be a category and let $\mathcal{I}$ be a class of morphisms in $\mathcal{C}$. Then for each object $X\in\mathcal{C}$ we have
\[
\mathcal{I}\boxslash X \implies \widehat{\mathcal{I}}\boxslash X .
\]
\end{lemma}

For a set $\mathcal I$ of morphisms in any of the presheaf
categories considered here, the following characterisations hold;
see~\cite{MH} and~\cite[Corollary~6.7]{AHS2007}:
\begin{equation*}
\overline{\mathcal{I}}
= \prescript{\boxslash}{}{(\mathcal{I}^{\boxslash})}, \quad\text{and}\quad
\widehat{\mathcal{I}}
= \{h \mid h\boxslash X
\text{ for every } \mathcal{I}\text{-injective object } X\}.
\end{equation*}
\begin{lemma}\label{lem:fibHelp2}
    Let $p\colon X\rightarrow Y$ and $f\colon A\rightarrow B$, $g\colon B\rightarrow C$ be morphisms in a category.
    If $gf\boxslash p$ and $g\boxslash Y$, then $f\boxslash p$.
\end{lemma}
\begin{proof}
Assume a lifting problem:
\begin{equation}\label{eq:lifProbHelp}
    \begin{tikzcd}[ampersand replacement=\&]
A\arrow[r,"h_1"]\arrow[d,"f"]\&X\arrow[d,"p"]\\
B\arrow[r,"h_2"]\&Y
    \end{tikzcd}
\end{equation}
By the assumptions, there is $h_3\colon C\rightarrow Y$, such that $h_3g=h_2$. Substituting into $ph_1=h_2f$,
we get $h_3gf=ph_1$. Since $gf\boxslash p$, there is $l$, such that $lgf=h_1$ and $pl=h_3$. It easily follows that
$lg$ is the desired solution of~\eqref{eq:lifProbHelp}.
\end{proof}
Limits and colimits in presheaf categories are computed pointwise.
In particular, this applies to $\SSet$ and $\eSSet$. We will use the familiar notation of intersection and union. The closure under pushout gives the following result.
\begin{lemma}\label{lem:pushOut}
    Let $X,Y\subseteq Z$ be $\varepsilon$-simplicial sets and $\mathcal{I}$ be a weakly saturated class. If $i\colon X\cap Y\hookrightarrow Y$ belongs to $\mathcal{I}$, then $j\colon X\hookrightarrow X\cup Y$ belongs to $\mathcal{I}$ as well.
\end{lemma}
\begin{proof}
    Since colimits in $\eSSet$ are computed pointwise, we see that $j$ is a pushout of $i$. Therefore, $j$ belongs to $\mathcal{I}$ whenever $i$ does.
\end{proof}
\begin{lemma}\label{lem:genOfMono}
    The class of monomorphisms in $\SSet$ is generated as a weakly saturated class by $\{\partial\Delta^n\hookrightarrow \Delta^n\mid n\geq 0 \}$ and the class of monomorphisms in $\eSSet$ is
    generated as a weakly saturated class by $\{\partial\Delta^n\hookrightarrow \Delta^n\mid n\geq 0 \}\cup \{\Delta^1\hookrightarrow \Sim^1\}$.
\end{lemma}
\begin{proof}
     Monomorphisms of presheaves are precisely the pointwise
    injective morphisms. The first statement is standard.
    For the second, identify a monomorphism of $\varepsilon$-simplicial sets
    $A\hookrightarrow B$ with an inclusion and factor it as
    $A\hookrightarrow C\hookrightarrow B$, where $C$ has
    the same underlying simplicial set as $B$ and
    $C_\varepsilon=A_\varepsilon$.
    The first inclusion is obtained by attaching simplices
    along their boundaries, and the second by attaching
    one copy of $\Sim^1$ along $\Delta^1$ for each missing
    $\varepsilon$-witness.
    Conversely, the generators are monomorphisms,
    and monomorphisms form a weakly saturated class.
\end{proof}

     \begin{theorem}[{\cite[Theorem~13]{Lachman2026}}]
        \label{thm:lifting1}
    Let $n\geq 1$ be an integer and $I\subseteq\{0,\ldots,n\}$ be such that $|I\cap \{0,n\}|=1$. Then
    \begin{align}\label{eq:genLift}
        \bigcup_{i\in I}\partial_{i}\Sim^n\hookrightarrow \Sim^n
    \end{align}
    belongs to $\overline{\eHorn}$.
\end{theorem}

To study lifting properties of internal homs, we use the cartesian closed structure of presheaf categories and the associated pushout-product.
We write $\mathrm{Fun}(X,Y)$ for the internal hom in the ambient
category, characterised by the natural bijection
\[
    \mathrm{Hom}(Z,\mathrm{Fun}(X,Y))
    \cong \mathrm{Hom}(Z\times X,Y).
\]
Products are computed pointwise; in particular, for
$\varepsilon$-simplicial sets,
$(X\times Y)_\varepsilon=X_\varepsilon\times Y_\varepsilon$.

By Yoneda, the simplices and witnesses of the internal hom
in $\eSSet$ are given by
\[
    \mathrm{Fun}(X,Y)_n
    \cong \mathrm{Hom}_{\eSSet}(\Delta^n\times X,Y),
    \qquad
    \mathrm{Fun}(X,Y)_\varepsilon
    \cong \mathrm{Hom}_{\eSSet}(\Sim^1\times X,Y).
\]
Since $\Delta^0$ has no $\varepsilon$-witnesses, the vertices
of $\mathrm{Fun}(X,Y)$ correspond to simplicial maps between
the underlying simplicial sets of $X$ and $Y$.
Such maps need not preserve $\varepsilon$-witnesses.

For morphisms $u\colon A\to B$ and $v\colon C\to D$, their
\emph{pushout-product} is the canonical morphism
\[
    u\Box v\colon
    (B\times C)\coprod_{A\times C}(A\times D)
    \longrightarrow B\times D.
\]
When $u$ and $v$ are inclusions, we identify its domain with
the union $(B\times C)\cup(A\times D)$ inside $B\times D$.
The internal-hom adjunction identifies the corresponding
lifting problems, giving
\[
    (u\Box v)\boxslash Y
    \quad\Longleftrightarrow\quad
    u\boxslash
    \bigl(\mathrm{Fun}(D,Y)\to\mathrm{Fun}(C,Y)\bigr).
\]

In Section~\ref{sec:6} we investigate the relationship between $\varepsilon$-simplicial sets and paracyclic sets through a chain of adjunctions.
We will need there the following lemma.
\begin{lemma}\label{lem:eqLifProb}
Let $P_!\dashv P^*\dashv P_*$ be adjunctions between
categories $\mathbf C$ and $\mathbf D$, with
$P^*\colon\mathbf D\to\mathbf C$.
Write $\eta\colon\mathrm{Id}\to P^*P_!$ and
$\epsilon'\colon P^*P_*\to\mathrm{Id}$ for the unit
and counit, respectively.
Let $m\colon A\to B$ and $p\colon X\to Y$ be morphisms
in $\mathbf C$, and assume that the pushout and pullback
below exist.
For corresponding maps under $P^*P_!\dashv P^*P_*$,
the following lifting problems have naturally bijective
sets of fillers:
\begin{equation}\label{diag:eqLifProb}
\begin{tikzcd}[column sep=3.6em, row sep=2.2em]
B\coprod_A P^*P_!(A)
\arrow[r,"{(f_1,f_2)}"]
\arrow[d,"{(\eta_B,P^*P_!(m))}"']
&
X
\arrow[d,"p"]
\\
P^*P_!(B)
\arrow[ru,"l"' description,dotted]
\arrow[r,"g"]
&
Y
\end{tikzcd}
\qquad
\begin{tikzcd}[column sep=3.6em, row sep=2.2em]
A
\arrow[r,"\bar f_2"]
\arrow[d,"m"']
&
P^*P_*(X)
\arrow[d,"{(P^*P_*(p),\epsilon'_X)}"]
\\
B
\arrow[ru,"\bar l"' description,dotted]
\arrow[r,"{(\bar g,f_1)}"]
&
P^*P_*(Y)\times_Y X .
\end{tikzcd}
\end{equation}
\end{lemma}
\begin{proof}
A morphism \(l\colon P^*P_!(B)\to X\) is a solution of the left-hand lifting problem if and only if it satisfies the equations in the left column below. Similarly, a morphism \(\bar l\colon B\to P^*P_*(X)\) is a solution of the right-hand lifting problem if and only if it satisfies the equations in the right column:
\[
\begin{alignedat}{3}
f_1 &= l\circ\eta_B
&\quad&\Longleftrightarrow\quad&
f_1 &= \epsilon'_X\circ\bar l,\\
f_2 &= l\circ P^*P_!(m)
&\quad&\Longleftrightarrow\quad&
\bar f_2 &= \bar l \circ m,\\
g &= p\circ l
&\quad&\Longleftrightarrow\quad&
\bar g &= P^*P_*(p)\circ\bar l .
\end{alignedat}
\]
Therefore \(l\) is a filler of the left-hand lifting problem if and only if its adjoint \(\bar l\) is a filler of the right-hand lifting problem.
\end{proof}

\section{Rotations}\label{sec:3}
In this section we establish rotation relations $\alpha$ and $\beta$
for an $\varepsilon$-simplicial set $X$ admitting fillers
of $\varepsilon$-horns. These relations need not determine
rotations uniquely. In Section~\ref{sec:6}, we show that,
for each simplex $\sigma$, coherent choices of all its iterated
rotations form a contractible Kan complex.

\begin{defn}
Let $X\in\eSSet$. We define binary relations $\alpha$ and $\beta$ on the simplices of $X$ as follows. Given two $n$-simplices $\sigma_1,\sigma_2\in X_n$, $n\geq 0$, we say that $\sigma_2$ is a rotation of $\sigma_1$, and we write
\begin{equation}\label{eq}
\beta\colon \sigma_1\relto \sigma_2
\quad\text{and}\quad
\alpha\colon\sigma_2\relto\sigma_1,
\end{equation}
if there exists a simplex $\hat{\sigma}\colon \Sim^{n+1}\to X$ such that $d_{n+1}\hat{\sigma}=\sigma_1$ and $d_0\hat{\sigma}=\sigma_2$.
\end{defn}
We will sometimes abbreviate~\eqref{eq} by writing $\sigma_2=\beta(\sigma_1)$ or $\sigma_1=\alpha(\sigma_2)$.
A useful property of $\eSSet$ is that lifting problems of the form $S\hookrightarrow\Delta^n$ can be rotated. We identify $S$ with a subset of $\partial_{n+1}\Sim^{n+1}$, and its rotation $\beta(S)$ with a subset of $\partial_0\Sim^{n+1}$, as follows.
First define $\rot{S}\subseteq\Sim^{n+1}$ by
\begin{equation}\label{eq:rotHat}
    \rot{S}=S\cup\bigcup\bigg\{\Sim^{I\cup\{0,n+1\}}\subseteq \Sim^{n+1}\mid I\subseteq \{1,\ldots,n\}\text{\ and\ }\Delta^{\{0\}\cup I}\subseteq S\bigg\}.
\end{equation}
Then we set
\begin{equation}\label{eq:rotBeta}
    \beta(S)=\rot{S}\cap\partial_0\Sim^{n+1}.
\end{equation}
In particular, for a vertex $\{i\}\subseteq \Delta^{n}$, we have $\beta(\{i\})=\{i\}$, if $1\leq i\leq n$, and $\beta(\{0\})=\{n+1\}$.
\begin{example}
    For $S=\partial\Delta^2\subseteq \Delta^2$, we have $\rot{S}=\partial_1\Sim^3\cup\partial_2\Sim^3\cup\Delta^{\{1,2\}}$ and $\beta(S)=\partial(\partial_0\Sim^3)$.
    See the following picture of the inclusions $S\hookrightarrow \rot{S}\hookleftarrow\beta(S)$:
\begin{equation*}
     \begin{tikzcd}[row sep={2.9em,between origins}, column sep={4em,between origins}, ampersand replacement=\&]
    \&{\ \ }
    \&\\[10pt]\&
    2
    \&\\
    0
    \arrow[rr, "{a}"{description}, bend right=5]
    \arrow[ru, "{c}"{description}]
    \&\&
    1\arrow[lu,"b"'description]
\end{tikzcd} \hookrightarrow
     \begin{tikzcd}[row sep={2.9em,between origins}, column sep={4em,between origins}, ampersand replacement=\&,execute at end picture={
\foreach \Nombre in  {A,B,C}
  {\coordinate (\Nombre) at (\Nombre.center);}
\fill[pattern=north east lines,opacity=0.3]
  (A) -- (B) -- (C) -- cycle;
  \foreach \Nombre in  {A,C,D}
  {\coordinate (\Nombre) at (\Nombre.center);}
\fill[pattern=vertical lines,opacity=0.3]
  (A) -- (C) -- (D) -- cycle;
}]
    \&
    |[alias=A]|3
    \&\\[10pt]\&
    |[alias=B]|2
    \arrow[u, "\beta(c)"{description}]
    \&\\
    |[alias=C]|0
    \arrow[ruu, LA=1.5pt, bend left=5]
    \arrow[rr, "{a}"{description}, bend right=5]
    \arrow[ru, "{c}"{description}]
    \&\&
    |[alias=D]|1\arrow[uul, "\beta(a)"{description}, bend right=5]\arrow[lu,"b"'description]
\end{tikzcd}
\hookleftarrow
\begin{tikzcd}[row sep={2.9em,between origins}, column sep={4em,between origins}, ampersand replacement=\&]
    {3}
    \&\\[10pt]
    2\arrow[u,"{\beta(c)}"'{description}]
    \&\\
    {\ }
    \&
    1\arrow[luu,"{\beta(a)}"'{description},bend right=5]\arrow[lu,"b"'description]
\end{tikzcd}
\end{equation*}
Now let $X$ be any $\varepsilon$-simplicial set which admits fillers of $\varepsilon$-horns. Any $\gamma\colon S\rightarrow X$ can be extended to $\hat{\gamma}\colon \rot{S}\rightarrow X$. One can see that $\gamma$ has a filler iff $\beta(\gamma):=\hat{\gamma}\restriction_{\beta(S)}$ has a filler.
\end{example}
\begin{example}\label{eq:rotOfHorns}
For $n\geq 1$ and $0\leq k\leq n$, let
$j\colon\Lambda_k^n\hookrightarrow\Delta^n$.
Then
\(
    \beta(j)\colon
    \beta(\Lambda_k^n)
    \hookrightarrow\Delta^{\{1,\ldots,n+1\}}
\)
is isomorphic to
$\Lambda_{k-1}^n\hookrightarrow\Delta^n$,
where the index is taken modulo $n+1$.\end{example}
\begin{proposition}\label{prop:cohRotation}
    Let $S,\beta(S),\rot{S}\subseteq\Sim^{n+1}$ have the meanings given above. Then the following inclusions belong to $\overline{\eHorn}$:
    \begin{enumerate}[(i)]
        \item $S\hookrightarrow\rot{S}\quad\text{and}\quad\beta(S)\hookrightarrow\rot{S}$;
        \item $\rot{S}\cup \partial_{n+1}\Sim^{n+1}\hookrightarrow \Sim^{n+1}\quad\text{and}\quad\rot{S}\cup \partial_{0}\Sim^{n+1}\hookrightarrow \Sim^{n+1}$.
    \end{enumerate}
    \end{proposition}

\begin{proof} In both cases (i-ii) we will consider only the first inclusions; the second ones are analogous.

    (i) First note that if $S$ does not contain the vertex $\Delta^{\{0\}}$, then $S=\rot{S}=\beta(S)$ and there is nothing to prove. Thus assume that $\Delta^{\{0\}}\subseteq S$, and consider the following chain of inclusions:
\begin{equation}
    S_{-1}:=S\subseteq S_0\subseteq S_1\subseteq\cdots\subseteq S_n=\rot{S},
\end{equation}
where $S_i$ contains all simplices of $S$ together with those simplices $\sigma$ of $\rot{S}$ for which the dimension of $\sigma\cap \partial_0\Sim^{n+1}$
is at most $i$. We prove by induction that all inclusions $S_i\subseteq
S_{i+1}$, $i<n$, belong to $\overline{\eHorn}$. First, $S_0=S\cup\Sim^{\{0,n+1\}}$ and $S\cap \Sim^{\{0,n+1\}}=\Delta^{\{0\}}$. Hence, we can apply
Lemma~\ref{lem:pushOut}.

Now consider the inclusion $S_i\subseteq S_{i+1}$ with $0\leq i<n$. First observe that
\begin{equation*}
    S_{i+1}=S_{i}\cup \bigcup\bigg\{\sigma=\Sim^{\{0\}\cup I\cup \{n+1\}}\mid I\subseteq\{1,\ldots,n\}, |I|=i+1, \Delta^{\{0\}\cup I}\subseteq S\bigg\}.
\end{equation*}

Let $\sigma$ be an $(i+2)$-dimensional cell of $S_{i+1}$ of the form $\sigma=\Sim^{\{0\}\cup I\cup \{n+1\}}$, where  $\sigma\cap \partial_{n+1}\Sim^{n+1}=\Delta^{\{0\}\cup I}\subseteq S\subseteq S_i$.
Moreover, for each $k\in I$, we have $\Sim^{\{0\}\cup I\setminus\{k\}\cup\{n+1\}}\subseteq S_i$. Hence all but one of the faces of $\sigma$ are contained in $S_i$ and $\sigma\cap S_i$ is isomorphic to $\ELambda^{i+2}_0$. Therefore $S_i\hookrightarrow S_i\cup \sigma$ belongs to $\overline{\eHorn}$ by Lemma~\ref{lem:pushOut}.
Since all simplices $\sigma$ under consideration intersect in $S_i$, we may attach them one by one.

We will handle (ii) using the same strategy as in part~(i).
First suppose that $\Delta^{\{0\}}$ does not belong to $S$.
Then $\rot{S}=S\subseteq \partial_{n+1}\Sim^{n+1}$, and the inclusion $\partial_{n+1}\Sim^{n+1}\hookrightarrow\Sim^{n+1}$ belongs to $\overline{\eHorn}$ by Theorem~\ref{thm:lifting1}.

Now suppose that $\Delta^{\{0\}}\subseteq S$.
For $0\leq i\leq n$, define
\[
    T_i
    =
    \rot{S}\cup\partial_{n+1}\Sim^{n+1}
    \cup
    \bigcup
    \big\{
    \Sim^{\{0\}\cup I\cup\{n+1\}}\mid
        I\subseteq\{1,\ldots,n\},
        |I|\leq i
    \big\}
\]
Since $\Sim^{\{0,n+1\}}\subseteq\rot{S}$, these
objects form a chain
\[
    \rot{S}\cup\partial_{n+1}\Sim^{n+1}
    =T_0\subseteq\cdots\subseteq T_n=\Sim^{n+1}.
\]
Moreover,
\[
    T_{i+1}
    =
    T_i\cup \bigcup
    \big\{
    \Sim^{\{0\}\cup I\cup\{n+1\}}\mid
        I\subseteq\{1,\ldots,n\},
        |I|\leq i+1
    \big\}.
\]
We show that each inclusion $T_i\hookrightarrow T_{i+1}$, $i<n$, belongs to $\overline{\eHorn}$. Let $\tau\in T_{i+1}\setminus T_i$ be of the form
   $\Sim^{\{0\}\cup I\cup\{n+1\}}$.
As in the previous part, we observe that $\tau\cap T_i$ is isomorphic to an $\varepsilon$-horn. Hence $T_i\subseteq T_{i+1}$ is a composition of pushouts of $\varepsilon$-horn inclusions.
\end{proof}
Let us describe the rotation of a lifting problem. Assume a lifting problem of the form
\begin{equation}
    \begin{tikzcd}
S\arrow[r,"f"]\arrow[d,"i",hook] & X\arrow[d,"p"]\\
{\Delta^n}\arrow[r,"g"] & Y
\end{tikzcd}
\end{equation}
where $p$ has the RLP with respect to $\varepsilon$-horns and $Y$ admits fillers of $\varepsilon$-horns.
By the assumptions, $(\Delta^n\hookrightarrow\Sim^{n+1})\boxslash Y$ and $(S\hookrightarrow\rot{S})\boxslash p$. Hence we can define $\rot{g}, \rot{f}$ and their restrictions $\beta(g), \beta(f)$ as in the following diagram:
\begin{equation}\label{diag:rotOfProb}
\begin{tikzcd}[row sep=1.3em, column sep=1.8em]
S
  \arrow[rr, "f"]
  \arrow[dd, "i"']
  \arrow[dr, hook]
&&
X
  \arrow[dd, "p"near start]
\\
&
{\rot{S}}
  \arrow[ur, "\rot{f}" description]
  \arrow[dd]
&&
{\beta(S)}
  \arrow[lu, "\beta(f)"']
  \arrow[dd, hook,"{\beta(i)}"]
  \arrow[ll,hook]
\\
\Delta^n
  \arrow[rr, "g" near start]
  \arrow[dr,hook]
&&
Y
\\
&
\Sim^{n+1}
  \arrow[ur, "\rot{g}" description]
&&
{\Delta^{\{1,\ldots,n+1\}}}
  \arrow[lu, "{\beta(g)}"']
  \arrow[ll,hook]
\end{tikzcd}
\end{equation}
We identify $\partial_0\Sim^{n+1}
=\Delta^{\{1,\ldots,n+1\}}$ with $\Delta^n$
via the order-preserving bijection.
\begin{proposition}\label{prop:cohRot}
    Let $X,Y$ be $\varepsilon$-simplicial sets admitting fillers of $\varepsilon$-horns and $p\colon X\rightarrow Y$ be a morphism having the RLP with respect to $\varepsilon$-horns.
    Assume the situation of diagram~\eqref{diag:rotOfProb} described above. Then the following left lifting problem has a solution if and only if the right lifting problem has a solution:
\begin{equation}\label{diag:twoLifProb}
     \begin{tikzcd}
S\arrow[r,"f"]\arrow[d,"i",hook] & X\arrow[d,"p"]\\
{\Delta^n}\arrow[r,"g"] & Y
\end{tikzcd}
\qquad\qquad
    \begin{tikzcd}
{\beta(S)}\arrow[r,"{\beta(f)}"]\arrow[d,"{\beta(i)}",hook] & X\arrow[d,"p"]\\
{\Delta^n}\arrow[r,"{\beta(g)}"] & Y
\end{tikzcd}
\end{equation}
Moreover, assume that $p$ has the unique RLP with respect to $\varepsilon$-horns. Then there is a bijection between the fillers of the two lifting problems in~\eqref{diag:twoLifProb}.
\end{proposition}
\begin{proof}
Fix the chosen extensions $\rot{g}$ and $\rot{f}$.
Let $L$ and $L'$ denote the sets of fillers of the
left-hand and right-hand problems, respectively, and put
\[
    \widehat{L}
    =
    \bigl\{
        h\colon\Sim^{n+1}\to X
        \mid h|_{\rot{S}}=\rot{f},\;
        ph=\rot{g}
    \bigr\}.
\]
Restriction to the faces $\partial_{n+1}\Sim^{n+1}$ and $\partial_{0}\Sim^{n+1}$, respectively, gives the span
\begin{equation}\label{eq:span}
    L
    \leftarrow
    \widehat{L}
    \rightarrow
    L'.
\end{equation}
Since
\[
    \rot{S}\cap\partial_{n+1}\Sim^{n+1}=S,
    \qquad
    \rot{S}\cap\partial_0\Sim^{n+1}=\beta(S),
\]
Proposition~\ref{prop:cohRotation}(ii) and the RLP
of $p$ imply that both restriction maps in~\eqref{eq:span} are surjective.
Thus $L$ is nonempty if and only if $L'$ is nonempty.

If $p$ has the unique RLP with respect to
$\varepsilon$-horns, it also has the unique RLP
with respect to $\overline{\eHorn}$.
Both restriction maps are then bijective,
yielding the desired bijection $L\cong L'$.
\end{proof}

\begin{example}
    Let $N(F)$ be an $\varepsilon$-nerve of a Frobenius monoid in $\Rel$. The associativity of $F$ corresponds to
    \[
    (\partial_0\Delta^3\cup\partial_2\Delta^3 \hookrightarrow \Delta^3)\boxslash N(F) \quad\text{and}\quad (\partial_1\Delta^3\cup\partial_3\Delta^3 \hookrightarrow \Delta^3)\boxslash N(F).
    \]
    See that the second lifting condition is a rotation of the first one. Hence, under the assumption of fillers of $\varepsilon$-horns, either one-sided associativity condition is equivalent to full associativity.
\end{example}
\begin{example}
    Let $X$ be an $\varepsilon$-simplicial set admitting fillers of $\varepsilon$-horns.
    It follows from Proposition~\ref{prop:cohRot} (set $Y=\Sim^0$) that
\[
\bigl(\exists k\in\{0,1,2,3\}:\;
(\Lambda^3_k\hookrightarrow\Delta^3)\boxslash X\bigr)
\implies
\bigl(\forall k\in\{0,1,2,3\}:\;
(\Lambda^3_k\hookrightarrow\Delta^3)\boxslash X\bigr).
\]
    Indeed, for $i\colon\Lambda^3_k\hookrightarrow\Delta^3$ the rotated embedding $\beta(i)\colon\beta(\Lambda^3_k)\hookrightarrow\Delta^3$ is isomorphic to $\Lambda^3_{k-1}\hookrightarrow\Delta^3$, where $k-1$ is taken modulo $4$.

    Recall Theorem~6 in \cite{Lachman2026}. It is shown there that for a monoid $\mathcal{M}$ in $\mathbf{Rel}$ and its nerve $N(\mathcal{M})\in\eSSet$, the condition $(\Lambda^3_k\hookrightarrow\Delta^3)\boxslash N(\mathcal{M})$, for $k\in\{0,1,2,3\}$,
     is equivalent to the algebraic properties \emph{left cancellation}, \emph{partiality}, \emph{partiality} (again), and \emph{right cancellation}, respectively. Our observation reflects the fact that
     in the case of Frobenius monoids in $\mathbf{Rel}$, the mentioned algebraic properties are equivalent.
    \end{example}
    \begin{example}
    Let $F$ be a Frobenius monoid in $\Rel$ which is a category. Then $F$ is a groupoid. Let $N(F)$ be the $\varepsilon$-simplicial nerve of $F$ constructed in~\cite{Lachman2026}.
    Then $N(F)$ admits (unique) fillers of $\varepsilon$-horns. Since $F$ is a category, $N(F)$ admits unique fillers of all inner horns $\Lambda^n_i\hookrightarrow \Delta^n$, $1\leq i \leq n-1$.
    By rotation also the outer horns $\Lambda_0^n$ and $\Lambda_n^n$, $n\geq 2$, admit unique fillers in $N(F)$.
    Consequently $F$ is a groupoid.
\end{example}
\begin{example}
    Assume an $\varepsilon$-simplicial set $X$ admits fillers of $\varepsilon$-horns. Then $X$ is a quasicategory iff it is a Kan complex. Indeed, rotation gives fillers for all outer horns
in dimensions $n\geq2$, while one-dimensional horns
have degenerate fillers. The converse is immediate.
\end{example}
\section{Epsilon-coherence conditions}
The $\varepsilon$-simplicial sets of interest have the following property:
Given a rhombus of some particular shape with three outer $\varepsilon$-edges, the last outer edge is an $\varepsilon$-edge as well.
\begin{defn}
    Let $X$ be an $\varepsilon$-simplicial set. We say that $X$ satisfies the \emph{$\varepsilon$-coherence conditions} if it has the right lifting property with respect to the following morphisms:

\begin{equation}\label{L1}
\tag{L1}
\begin{aligned}
    \begin{tikzcd}[row sep={1.8em,between origins}, column sep={2.4em,between origins}, ampersand replacement=\&,execute at end picture={
\foreach \Nombre in  {A,B,C}
  {\coordinate (\Nombre) at (\Nombre.center);}
\fill[pattern=vertical lines,opacity=0.3]
  (A) -- (B) -- (C) -- cycle;
  \foreach \Nombre in  {A,C,D}
  {\coordinate (\Nombre) at (\Nombre.center);}
\fill[pattern=vertical lines,opacity=0.3]
  (B) -- (C) -- (D) -- cycle;
}]
    \&
    |[alias=C]|2\arrow[rr,"e_3"{description}] \&\& |[alias=D]|3
    \\[10pt]\&
    \&\&\\
    |[alias=A]|0
    \arrow[ruu,LA=1.5pt, "e_1"{description}]
    \arrow[rr,LA=1.5pt,"e_4"{description}]
    \&\&
    |[alias=B]|1\arrow[uul, "a"{description}]\arrow[ruu,LA=1.5pt,"e_2"{description}]\&
\end{tikzcd}
\quad &\hookrightarrow\quad
\begin{tikzcd}[row sep={1.8em,between origins}, column sep={2.4em,between origins}, ampersand replacement=\&,execute at end picture={
\foreach \Nombre in  {A,B,C}
  {\coordinate (\Nombre) at (\Nombre.center);}
\fill[pattern=vertical lines,opacity=0.3]
  (A) -- (B) -- (C) -- cycle;
  \foreach \Nombre in  {A,C,D}
  {\coordinate (\Nombre) at (\Nombre.center);}
\fill[pattern=vertical lines,opacity=0.3]
  (B) -- (C) -- (D) -- cycle;
}]
    \&
    |[alias=C]|2\arrow[rr,LA=1.5pt,"e_3"{description}] \&\& |[alias=D]|3
    \\[10pt]\&
    \&\&\\
    |[alias=A]|0
    \arrow[ruu,LA=1.5pt, "e_1"{description}]
    \arrow[rr, LA=1.5pt,"e_4"{description}]
    \&\&
    |[alias=B]|1\arrow[uul, "a"{description}]\arrow[ruu,LA=1.5pt,"e_2"{description}]\&
\end{tikzcd}
\end{aligned}
\end{equation}

\begin{equation}\label{L2}
\tag{L2}
\begin{aligned}
    \begin{tikzcd}[row sep={1.8em,between origins}, column sep={2.4em,between origins}, ampersand replacement=\&,execute at end picture={
\foreach \Nombre in  {A,B,C}
  {\coordinate (\Nombre) at (\Nombre.center);}
\fill[pattern=vertical lines,opacity=0.3]
  (A) -- (B) -- (C) -- cycle;
  \foreach \Nombre in  {A,C,D}
  {\coordinate (\Nombre) at (\Nombre.center);}
\fill[pattern=vertical lines,opacity=0.3]
  (B) -- (C) -- (D) -- cycle;
}]
    \&
    |[alias=C]|2\arrow[rr,LA=1.5pt,"e_3"{description}] \&\& |[alias=D]|3
    \\[10pt]\&
    \&\&\\
    |[alias=A]|0
    \arrow[ruu,LA=1.5pt, "e_1"{description}]
    \arrow[rr,"e_4"{description}]
    \&\&
    |[alias=B]|1\arrow[uul, "a"{description}]\arrow[ruu,LA=1.5pt,"e_2"{description}]\&
\end{tikzcd}
\quad &\hookrightarrow\quad
\begin{tikzcd}[row sep={1.8em,between origins}, column sep={2.4em,between origins}, ampersand replacement=\&,execute at end picture={
\foreach \Nombre in  {A,B,C}
  {\coordinate (\Nombre) at (\Nombre.center);}
\fill[pattern=vertical lines,opacity=0.3]
  (A) -- (B) -- (C) -- cycle;
  \foreach \Nombre in  {A,C,D}
  {\coordinate (\Nombre) at (\Nombre.center);}
\fill[pattern=vertical lines,opacity=0.3]
  (B) -- (C) -- (D) -- cycle;
}]
    \&
    |[alias=C]|2\arrow[rr,LA=1.5pt,"e_3"{description}] \&\& |[alias=D]|3
    \\[10pt]\&
    \&\&\\
    |[alias=A]|0
    \arrow[ruu,LA=1.5pt, "e_1"{description}]
    \arrow[rr, LA=1.5pt,"e_4"{description}]
    \&\&
    |[alias=B]|1\arrow[uul, "a"{description}]\arrow[ruu,LA=1.5pt,"e_2"{description}]\&
\end{tikzcd}
\end{aligned}
\end{equation}
If $X$ satisfies the $\varepsilon$-coherence conditions and moreover admits fillers of $\varepsilon$-horns, then we say $X$ is \emph{$\varepsilon$-fibrant}.
\end{defn}
The $\varepsilon$-coherence conditions assure that $\varepsilon$-edges of an $\varepsilon$-fibrant $X$ are closed under double-rotation.
Note that the two kinds of $\varepsilon$-horns $\ELambda^n_n\rightarrow \Sim^n$ and $\ELambda^n_0\rightarrow \Sim^n$
as well as (L1) and (L2) are mutually connected by reversing arrows. Hence, the following duality principle holds (the operation $X\mapsto X^{op}$ reverses direction of arrows and leaves $X_\varepsilon$ unchanged):
\begin{lemma}\label{lem:dual}
    An embedding $A\hookrightarrow B$ belongs to the weakly saturated class generated by $\varepsilon$-horns and (L1--2), iff so does $A^{op}\hookrightarrow B^{op}$.
\end{lemma}
It turns out that several rhombi in an $\varepsilon$-fibrant simplicial set satisfy the $4$-for-$3$ conditions:
\begin{theorem}\label{thm:4for3}
    Let $E$ be any of the three $\varepsilon$-simplicial sets in~\eqref{eq:4for3} and $A\subset E$ arise by omitting exactly one element of $E_\varepsilon$:
    \begin{equation}\label{eq:4for3}
           \begin{tikzcd}[row sep={1.8em,between origins}, column sep={2.4em,between origins}, ampersand replacement=\&,execute at end picture={
\foreach \Nombre in  {A,B,C}
  {\coordinate (\Nombre) at (\Nombre.center);}
\fill[pattern=vertical lines,opacity=0.3]
  (A) -- (B) -- (C) -- cycle;
  \foreach \Nombre in  {B,C,D}
  {\coordinate (\Nombre) at (\Nombre.center);}
\fill[pattern=vertical lines,opacity=0.3]
  (B) -- (C) -- (D) -- cycle;
  \foreach \Nombre in  {AA,BB,CC}
  {\coordinate (\Nombre) at (\Nombre.center);}
\fill[pattern=vertical lines,opacity=0.3]
  (AA) -- (BB) -- (CC) -- cycle;
  \foreach \Nombre in  {BB,CC,DD}
  {\coordinate (\Nombre) at (\Nombre.center);}
\fill[pattern=vertical lines,opacity=0.3]
  (BB) -- (CC) -- (DD) -- cycle;
    \foreach \Nombre in  {AAA,BBB,CCC}
  {\coordinate (\Nombre) at (\Nombre.center);}
\fill[pattern=vertical lines,opacity=0.3]
  (AAA) -- (BBB) -- (CCC) -- cycle;
  \foreach \Nombre in  {BBB,CCC,DDD}
  {\coordinate (\Nombre) at (\Nombre.center);}
\fill[pattern=vertical lines,opacity=0.3]
  (BBB) -- (CCC) -- (DDD) -- cycle;
}]
    \&
    |[alias=C]|2\arrow[rr,LA=1.5pt,"e_3"{description}] \&\& |[alias=D]|3 \&
    \&%
    |[alias=CC]|3 \&\& |[alias=DD]|1\arrow[ll,LA=1.5pt,"f_4"{description}]\arrow[ldd,LA=1.5pt,"f_3"{description}]\&
    \&%
    |[alias=CCC]|1\arrow[rr,LA=1.5pt,"f_3"{description}]   \arrow[ldd,LA=1.5pt, "f_4"{description}]\&\& |[alias=DDD]|2
    \\
    [10pt]\&\&\&\&\&\&\&\&\&\&\&
    \\
    |[alias=A]|0
    \arrow[ruu,LA=1.5pt, "e_1"{description}]
    \arrow[rr, LA=1.5pt,"e_4"{description}]
    \&\&
    |[alias=B]|1\arrow[uul, "a"{description}]\arrow[ruu,LA=1.5pt,"e_2"{description}]\&
    \&%
    |[alias=AA]|0
    \arrow[ruu,LA=1.5pt, "f_2"{description}]
    \arrow[rr, LA=1.5pt,"f_1"{description}]
    \&\&
    |[alias=BB]|2\arrow[uul, "b"{description}]\&
    \&%
     |[alias=AAA]|3
    \&\&
    |[alias=BBB]|0\arrow[uul, "c"{description}]\arrow[ruu,LA=1.5pt,"f_1"{description}]\arrow[ll, LA=1.5pt,"f_2"{description}]\&
\end{tikzcd}
    \end{equation}
    Then each $\varepsilon$-fibrant $X\in\eSSet$ satisfies $(A\hookrightarrow E) \boxslash X$.
\end{theorem}
\begin{proof}
    Denote the three $\varepsilon$-simplicial sets in~\eqref{eq:4for3} as $E_1,E_2,E_3$, respectively.
    We first transform the case of $E_3$ to $E_2$.
    By symmetry, we can assume that the edge $f_2$ of $E_3$ is an $\varepsilon$-edge in $A\subset E_3$.
    Since $X$ admits fillers of $\varepsilon$-horns, we can extend $A$ to the whole $3$-simplex~\eqref{diag:new3sim}. Consider the two newly obtained $2$-faces. These form precisely the shape of $E_2$.
    If the statement is true for $E_2$, Lemma~\ref{lem:fibHelp} applies. (Alternatively, we can apply Lemma~\ref{lem:dual}.)
\begin{equation}\label{diag:new3sim}
    \begin{tikzcd}[row sep={1.8em,between origins}, column sep={2.4em,between origins}, ampersand replacement=\&,execute at end picture={
\foreach \Nombre in  {A,B,C,D}
  {\coordinate (\Nombre) at (\Nombre.center);}
\fill[,opacity=0.1]
  (C) -- (D) -- (B) -- (A) -- cycle;
}]
     \&|[alias=C]|1\arrow[from=ddr]\arrow[rr,LA=1.5pt,"f_3"{description}]   \arrow[ldd,LA=1.5pt, "f_4"{description}]\&\& |[alias=D]|2\arrow[llldd, crossing over]
    \\
    [10pt]\&\&\&
    \\
     |[alias=A]|3
    \&\&
    |[alias=B]|0\arrow[ruu,LA=1.5pt,"f_1"{description}]\arrow[ll, LA=1.5pt,"f_2"{description}]\&
\end{tikzcd}
\end{equation}
    Next consider $A\subset E_2$. Because of the symmetry we can assume the missing $\varepsilon$-edge is $f_1$ or $f_2$.
    It is easy to show that for each $A\rightarrow X$, there is in $X$ an extra vertex $\{4\}\rightarrow X$ and a simplex $\sigma\colon \Delta^{\{2,3,4\}}\rightarrow X$ such that  $d_0(\sigma),d_1(\sigma)$ are $\varepsilon$-edges in $X$. Such a $\sigma$, together with the image of $\Delta^{\{0,2,3\}}$, forms the shape of $E_1$. Hence, we have arrived at the case of $E_1$.
This is the last case to be resolved.

Finally assume $A\subset E_1$ and $A$ misses the witness of the $\varepsilon$-edge $e_1$ (the case with $e_2$ is analogous).
Consider the embedding $i\colon A \hookrightarrow B$ pictured in~\eqref{eq:thirdCase}, where $B$ is a union of two $3$-simplices. We will prove $(A\hookrightarrow B)\boxslash X$.
Starting with $A$, we firstly construct $\Delta^{\{0,1,2,x\}}$ with the $\varepsilon$-edge $\Sim^{\{0,x\}}$. Then we add the edge $\Delta^{\{x,3\}}$ and we make it an $\varepsilon$-edge using~\eqref{L1}.
Finally we can fill the $\varepsilon$-horn to obtain a $3$-simplex $\Delta^{\{1,2,x,3\}}$ with the depicted $\varepsilon$-edges. This gives $B$.
Now focus on the faces $\Delta^{\{0,2,x\}}\cup\Delta^{\{2,x,3\}}$. Due to (L2), we can make $e_1$ an $\varepsilon$-edge. The resulting $B'$ satisfies $A\subseteq E_1\subseteq B'$. Hence Lemma~\ref{lem:fibHelp} applies.
\begin{equation}\label{eq:thirdCase}
     \begin{tikzcd}[row sep={2.5em,between origins}, column sep={2.4em,between origins}, ampersand replacement=\&,execute at end picture={
\foreach \Nombre in  {A,B,C}
  {\coordinate (\Nombre) at (\Nombre.center);}
\fill[opacity=0.1]
  (A) -- (B) -- (C) -- cycle;
  \foreach \Nombre in  {A,C,D}
  {\coordinate (\Nombre) at (\Nombre.center);}
\fill[opacity=0.1]
  (B) -- (C) -- (D) -- cycle;
}]
    \&
    |[alias=C]|2\arrow[rr,LA=1.5pt,"e_3"{description}]\arrow[rd] \&\& |[alias=D]|3
    \\ \&
    \& x\arrow[ru,LA=1.5pt] \&\\
    |[alias=A]|0\arrow[rru,LA=1.5pt]
    \arrow[ruu, "e_1"{description}]
    \arrow[rr,LA=1.5pt,"e_4"{description}]
    \&\&
    |[alias=B]|1\arrow[uul]\arrow[ruu,LA=1.5pt,"e_2"']\arrow[u]\&
\end{tikzcd}
\end{equation}
\end{proof}
\begin{example}\label{ex:orthomodular}
    An effect test space is a pair $\mathcal{T}=(X,\mathcal{T})$, where $X$ is a set, and $\mathcal{T}$ a collection of multisubsets of $X$ (represented as finitely supported mappings $X\rightarrow \mathbb{N}_0$), the elements of which are called tests.
    The tests are required to cover $X$ and be irredundant (each pair of tests $T,S$ satisfies $T\not\subset S$).
    The multi-subsets of tests are called events $\mathcal{E}(\mathcal{T})=\{A\leq T\mid T\in\mathcal{T}\}$. Given two events $A$ and $B$, we say $B$ is a local orthocomplement to $A$ if $A + B$ is a test. Further, we define a binary relation $\sim$ on events as $A\sim B$ iff $A$ and $B$ admit a common local orthocomplement $C$.
    The effect test space is called algebraic if, for every
quadruple of events $A,B,C,D$, we have \begin{align}\label{eq:algebraicity}
        A +  B\in\mathcal{T}\ \&\ B + C\in\mathcal{T}\ \&\ C + D\in\mathcal{T}\implies A + D\in\mathcal{T}.
    \end{align}
In this case, $\sim$ is an equivalence relation, and
$\mathcal{E}(\mathcal{T})/{\sim}$ carries a natural
effect algebra structure.
See~\cite{Gudder1997EffectTestSpacesAndAlgebras}
for more details.

Note that we can organize an effect test space $\mathcal{T}$ as an $\varepsilon$-simplicial set: For $n\geq 0$, $n$-simplices correspond to $n$-tuples of  events $A_1,\ldots,A_n$ such that $A_1+\cdots+A_n\leq T$ for some test $T$.
The face maps and degeneracy maps are defined in the obvious way. For $\varepsilon$-edges, we take all tests.
The algebraic property~\eqref{eq:algebraicity} then corresponds to the right lifting property with respect to the following morphism:
    \begin{equation}\label{eq:alg}
      \begin{tikzcd}[row sep={3em,between origins}, column sep={3em,between origins}, ampersand replacement=\&,execute at end picture={
\foreach \Nombre in  {x,y,w,o,z}
  {\coordinate (\Nombre) at (\Nombre.center);}
\fill[opacity=0.3,pattern = vertical lines]
      (x) -- (y)-- (z) -- (w)--(o)--cycle;}]
    |[alias=y]|3 \&  \&|[alias=z]|0\arrow[ll,LA=1.5pt]\arrow[dd,LA=1.5pt]\arrow[ld,"C"{description}]\\
    \&|[alias=o]|2\arrow[lu,"B"{description}]\arrow[rd,"D"{description}]\&\\
    |[alias=x]|1\arrow[ru,"A"{description}]\arrow[uu,LA=1.5pt]\&\&|[alias=w]|4
\end{tikzcd} \quad\hookrightarrow\quad
      \begin{tikzcd}[row sep={3em,between origins}, column sep={3em,between origins}, ampersand replacement=\&,execute at end picture={
\foreach \Nombre in  {x,y,w,o,z}
  {\coordinate (\Nombre) at (\Nombre.center);}
\fill[opacity=0.3,pattern = vertical lines]
      (x) -- (y)-- (z) -- (w)--cycle;}]
    |[alias=y]|3 \&  \&|[alias=z]|0\arrow[ll,LA=1.5pt]\arrow[dd,LA=1.5pt]\arrow[ld,"C"{description}]\\
    \&|[alias=o]|2\arrow[lu,"B"{description}]\arrow[rd,"D"{description}]\&\\
    |[alias=x]|1\arrow[rr,LA=1.5pt]\arrow[ru,"A"{description}]\arrow[uu,LA=1.5pt]\&\&|[alias=w]|4
\end{tikzcd}
\end{equation}
\end{example}
\begin{lemma}
    Let $X$ be an $\varepsilon$-fibrant $\varepsilon$-simplicial set and $f\colon A\hookrightarrow B$ be the morphism~\eqref{eq:alg}. Then $f\boxslash X$.
\end{lemma}
\begin{proof}
    We will even prove $$\big(A\hookrightarrow \Sim^{\{0,1,2,3,4\}}\cup\Sim^{\{0,3\}}\cup\Sim^{\{1,3\}}\cup\Sim^{\{1,4\}}\big)\boxslash X.$$ Lemma~\ref{lem:fibHelp} then gives the desired result.
    As the first step we fill $\Sim^{\{0,1,2,3\}}$. Then we can fill $\Sim^{\{0,1,2,4\}}$ as well. At this point we have $\partial_3\Sim^4\cup\partial_4\Sim^4$. By Theorem~\ref{thm:lifting1}, the inclusion
    \begin{eqnarray*}
        A\hookrightarrow \Sim^{\{0,1,2,3,4\}}\cup\Sim^{\{0,3\}}\cup\Sim^{\{1,3\}}
    \end{eqnarray*}
    belongs to $\overline{\eHorn}$. The last witness of the $\varepsilon$-edge $\Sim^{\{1,4\}}$ is obtained by the $4$-for-$3$ condition (Theorem~\ref{thm:4for3}).
\end{proof}
\begin{lemma}\label{lem:defOfEtaEdge}
    Let $a$ be an edge of an $\varepsilon$-fibrant $X$. Then the following are equivalent
    \begin{enumerate}
        \item[(i)] There is a $2$-simplex $\sigma\colon \Delta^2\rightarrow X$ such that $d_2\sigma = a$ and $d_1\sigma, d_0\sigma$ are $\varepsilon$-edges in $X$.
        \item[(ii)] There is a $2$-simplex $\sigma\colon \Delta^2\rightarrow X$ such that $d_0\sigma = a$ and $d_1\sigma, d_2\sigma$ are $\varepsilon$-edges in $X$.
    \end{enumerate}
\end{lemma}
\begin{proof}
    Assume (i) holds. Since $\eHorn\boxslash X$, there is a $2$-simplex $\tau\colon \Sim^2\rightarrow X$, such that $d_0\tau = a = d_2\sigma$.
    The two simplices $\sigma$ and $\tau$ form a rhombus, with three $\varepsilon$-edges.
By (L2), the edge $d_2\tau$ is an $\varepsilon$-edge,
so $\tau$ witnesses (ii). The converse follows by duality.
\end{proof}
\begin{defn}\label{def:etaEdge}
    Let $a$ be an edge of an $\varepsilon$-simplicial set $X$. We say $a$ is an $\eta$-edge if it satisfies one of the conditions (i--ii) in Lemma~\ref{lem:defOfEtaEdge}.
\end{defn}
\begin{theorem}\label{thm:etaHorn}
    Let $X$ be $\varepsilon$-fibrant and $\Lambda^n_i\rightarrow X$, $n\geq 2$, be a horn where $\Delta^{\{i,i+1\}}\hookrightarrow\Lambda^n_i\rightarrow X$
    is an $\eta$-edge (if $i<n$) or $\Delta^{\{i-1,i\}}\hookrightarrow\Lambda^n_i\rightarrow X$
    is an $\eta$-edge (if $0< i$). Then the horn admits a filling in $X$.
\end{theorem}
\begin{proof}
    Using Proposition~\ref{prop:cohRot} and Lemma~\ref{lem:defOfEtaEdge}, we can repeatedly rotate the horn until the involved $\eta$-edge rotates to an $\varepsilon$-edge. The rotation is isomorphic to an $\varepsilon$-horn $\ELambda^n_n\rightarrow X$ or $\ELambda^n_0\rightarrow X$.
\end{proof}
\begin{lemma}\label{lem:3for2}
    Let $X\in\eSSet$ be $\varepsilon$-fibrant and $a,b,c$ be edges of a $2$-simplex in $X$. If any two of the edges $a,b,c$ are $\eta$-edges, the last edge is an $\eta$-edge as well.
\[
    \begin{tikzcd}[row sep={1.8em,between origins}, column sep={2.4em,between origins}, ampersand replacement=\&]
    \&
    2\&
    \\[5pt]
    \&\&\\
    0
    \arrow[ruu, "c"{description}]
    \arrow[rr,"a"{description}]
    \&\&
    1\arrow[uul, "b"{description}]
\end{tikzcd}
\]
\end{lemma}
\begin{proof}
    First assume $a,b$ are $\eta$-edges. We have the following configuration of three $2$-simplices:
\[
          \begin{tikzcd}[row sep={2.4em,between origins}, column sep={2.6em,between origins}, ampersand replacement=\&]
    \&
    {3}
    \&\\[10pt]\&
    {2}\arrow[u,LA=1.5pt]
    \&\\
    {0}
    \arrow[rr,"a"{description}, dotted, bend right=5]
    \arrow[ru, "c"{description}]
    \&\&
    {1}\arrow[ul, "b"{description},dotted]\arrow[luu,LA=1.5pt]\\
    \&{x}\arrow[ul,LA=1.5pt]\arrow[ur,LA=1.5pt]\&
\end{tikzcd}
\]
Applying Theorem~\ref{thm:etaHorn} twice, first to the horn
on $\langle0,1,3\rangle$ and then to $\Lambda^3_1$,
we obtain a $3$-simplex on $\langle0,1,2,3\rangle$.
Theorem~\ref{thm:4for3}, applied to
$\langle x,0,1\rangle$ and $\langle0,1,3\rangle$,
shows that $\langle0,3\rangle$ is an $\varepsilon$-edge.
Thus $\langle0,2,3\rangle$ witnesses that $c$ is an $\eta$-edge.

Next assume that $a,c$ are $\eta$-edges. Then we have the following configuration of three $2$-simplices:
\[
          \begin{tikzcd}[row sep={2.4em,between origins}, column sep={2.6em,between origins}, ampersand replacement=\&]
    \&
    {3}
    \&\\[10pt]\&
    {2}\arrow[u,LA=1.5pt]
    \&\\
    {0}
    \arrow[ruu,LA=1.5pt, bend left=5]
    \arrow[rr,"a"{description}, dotted, bend right=5]
    \arrow[ru, "c"{description},dotted]
    \&\&
    {1}\arrow[ul, "b"{description}]\\
    \&{x}\arrow[ul,LA=1.5pt]\arrow[ur,LA=1.5pt]\&
\end{tikzcd}
\]
Filling two $\varepsilon$-horns, first on
$\langle0,1,3\rangle$ and then on $\langle0,1,2,3\rangle$,
we obtain a $3$-simplex.
Applying $\mathrm{(L1)}\boxslash X$ to
$\langle x,0,1\rangle$ and $\langle0,1,3\rangle$,
we see that $\langle1,3\rangle$ is an $\varepsilon$-edge.
Thus $\langle1,2,3\rangle$ witnesses that $b$ is an $\eta$-edge.

The case where $b,c$ are $\eta$-edges follows by duality.
\end{proof}
\begin{corollary}
    Let $X$ be $\varepsilon$-fibrant and $K$ its sub-simplicial set of all the simplices whose edges are $\eta$-edges. Then $K$ is a Kan complex.
\end{corollary}
\begin{proof}
    Horns $\Lambda^n_i\rightarrow K$, $n\geq 2$, have fillers by Theorem~\ref{thm:etaHorn}. In the case $n=2$, Lemma~\ref{lem:3for2} guarantees that the fillers belong to $K$.
\end{proof}
For a quasi-category $C$, the homotopy category $h(C)$
admits an explicit description in terms of homotopy classes
of edges; see \cite[Section~1.6]{Cisinski2019}. Its objects are the vertices
of $C$, and its morphisms are homotopy classes of edges.
Two edges $a,b$ with the same source and target are homotopic,
written $a\simeq b$, if there is a $2$-simplex of the following
form, where the edge labelled $\mathrm{id}$ is degenerate:
 \begin{equation}
        \begin{tikzcd}[row sep = 0.3em, column sep=3.6em]
            &{1}\\
            {0}\arrow[ur,"b"{description}]\arrow[dr,"a"{description}]\\
            &{1}\arrow[uu,"\mathrm{id}"{description}]
        \end{tikzcd}
    \end{equation}
One can prove that the relation $\simeq$ defined in this way is an equivalence relation which preserves composition of arrows.

The following lemma introduces an equivalence relation on edges
that may play a similar role in the study of Frobenius monoids
via $\varepsilon$-simplicial sets. Here we only establish that
it is an equivalence relation; its interaction with the algebraic
structure is left for future work.
\begin{lemma}
    Let $X\in\eSSet$ be $\varepsilon$-fibrant. Define a binary relation $\sim$ on the edges $X_1$ as follows. Given two edges $a,b\in X_1$ we set $a\sim b$ if and only if in $X$ we have the following pair of $2$-simplices
    \begin{equation}
        \begin{tikzcd}[row sep = 1.8em, column sep=3.6em]
            {1}\arrow[r,"b"{description}]&{3}\\
            {0}\arrow[u,dotted,"n_1"'{description}]\arrow[r,"a"{description}]\arrow[ru]&{2}\arrow[u,dotted,"n_2"{description}]
        \end{tikzcd}
    \end{equation}
    where $n_1,n_2$ are $\eta$-edges.
Then $\sim$ is an equivalence relation on $X_1$.
\end{lemma}
\begin{proof}
    Reflexivity is easy, since degenerate edges in $X$ are $\eta$-edges. Using Theorem~\ref{thm:etaHorn} and Lemma~\ref{lem:3for2}, it is straightforward to show transitivity and symmetry.
\end{proof}
\begin{example}
Recall Example~\ref{ex:orthomodular}. The quotient construction
identifies two events $A,B$ if and only if they admit a common
local orthocomplement $C$. This is represented by the following
pair of $2$-simplices:
    \[
          \begin{tikzcd}[row sep={1.8em,between origins}, column sep={3.2em,between origins}, ampersand replacement=\&]
    \&
    {3}
    \&\\[10pt]\&
    {2}
    \&\\
    {0}
    \arrow[rr,"C"{description}, bend right=5]
    \arrow[ruu, LA=1.5pt]\arrow[ru,LA=1.5pt]
    \&\&
    {1}\arrow[ul, "A"{description}]\arrow[luu,"B"{description}]
\end{tikzcd}
\]
It is easy to see that in an $\varepsilon$-fibrant $X$, such a configuration implies
$A\sim B$.

Note that the nerve of an effect test space with two
distinct tests $S,T$ is not $\varepsilon$-fibrant:
a filler of the two-dimensional $\varepsilon$-horn
with edges $\langle0,2\rangle$ and $\langle0,1\rangle$
encoding $T$ and $S$, respectively, would require
an event $E$ with $S+E=T$, contrary to irredundancy.
In any $\varepsilon$-fibrant extension, a filler of this horn provides a triangle
whose third edge is an $\eta$-edge.
This suggests a homotopical interpretation of the
identifications involved in the effect algebra quotient.
\end{example}

\section{Cycles}
Let us introduce some notation. For $k\geq 0$, let $C_\infty[k]$ be the $\varepsilon$-simplicial set
\begin{align}\label{eq:infCyc}
    C_\infty[k]&=\bigcup_{i\in \mathbb{Z}}\Sim^{\{i,i+1,\ldots,i+k+1\}}
\end{align}
By $\partial C_\infty[k]\subset C_\infty[k]$ we denote the $\varepsilon$-simplicial set
\begin{equation}\label{eq:infCycBoundary}
    \partial C_\infty[k]=\bigcup_{j\in \mathbb{Z},1\leq i\leq k}\Sim^{\{j,j+1,\ldots\widehat{j+i},\ldots,j+k+1\}}.
\end{equation}
Figure~\ref{fig:snake} illustrates (part of) the $\varepsilon$-simplicial set $C_\infty[1]$. The marked horizontal arrows form $\partial C_\infty[1]$.
\begin{figure}[H]
     \caption{$C_\infty[1]$}
     \label{fig:snake}
\begin{equation*}
\cdots\begin{tikzcd}[row sep=2.2em, column sep=1.6em]
&0\arrow[rr,LA=1.5pt]\arrow[rd] &&2\arrow[rr,LA=1.5pt]\arrow[rd] &&4\arrow[rr,LA=1.5pt]\arrow[rd] &&6\arrow[rr,LA=1.5pt]\arrow[rd] &&8\\
{-1}\arrow[rr,LA=1.5pt]\arrow[ru]&&{1}\arrow[rr,LA=1.5pt]\arrow[ru]&&{3}\arrow[rr,LA=1.5pt]\arrow[ru]&&{5}\arrow[rr,LA=1.5pt]\arrow[ru]&&{7}\arrow[ru]&
\end{tikzcd}\cdots
\end{equation*}
\end{figure}

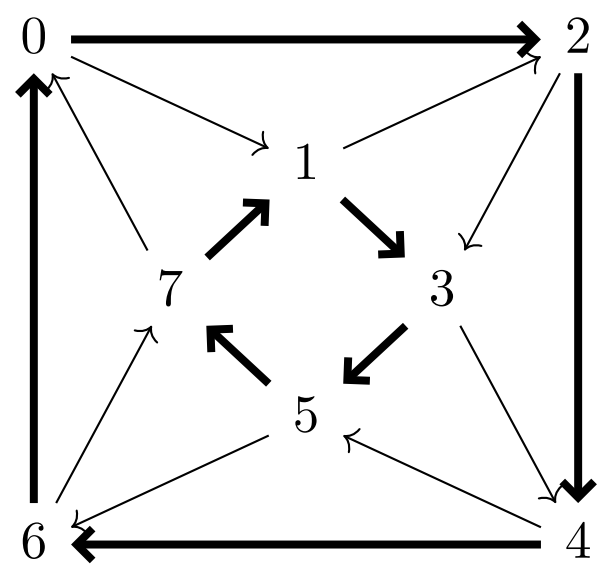
\begin{figure}
    \centering
    \begin{equation*}
\begin{tikzcd}[row sep=1em, column sep=1.2em]
    {0}\arrow[rrrr,LA=1.5pt]\arrow[rrd] &&&&{2}\arrow[dddd,LA=1.5pt]\arrow[ldd]\\
    &&{1}\arrow[rru]\arrow[rd,LA=1.5pt]&&\\
    &{7}\arrow[ru,LA=1.5pt]\arrow[luu]&&{3}\arrow[rdd]\arrow[ld,LA=1.5pt]&\\
     &&{5}\arrow[lu,LA=1.5pt]\arrow[lld]&&\\
     {6}\arrow[uuuu,LA=1.5pt]\arrow[ruu]&&&&{4}\arrow[llll,LA=1.5pt]\arrow[llu]
\end{tikzcd}
\end{equation*}
    \caption{$\varepsilon$-simplicial set $C_8[1]$}
    \label{fig:cyc}
\end{figure}

Let $l\in\mathbb{N}$ and $\sigma,\sigma'\colon\Delta^d\rightarrow C_\infty[k]$ be two $d$-simplices of $C_\infty[k]$
\begin{equation}
    \sigma=\langle i_0,i_1,\ldots,i_d\rangle\quad\sigma'=\langle i'_0,i'_1,\ldots,i'_d\rangle
\end{equation}
where $i_0\leq i_1\leq\cdots\leq i_d\leq i_0+k+1$ and  $i'_0\leq i'_1\leq\cdots\leq i'_d\leq i'_0+k+1$. Write $\sigma\equiv\sigma' \pmod{l}$,
if there is $r\in\mathbb{Z}$ such that $i'_j= i_j+rl$ for each $j=0,\ldots,d$.
This defines an equivalence relation on each simplicial level. Since it is compatible with the face and degeneracy maps, the levelwise quotient defines an $\varepsilon$-simplicial set
\begin{eqnarray}
    C_l[k]=C_\infty[k]/\equiv.
\end{eqnarray}
We similarly quotient the set of $\varepsilon$-witnesses by the equivalence relation induced by the composite $C_\infty[k]_\varepsilon\rightarrow C_\infty[k]_1\rightarrow C_l[k]_1$.
See Figure~\ref{fig:cyc} for an illustration of $C_8[1]$.

Moreover, we define the $\varepsilon$-simplicial set $\partial C_l[k]$ as
\[\partial C_l[k]=\mathrm{Im}(\partial C_\infty[k]\rightarrow C_l[k]).\]
Note that two simplices $\sigma=[\langle i_0,\ldots,i_d\rangle]$ and $\sigma'=[\langle i'_0,\ldots,i'_d\rangle]$ in $C_l[k]$ may be distinct even if
$i_j\equiv i'_j \pmod{l}$ for each $j=0,\ldots,d$. Indeed, the definition requires a single $r$ such that $i'_j= i_j+rl$ for all $j$. For example, $C_1[k]$ has only one vertex and only one nondegenerate $k+1$-simplex.
However, for $0<d<k+1$, it may contain distinct nondegenerate $d$-simplices.

For some $l$ and $k$ the morphism $C_\infty[k]\rightarrow C_l[k]$ creates a lot of identifications.
However, it locally restricts to a mono on simplices with large enough dimension. This is the content of the
following lemma which will be useful in Section~\ref{sec:stability}.
\begin{lemma}\label{lem:PulbackAdPushout}
    Let $l\geq 1$ and $k\geq 0$ be integers and consider the morphism
    \begin{equation}
        \pi\colon \Sim^{k+1}\hookrightarrow C_\infty[k]\rightarrow C_l[k].
    \end{equation}
    \begin{enumerate}
        \item[(i)]   The following diagram is a pullback:
         \begin{equation}
            \begin{tikzcd}
                \ELambda^{k+1}_{0,k+1}\arrow[d,hook]\arrow[r]&\partial C_l[k]\arrow[d,hook]\\
                \Sim^{k+1}\arrow[r,"\pi"]&C_l[k]
            \end{tikzcd}
        \end{equation}
    where
    $$\ELambda^{k+1}_{0,k+1}=\bigcup_{1\leq j\leq k}\partial_j\Sim^{k+1}\cup (\partial_0\Sim^{k+1}\cap\partial_{k+1}\Sim^{k+1})$$
    is obtained from $\partial\Sim^{k+1}$ by deleting the interiors of the two faces $\partial_0\Sim^{k+1}$ and $\partial_{k+1}\Sim^{k+1}$.

  \item[(ii)]  If $l\geq 2$, then $\pi$ restricts to a levelwise embedding:
   \begin{equation}
            \Sim^{k+1}\setminus\ELambda^{k+1}_{0,k+1}\hookrightarrow C_l[k]\setminus\partial C_l[k].
        \end{equation}
    \item[(iii)] For each $l\geq 1$, the morphism $\pi$ restricts to a levelwise embedding:
    $$\Sim^{k+1}\setminus\ELambda^{k+1}_0\hookrightarrow C_l[k]$$
    \end{enumerate}
\end{lemma}
\begin{proof}
Let
\(
A=\pi^*\partial C_l[k]\subseteq\Sim^{k+1}.
\)
By the definition of $\partial C_l[k]$, we have
\[
\bigcup_{1\leq r\leq k}\partial_r\Sim^{k+1}
\cup
\bigl(\partial_0\Sim^{k+1}\cap\partial_{k+1}\Sim^{k+1}\bigr)
\subseteq A.
\]
Conversely, a simplex
\(
\tau=\langle i_0,\ldots,i_d\rangle
\notin\ELambda^{k+1}_{0,k+1}
\)
contains each of the vertices $1,\ldots,k$ as well as at least one of the vertices $0$ or $k+1$. By the definition of
$\partial C_l[k]$, this implies that $\pi(\tau)$ does not belong to
$\partial C_l[k]$. Hence
\[
A=\ELambda^{k+1}_{0,k+1},
\]
which proves (i).

It remains to consider possible identifications under $\pi$.
Let
\[
\tau,\rho\in
\Sim^{k+1}\setminus\ELambda^{k+1}_{0,k+1}
\]
and suppose that $\pi(\tau)=\pi(\rho)$. Both $\tau$ and $\rho$
contain all vertices $1,\ldots,k$. A nontrivial identification in
$C_l[k]$ is given by translation by a nonzero multiple of $l$.
Up to interchanging $\tau$ and $\rho$, suppose that $\rho$ is obtained
from $\tau$ by translation by $ql$, where $q\geq1$.
Since $\rho$ contains at least $k+1$ distinct vertices
in $\{0,\ldots,k+1\}$, its least vertex is at most $1$.
As $\rho$ is obtained from $\tau$ by translation
by $ql$, we obtain $ql\leq1$. Thus necessarily
$q=l=1$. In particular, for $l\geq2$ no nontrivial identification is
possible, proving (ii).

For $l=1$, the possible identifications are of the form
\[
\langle0,\ldots,k\rangle
\equiv
\langle1,\ldots,k+1\rangle \pmod 1.
\]
If moreover
\(
\tau,\rho\in
\Sim^{k+1}\setminus\ELambda^{k+1}_0,
\)
then both simplices contain the vertex $k+1$, so this exceptional
identification cannot occur. This proves (iii).
\end{proof}
We conclude the section by introducing one additional piece of notation. Set $$C_\infty[\varepsilon]=C_\infty[1]\cup \bigcup_{i}\Sim^{\{2i,2i+1\}}.$$
In $C_\infty[\varepsilon]$, the edge $\langle0,1\rangle$ is an $\varepsilon$-edge and so are all its even rotations, while odd rotations are $\eta$-edges. Similarly set
$C_l[\varepsilon]$ to be an $\varepsilon$-simplicial set $C_l[1]$ with additional $\varepsilon$-edges $[\langle 2i,2i+1\rangle]$, $i\in\mathbb Z$.
\section{Relation to paracyclic sets}\label{sec:6}
In this section we introduce the category of paracyclic sets $\mathrm{Set}_{\Lambda_\infty}=\mathrm{Fun}(\Lambda_\infty^{\mathrm{op}},\mathrm{Set})$ and explore its relation
 to $\eSSet$ through a chain of adjunctions $P_!\dashv P^*\dashv P_*$, where $P^*\colon \mathrm{Set}_{\Lambda_\infty}\rightarrow\eSSet$.
In particular, we identify the $\varepsilon$-simplicial
  sets $C_l[n]$ introduced in the previous section with the restrictions under $P^*$ of suitable quotients of representable paracyclic sets.

  By $\Lambda_\infty$ we denote the \emph{parasimplex} category with the same objects as $\Delta$ but morphisms from $[m]$ to $[n]$ correspond to maps $f\colon\mathbb{Z}\rightarrow\mathbb{Z}$ such that
\begin{enumerate}[(i)]
    \item if $i\leq j$ then $f(i)\leq f(j)$ (order-preserving),
    \item $f(i+m+1)=f(i)+n+1$.
\end{enumerate}
Let $f\colon\mathbb{Z}\rightarrow\mathbb{Z}$ satisfy (i--ii), so that it represents a morphism $[m]\rightarrow [n]$. By (ii), $f$ is fully determined by its values on any $m+1$-tuple of consecutive integers, in particular by $f(0),\ldots,f(m)$. Moreover, due to (ii-i) we have
\begin{equation}\label{eq:paracycmorii}
    f(0)\leq f(1)\leq\cdots\leq f(m)\leq f(m+1)=f(0)+n+1.
\end{equation}
Conversely, each $(m+1)$-tuple of integers $$\langle f(0),\ldots,f(m)\rangle$$ satisfying~\eqref{eq:paracycmorii} determines a unique morphism $f\colon[m]\rightarrow [n]$.
We regard $\Delta$ as a subcategory of $\Lambda_\infty$ by identifying
a morphism $[m]\to[n]$ in $\Delta$ with the unique morphism
$f\colon\mathbb Z\to\mathbb Z$ satisfying {\rm (i)--(ii)} and
\[
f(0),\ldots,f(m)\in\{0,\ldots,n\}.
\]

It is easy to see that $\Lambda_\infty$ is generated by $\Delta$ and morphisms $\tau_n,\tau_n^{-1}\colon[n]\rightarrow [n]$, $n\geq 0$, given by $\tau_n\colon i\mapsto i+1$ and $\tau_n^{-1}\colon i \rightarrow i-1$. Moreover, for each $n\geq 0$, the following morphism is important
\begin{equation}\label{eq:seqOfDegeneracy}
    \sigma^{n}_{n+1}\colon[n+1]\rightarrow[n],\quad \text{given\ by}\quad \langle 0,1,2,\ldots,n,n+1\rangle.
\end{equation}
It does not belong to $\Delta$, since its last value $n+1$
lies outside $[n]$.
A \emph{paracyclic set} is a presheaf $X\colon\Lambda_\infty^{op}\rightarrow\mathrm{Set}$. Equivalently, we can define a paracyclic set as a simplicial set $X$ together with a system of bijections $t_n\colon X_n\rightarrow X_n$, $n\geq 0$,
such that the following identities hold:
\begin{equation}\label{eq:parBound}
    d_i t_n =
\begin{cases}
t_{n-1} d_{i+1}& 0\leq i < n,\\
d_0 & i=n,
\end{cases}
\end{equation}

\begin{equation}\label{eq:parDeg}
s_i t_n =
\begin{cases}
t_{n+1}s_{i+1} & 0\leq i < n,\\
t_{n+1}^2s_0  & i=n.
\end{cases}
\end{equation}
Moreover, we denote by $$s_{n+1}^n\colon X_n\rightarrow X_{n+1}$$ the extra degeneracy map induced by~\eqref{eq:seqOfDegeneracy}. We can express the rotations using the extra degeneracy and vice versa:
\begin{equation}
    t_n=d_0s_{n+1}^n\quad\text{and}\quad s_{n+1}^n=t_{n+1}s_0^n.
\end{equation}
We say that an $n$-simplex of a paracyclic set $X$ is \emph{degenerate}, if it is in the image of any $s_i^{n-1}$, for $0\leq i\leq n$. Note that the rotations of degenerate simplices are degenerate.
Consider the functor
\[
P\colon\Sim\longrightarrow\Lambda_\infty
\]
whose restriction to $\Delta\subseteq\Sim$ is the inclusion
$\Delta\subseteq\Lambda_\infty$, with $P([\epsilon])=[0]$, and which
sends the distinguished morphism $[1]\to[\epsilon]$ to
$\sigma^0_1\colon[1]\to[0]$.

The restriction functor $P^*$ admits both adjoints,
\[
P_!\dashv P^*\dashv P_*,
\]
given by left and right Kan extension along $P^{op}$.

Using notation from the previous section, we can represent the units $\eta_{\Delta^n}$ and $\eta_{\partial\Delta^n}$ in concrete terms:
\begin{lemma}\label{lem:weNeed}
    For each $n\geq 0$, there are isomorphisms of $\varepsilon$-simplicial sets
    $P^*P_!(\Delta^n)\cong C_\infty[n]$ and $P^*P_!(\partial\Delta^n)\cong \partial C_\infty[n]$ which fit into the following commutative diagram:
\begin{equation}
    \begin{tikzcd}[column sep=4.5em, row sep=2.4em]
        {\partial\Delta^n}\arrow[rr,"{\eta_{\partial\Delta^n}}",bend left = 30]\arrow[r,hook]\arrow[d,hook,"i"]&{\partial C_\infty[n]}\arrow[r,"\cong"]\arrow[d,hook]&{P^*P_!(\partial\Delta^n)}\arrow[d,"{P^*P_!(i)}"]\\
        {\Delta^n}\arrow[rr,"\eta_{\Delta^n}",bend right =30]\arrow[r,hook]&{C_\infty[n]}\arrow[r,"\cong"]&{P^*P_!(\Delta^n)}
    \end{tikzcd}
\end{equation}
\end{lemma}

\begin{proof}
Since
\[
P_!\circ y_{\Sim}=y_{\Lambda_\infty}\circ P,
\]
for every $k\geq 0$ we have natural isomorphisms
\begin{align}\label{eq:tractUnit}
\mathrm{Hom}(\Delta^k,P^*P_!(\Delta^n))
&\cong
\mathrm{Hom}(P_!(\Delta^k),P_!(\Delta^n))\\
&\cong
\mathrm{Hom}_{\Lambda_\infty}([k],[n])\\
&\cong
\mathrm{Hom}(\Delta^k,C_\infty[n]).
\end{align}
The last bijection sends a morphism $f\colon[k]\to[n]$, represented by
\[
f(0)\leq\cdots\leq f(k)\leq f(0)+n+1,
\]
to the corresponding simplex
\[
\langle f(0),\ldots,f(k)\rangle
\]
of $C_\infty[n]$.
Since $P([1]\rightarrow[\varepsilon])=\sigma_1^0\colon [1]\rightarrow [0]$, we can trace the restriction morphism $\mathrm{Hom}(\Sim^1,P^*P_!(\Delta^n))\rightarrow \mathrm{Hom}(\Delta^1,P^*P_!(\Delta^n))$
through the above isomorphisms to
$$(\sigma^0_1)^*\colon\mathrm{Hom}_{\Lambda_\infty}([0],[n])\rightarrow \mathrm{Hom}_{\Lambda_\infty}([1],[n]).$$
In the image are exactly the morphisms $f\colon [1]\rightarrow [n]$ such that $f(1)=f(0)+n+1$. These correspond precisely to $\varepsilon$-edges $\langle j,j+n+1\rangle$ of $C_\infty[n]$.

It follows by Yoneda that
\[
P^*P_!(\Delta^n)\cong C_\infty[n].
\]
Under this isomorphism, the unit $\eta_{\Delta^n}$ corresponds to
$\mathrm{id}_{[n]}$, hence to the canonical inclusion
\[
\Delta^n\hookrightarrow C_\infty[n].
\]

For each nonempty subset $I\subseteq[n]$, the corresponding
face inclusion induces a monomorphism
$y[|I|-1]\rightarrow y[n]$ whose image consists of morphisms
with values modulo $n+1$ in $I$.
These images intersect according to intersections of
the subsets $I$.
Thus, since $P_!$ preserves colimits, we have
\[
P_!(\partial\Delta^n)
=
\bigcup_{i=0}^n P_!(\partial_i\Delta^n)
\subseteq y[n].
\]
Under the above identification, its $k$-simplices correspond precisely
to those morphisms $f\colon[k]\to[n]$ which factor through a simplicial
coface. Equivalently, the values of $f$ omit one residue class modulo
$n+1$, or, equivalently, do not contain $n+1$ consecutive integers.
These are exactly the $k$-simplices of $\partial C_\infty[n]$. Hence
\[
P^*P_!(\partial\Delta^n)\cong\partial C_\infty[n].
\]
The two isomorphisms are compatible with the inclusions and with the
units, giving the commutative diagram in the statement.
\end{proof}

Recall that every simplicial set can be obtained as a transfinite
composition of pushouts of the boundary inclusions
\[
\partial\Delta^n\hookrightarrow\Delta^n.
\]
There is an analogous cellular decomposition for paracyclic sets.
The paracyclic category belongs to the Eilenberg--Zilber framework of
Berger and Moerdijk; see \cite[Section~6]{BM2011}. More generally,
cellular decompositions of this kind are studied for elegant generalized
Reedy categories; see \cite[Theorem~5.46]{CaSa2025}. We will only need
the following explicit description in the paracyclic case.

For $[n]\in\Lambda_\infty$ and $l\geq1$, the automorphism
$\tau_n^l$ induces an action of the subgroup
\[
\langle\tau_n^l\rangle
\leq
\operatorname{Aut}_{\Lambda_\infty}([n])
\]
on the representable paracyclic set $y[n]$. We denote the corresponding
orbit presheaf by
\[
y[n]/l
:=
y[n]/\langle\tau_n^l\rangle.
\]
Since colimits of presheaves are computed pointwise, two $d$-simplices
of $y[n]$ represented by
\[
\sigma_1=\langle f(0),\ldots,f(d)\rangle,
\qquad
\sigma_2=\langle g(0),\ldots,g(d)\rangle
\]
are identified in $y[n]/l$ if and only if
\[
\sigma_1\equiv\sigma_2\pmod l.
\]

Similarly, we set
\[
\partial y[n]/l
:=
(\partial y[n])/\langle\tau_n^l\rangle.
\]
Under the identifications from the previous section, we then have
\begin{equation}\label{eq:repOfCycles}
P^*(y[n]/l)\cong C_l[n],
\qquad
P^*(\partial y[n]/l)\cong\partial C_l[n].
\end{equation}

\begin{lemma}\label{lem:genOfMonoInPar}
Let $m\colon X\hookrightarrow Y$ be a monomorphism of paracyclic sets.
Then $m$ is a transfinite composition of pushouts of monomorphisms of
the form
\[
\partial y[n]\hookrightarrow y[n]
\qquad\text{and}\qquad
\partial y[n]/l\hookrightarrow y[n]/l,
\]
where $n\geq 0$ and $l\geq 1$.
\end{lemma}
\begin{proof}
Identify $X$ with a subobject of $Y$. We construct $Y$ from $X$
by attaching orbits of nondegenerate simplices in increasing
dimension. Let $Z\subseteq Y$ be the current stage, containing
all nondegenerate simplices of dimension less than $n$, and
choose a nondegenerate $x\in Y_n\setminus Z_n$.

Let $H\leq\operatorname{Aut}_{\Lambda_\infty}([n])$ be the
stabilizer of $x$. Regarding $x$ as a morphism $y[n]\to Y$,
its component in degree $k$ is given by
\[
x_k(\rho)=Y(\rho)(x)=:\rho^*x,
\qquad \rho\colon[k]\to[n].
\]
It induces a map $\bar x\colon y[n]/H\to Y$,
which sends $\partial y[n]/H$ into $Z$.
We show that the canonical map
\[
Z\coprod_{\partial y[n]/H}y[n]/H
\longrightarrow Z\cup\operatorname{Im}(x)
\]
is an isomorphism.

An element of $y[n]_k$ outside $\partial y[n]_k$ is represented
by a morphism $\rho\colon[k]\to[n]$ whose underlying map
$\mathbb Z\to\mathbb Z$ is surjective. Choosing preimages
compatibly with the order and periodicity gives a morphism
$\delta\colon[n]\to[k]$ such that $\rho\delta=\mathrm{id}_{[n]}$.

For two such morphisms $\rho,\rho'$, Proposition~6.9
of~\cite{BM2011}, applied to the nondegenerate simplex $x$,
gives
\[
\rho^*x=(\rho')^*x
\quad\Longleftrightarrow\quad
\rho'=h\rho \text{ for some }h\in H.
\]
Thus $\bar x$ is levelwise injective on the complement
of $\partial y[n]/H$.
Furthermore, if $\rho^*x$ belonged to $Z_k$, then
\[
x=(\rho\delta)^*x=\delta^*(\rho^*x)
\]
would belong to $Z_n$, since $Z$ is a subpresheaf.
This contradicts the choice of $x$.

Consequently, the complement of $\partial y[n]/H$ in $y[n]/H$
maps bijectively, in every degree, onto the complement of $Z$
in $Z\cup\operatorname{Im}(x)$.
Since pushouts of presheaves are computed pointwise,
the claimed map is an isomorphism.

Finally, $\operatorname{Aut}_{\Lambda_\infty}([n])
=\langle\tau_n\rangle\cong\mathbb Z$, so $H$ is either trivial
or equal to $\langle\tau_n^l\rangle$ for some $l\geq1$.
The attaching map is therefore one of the stated forms.
Well-ordering the missing orbits in each dimension and
proceeding through all dimensions gives the required
transfinite composition.
\end{proof}

We now turn to the right adjoint $P_*$. Denote the counit of
$P^*\dashv P_*$ by
\[
\epsilon'_X\colon P^*P_*X\longrightarrow X.
\]
Observe that the $\varepsilon$-simplicial set $C_\infty[n]$ is naturally identified with $P^*\Lambda_\infty[n]$.
Since $y_{\Lambda_\infty}\circ P=P_!\circ y_{\Sim}$, we have $C_\infty[n]\cong P^*\Lambda_\infty[n]=P^*P_!\Delta^n$. It follows that
\begin{align}
    (P^*P_*(X))_n&=\mathrm{Hom}(\Delta^n,P^*P_*(X))\cong\mathrm{Hom}(P_!(\Delta^n),P_*(X))\\
    &\cong \mathrm{Hom}(P^*P_!(\Delta^n),X)\cong\mathrm{Hom}(C_\infty[n],X).\label{eq:mimicThis}
\end{align}
Let $f\colon Y\rightarrow P^*P_*(X)$ be a morphism for $X,Y\in\eSSet$. By adjunction, $f$ corresponds to a morphism $P_!(Y)\rightarrow P_*(X)$, which in turn corresponds to a morphism $\bar{f}\colon P^*P_!(Y)\rightarrow X$.
We have the following commutative diagram
\begin{equation}
    \begin{tikzcd}
        Y\arrow[r,"f"]\arrow[d, "\eta_Y"]&P^*P_*(X)\arrow[d,"\epsilon'_X"]\\
        P^*P_!(Y)\arrow[r,"\bar{f}"]&X
    \end{tikzcd}
\end{equation}

It follows that
\begin{equation}
\begin{tikzcd}
P^*P_*(X)_n
  \arrow[r, "\cong"]
  \arrow[d, "(\epsilon'_X)_n"']
&
\mathrm{Hom}(\Delta^n,P^*P_*(X))
  \arrow[r, "\cong"]
  \arrow[d, "{\mathrm{Hom}(\Delta^n,\epsilon'_X)}"']
&
\mathrm{Hom}(P^*P_!(\Delta^n),X)
  \arrow[d, "{\mathrm{Hom}(\eta_{\Delta^n},X)}"']
\\
X_n
  \arrow[r, "\cong"]
&
\mathrm{Hom}(\Delta^n,X)
  \arrow[r, "="]
&
\mathrm{Hom}(\Delta^n,X)
\end{tikzcd}
\end{equation}
Thus, under the identification provided by Lemma~\ref{lem:weNeed} the counit is given by restriction along the canonical inclusion
\begin{equation}\label{eq:counitDesc}
    (\epsilon'_X)_n\colon \mathrm{Hom}(C_\infty[n],X)\rightarrow \mathrm{Hom}(\Delta^{\{0,\ldots,n\}},X).
\end{equation}

We may therefore think of an $n$-simplex of $P^*P_*X$ as an
$n$-simplex of $X$ equipped with a coherent choice of all its rotations. Indeed, each $n$-simplex
\(a\in (P^*P_*X)_n\) admits a canonical choice of its rotation $t_n(a)$, where we can represent $t_n$ as follows:
\begin{equation}\label{diag:canonicalRot}
    \begin{tikzcd}
        (P^*P_*(X))_n\arrow[r,"\cong"]\arrow[d,"{t_n}"']&\mathrm{Hom}(C_\infty[n],X)\arrow[d,"{\mathrm{Hom}(\tilde{\tau}_n,X)}"]\\
        (P^*P_*(X))_n\arrow[r,"\cong"]&\mathrm{Hom}(C_\infty[n],X)
    \end{tikzcd}
\end{equation}
where $\tilde{\tau}_n\colon C_\infty[n]\rightarrow C_\infty[n]$ is the isomorphism of $\varepsilon$-simplicial sets given by a shift of indices $i\mapsto i+1$.

The following lemma shows that the rotations of a simplex $a\in (P^*P_*X)_n$ form an orbit in a canonical way and in particular we have $a \beta t_n(a)$.
\begin{equation}
    \begin{tikzcd}
        {C_\infty[n]\cong P^*P_!\Delta^n}\arrow[r,"P^*(\tilde{a})"]&{P^*P_*X}\\
        {\Delta^n}\arrow[u,"\eta_{\Delta^n}"]\arrow[ru,"a"]&
    \end{tikzcd}
\end{equation}

\begin{lemma}\label{lem:canonicalOrbits}
Let \(X\in\eSSet\), and set \(A=P^*P_*X\).
    Every \(n\)-simplex \(a\in A_n\) determines a canonical morphism
    \[
        \rho_a\colon C_\infty[n]\longrightarrow A
    \]
    such that, for every \(i\in\mathbb Z\), its restriction to
    \(\Delta^{\{i,\ldots,i+n\}}\subseteq C_\infty[n]\)
    is the \(i\)-th canonical rotation \(t_n^i(a)\).
\end{lemma}
\begin{proof}
    Let
\(
    a\colon \Delta^n\longrightarrow A=P^*P_*X
\)
be an \(n\)-simplex. Under the adjunction \(P_!\dashv P^*\),
the morphism \(a\) corresponds to
\[
    \tilde{a}\colon P_!(\Delta^n)\longrightarrow P_*X.
\]
Applying \(P^*\) and using the canonical isomorphism
\(
    P^*P_!(\Delta^n)\cong C_\infty[n],
\)
we obtain the morphism
\[
    \rho_a:=P^*(\tilde{a})\colon C_\infty[n]\longrightarrow A.
\]
By the unit identity,
\(
    P^*(\tilde{a})\eta_{\Delta^n}=a,
\)
so the restriction of \(\rho_a\) to the distinguished copy
\(\Delta^{\{0,\ldots,n\}}\subset C_\infty[n]\) is \(a\).

For \(i\in\mathbb Z\), let
\[
    j_i\colon\Delta^n\longrightarrow C_\infty[n]
\]
denote the inclusion of
\(\Delta^{\{i,\ldots,i+n\}}\).
Under the identification
\(C_\infty[n]\cong P^*P_!(\Delta^n)\), we have
\(
    j_i=P^*(\tau_n^i)\eta_{\Delta^n}.
\)
Hence the restriction \(b:=\rho_a j_i\) is
\[
    b
    =P^*(\tilde{a})P^*(\tau_n^i)\eta_{\Delta^n}
    =P^*(\tilde{a}\circ\tau_n^i)\eta_{\Delta^n}.
\]
Therefore, under the adjunction \(P_!\dashv P^*\), we have
\(
    \tilde{b}=\tilde{a}\circ\tau_n^i.
\)
Applying the adjunction \(P^*\dashv P_*\), we obtain
\[
    \bar{b}
    =\bar{a}\circ P^*(\tau_n^i).
\]
Identifying \(P^*(\tau_n^i)\) with
\(\widetilde\tau_n^i\colon C_\infty[n]\to C_\infty[n]\),
diagram~\eqref{diag:canonicalRot} gives
\[
    b=t_n^i(a).
\]
Thus \(\rho_a\) represents the whole canonical \(t_n\)-orbit of \(a\).

\end{proof}
Later on we show that for $\varepsilon$-fibrant $X$ the counit $\epsilon'_X\colon P^*P_*X\rightarrow X$ is a trivial Kan fibration. In particular, it has the right lifting property with respect to all monomorphisms of simplicial sets.
We can extend it to a morphism $\hat{\epsilon}_X\colon\widehat{X}\rightarrow X$ having the right lifting property with respect to all monomorphisms of $\varepsilon$-simplicial sets.
We define $\widehat{X}$ by mimicking the above formulas for $P^*P_*X$.

\begin{defn}
For each $X\in\eSSet$ define an
$\varepsilon$-simplicial set $\widehat X$ whose underlying simplicial
set agrees with that of $P^*P_*X$, while its set of
$\varepsilon$-witnesses is defined by
\[
\widehat X_\varepsilon
=
\mathrm{Hom}(C_\infty[\varepsilon],X).
\]
The structure map
\(
\widehat X_\varepsilon\longrightarrow\widehat X_1
\)
is given by restriction along
\[
C_\infty[1]\hookrightarrow C_\infty[\varepsilon].
\]
Moreover, define
\[
\widehat\epsilon_X\colon\widehat X\rightarrow X.
\]
On the underlying simplicial set it agrees with $\epsilon'_X$, while
on $\varepsilon$-witnesses it is given, via Yoneda, by restriction along
\(
\Sim^{\{0,1\}}
\hookrightarrow
C_\infty[\varepsilon]
\).
\end{defn}
We first show that the identity on $P^*P_*X$ extends to an embedding
\begin{equation}\label{eq:embJ}
    j\colon P^*P_*X \hookrightarrow \widehat{X}.
\end{equation}
\begin{lemma}\label{lem:defOfq}
    Under the identifications $P^*P_!\Delta^1\cong C_\infty[1]$ and $P^*P_!\Sim^1\cong C_\infty[0]$, the morphism $P^*P_!(\Delta^1\hookrightarrow\Sim^1)$ becomes the folding morphism
 \begin{equation}\label{eq:defOfq}
     q\colon C_\infty[1]\rightarrow C_\infty[0],\qquad 2k\mapsto k,\quad 2k+1\mapsto k+1.
 \end{equation}
\end{lemma}
\begin{proof}
    Under the canonical identification
\(P_!y\simeq yP\), the morphism
\(P_!(\Delta^1\hookrightarrow\Sim^1)\) is induced by
\[
P([1]\to[\varepsilon])=\sigma_1^0\colon[1]\to[0].
\]
According to~\eqref{eq:seqOfDegeneracy}, the latter is represented
by the map \(f\colon\mathbb Z\to\mathbb Z\) given by
\[
f(2i)=i,\qquad f(2i+1)=i+1.
\]
After applying \(P^*\), this is precisely the morphism~\eqref{eq:defOfq}.
\end{proof}
Moreover, denote by $q_\varepsilon\colon C_\infty[\varepsilon]\rightarrow C_\infty[0]$ the obvious extension of $q$ so that
\[
q=q_\varepsilon\iota.
\]
Now we have the following commutative diagram which establishes~\eqref{eq:embJ} on $\varepsilon$-witnesses.
\begin{equation}\label{eq:embeddingq}
    \begin{tikzcd}[column sep=1.5em]
{(P^*P_*X)_1}\arrow[r,"\cong"]&{\mathrm{Hom}(\Delta^1,P^*P_*X)}\arrow[r,"\cong"]&{\mathrm{Hom}(C_\infty[1],X)}\arrow[r,"\cong"]&{\widehat{X}_1}\\
{(P^*P_*X)_\varepsilon}\arrow[u]\arrow[r,"\cong"]&{\mathrm{Hom}(\Sim^1,P^*P_*X)}\arrow[u]\arrow[r,"\cong"]&{\mathrm{Hom}(C_\infty[0],X)}\arrow[u,"q^*"]\arrow[r,"q^*_\varepsilon"]&{\mathrm{Hom}(C_\infty[\varepsilon],X)}\arrow[u]\arrow[lu,"\iota^*"']\\
    \end{tikzcd}
\end{equation}
Since \(q_\varepsilon\) is levelwise surjective, precomposition
with \(q_\varepsilon\) is injective. Hence \(j_\varepsilon\),
and therefore \(j\), is a monomorphism.

Similarly, \(\varepsilon\)-witnesses admit a canonical double rotation.
Denote by $$\tilde{\tau}^2_\varepsilon\colon C_\infty[\varepsilon]\rightarrow C_\infty[\varepsilon]$$ the isomorphism given by a shift of indices $i\mapsto i+2$.
It defines a canonical action of double rotation $t^2_\varepsilon$ on $(\widehat{X})_\varepsilon$ given by the following commutative diagram
\begin{equation}\label{diag:canonicalEpsRot}
    \begin{tikzcd}
        (\widehat{X})_\varepsilon\arrow[r,"\cong"]\arrow[d,"{t_\varepsilon^2}"']&\mathrm{Hom}(C_\infty[\varepsilon],X)\arrow[d,"{\mathrm{Hom}(\tilde{\tau}^2_\varepsilon,X)}"]\\
        (\widehat{X})_\varepsilon\arrow[r,"\cong"]&\mathrm{Hom}(C_\infty[\varepsilon],X)
    \end{tikzcd}
\end{equation}
Here the superscript $2$ indicates the shift by two; no single rotation $t_\varepsilon$ is defined.
\begin{lemma}\label{lem:canonicalOrbitsEps}
   Every \(\varepsilon\)-witness
    \(e\in\widehat X_\varepsilon\) determines a canonical morphism
    \[
        \rho_e\colon C_\infty[\varepsilon]\longrightarrow\widehat X.
    \]
    If \(a\in\widehat X_1\) is the underlying edge of \(e\), then
    \(
        \rho_e|_{C_\infty[1]}=\rho_a.
    \)
    Moreover, the image
\(
    j_\varepsilon\bigl((P^*P_*X)_\varepsilon\bigr)
    \subseteq\widehat X_\varepsilon
\)
is invariant under \(t_\varepsilon^2\).
\end{lemma}

\begin{proof}
Let
\(
e\in\widehat X_\varepsilon
=\operatorname{Hom}(C_\infty[\varepsilon],X)
\)
and \(a\in\widehat X_1\) be its underlying edge. Denote by
\[
    \iota\colon C_\infty[1]\hookrightarrow C_\infty[\varepsilon]
\]
the canonical inclusion. Under the identification
\(
    \widehat{X}_1\cong\operatorname{Hom}(C_\infty[1],X),
\)
the edge \(a\) is represented by
\(
    \bar a=e\circ\iota.
\)

By Lemma~\ref{lem:canonicalOrbits}, \(a\) determines the canonical orbit morphism
\[
    \rho_a\colon C_\infty[1]\longrightarrow P^*P_*X
    \hookrightarrow\widehat X,
\]
whose restriction to the edge \(\langle i,i+1\rangle\) is \(t_1^i(a)\). We show that $\rho_a$ extends to
the desired $\rho_e$. It is enough to check that for each $i\in\mathbb{Z}$, the underlying edge of $(t^2_\varepsilon)^i(e)$ equals $t_1^{2i}(a)$.

Recall that \(t^2_\varepsilon\) is induced by precomposition
with the shift
\[
    \widetilde\tau_\varepsilon^{\,2}\colon
    C_\infty[\varepsilon]\longrightarrow C_\infty[\varepsilon],
    \qquad j\longmapsto j+2.
\]
For each \(i\in\mathbb Z\), we have a commutative diagram
\[
\begin{tikzcd}
C_\infty[1]
    \arrow[r,hook,"\iota"]
    \arrow[d,"\widetilde\tau_1^{\,2i}"']
&
C_\infty[\varepsilon]
    \arrow[d,"(\widetilde\tau_\varepsilon^{\,2})^i"]
\\
C_\infty[1]
    \arrow[r,hook,"\iota"]
&
C_\infty[\varepsilon].
\end{tikzcd}
\]
Hence the underlying edge of
\((t_\varepsilon^{2})^i(e)\) is represented by
\[
e\circ(\widetilde\tau_\varepsilon^{\,2})^i\circ\iota
=
e\circ\iota\circ\widetilde\tau_1^{\,2i}
=
\bar a\circ\widetilde\tau_1^{\,2i},
\]
and therefore equals \(t_1^{2i}(a)\) by Lemma~\ref{lem:canonicalOrbits}. Thus these witnesses, together with \(\rho_a\), define the desired
morphism \(\rho_e\colon C_\infty[\varepsilon]\to\widehat X\).

In order to verify the moreover part, it is enough to consider the following commuting square
\[
\begin{tikzcd}
C_\infty[\varepsilon]
    \arrow[r,"q_\varepsilon"]
    \arrow[d,"\widetilde\tau_\varepsilon^2"']
&
C_\infty[0]
    \arrow[d,"\widetilde\tau_0"]
\\
C_\infty[\varepsilon]
    \arrow[r,"q_\varepsilon"]
&
C_\infty[0].
\end{tikzcd}
\]
where $\tilde{\tau}_0$ is the shift $j\mapsto j+1$. It follows that if an $\varepsilon$-edge $e\in \widehat{X}_\varepsilon=\mathrm{Hom}(C_\infty[\varepsilon],X)$ is in the image of $q_\varepsilon^*$, then so is its rotation $t^2_\varepsilon(e)$.
\end{proof}

\begin{proposition}\label{prop:trivKanFibHelp}
   Let $A\hookrightarrow B$ be an embedding in $\eSSet$. If
\(
A_\varepsilon\longrightarrow B_\varepsilon
\)
is a bijection, then
\[
B\coprod_A P^*P_!A\rightarrow P^*P_!B
\]
belongs to $\overline{\eHorn}$.
In particular, if $B$ has no $\varepsilon$-witnesses, then
\(
\eta_B\colon B\rightarrow P^*P_!B
\)
belongs to $\overline{\eHorn}$.
\end{proposition}
Before proving the general situation, let us examine the case $A=\partial\Delta^1$, $B=\Delta^1$. Then the morphism in question is essentially $\Delta^{\{0,1\}}\cup\partial C_\infty[1]\hookrightarrow C_\infty[1]$ whose domain is depicted in the following diagram:
\begin{equation*}
\cdots\begin{tikzcd}[row sep=2.2em, column sep=1.6em]
&{-3}\arrow[rr,LA=1.5pt]&&{-1}\arrow[rr,LA=1.5pt] &&{1}\arrow[rr,LA=1.5pt]&&{3}\arrow[rr,LA=1.5pt] &&{5}\\
{-4}\arrow[rr,LA=1.5pt]&&{-2}\arrow[rr,LA=1.5pt]&&{0}\arrow[ru]\arrow[rr,LA=1.5pt]&&{2}\arrow[rr,LA=1.5pt]&&{4}&
\end{tikzcd}\cdots
\end{equation*}
One can see that the whole $C_\infty[1]$ is obtained by successively
attaching $2$-simplices along $\varepsilon$-horns isomorphic to
$\ELambda^2_0$ on one side and to $\ELambda^2_2$ on the other.
\begin{proof}
We obtain the particular case from the general one by setting $m\colon\emptyset\rightarrow B$.

    Denote by $\mathcal{M}$ the class of morphisms $m\colon A\rightarrow B$, having the property in question, that is, each lifting problem~\eqref{diag:liftProb01} with $p\in{\eHorn}^\boxslash$ has a solution $l$:
    \begin{equation}\label{diag:liftProb01}
        \begin{tikzcd}[column sep=3.6em, row sep=2.2em]
            {B\coprod_{A}P^*P_!(A)}\arrow[r,"{(f_1,f_2)}"]\arrow[d,"{(\eta_B,P^*P_!(m))}"']&{X}\arrow[d,"p"]\\
       {P^*P_!(B)}\arrow[ru,"l"' description,dotted]\arrow[r,"g"]&Y
        \end{tikzcd}
    \end{equation}
Applying Lemma~\ref{lem:eqLifProb}, the lifting problem~\eqref{diag:liftProb01} is equivalent to the right-hand lifting problem in~\eqref{diag:eqLifProb}.
It follows that $\mathcal{M}$ consists of precisely the morphisms having the LLP with respect to the morphisms having the form of the right-hand leg in~\eqref{diag:eqLifProb}. Hence $\mathcal{M}$ is weakly saturated.

It is enough to show that $\mathcal{M}$ contains $\partial\Delta^n\hookrightarrow\Delta^n$ for all $n\geq 0$.
To this end, note that we can represent the embedding $$\Delta^n\coprod_{\partial\Delta^n}P^*P_!(\partial\Delta^n)\rightarrow P^*P_!(\Delta^n)$$ as follows:
\begin{equation}\label{eq:fillThat}
    \Delta^{\{0,\ldots,n\}}\cup\partial C_\infty[n]\hookrightarrow C_\infty[n].
\end{equation}
Denote the inclusion~\eqref{eq:fillThat} by $A\subseteq B$. Similarly to the situation described before the proof, it is easy to show that $A\subset B$ is a composition of pushouts of $\varepsilon$-horns.
Indeed, in the first step we attach $\Sim^{\{0,\ldots,n+1\}}$ along $A\cap \Sim^{\{0,\ldots,n+1\}}\cong \ELambda^{n+1}_{0}$.
Then we proceed to $\Sim^{\{1,\ldots,n+2\}}$, and so on.  Analogously we can fill all the simplices $\Sim^{\{k,\ldots,k+n+1\}}$ with $k<0$ on the opposite side.
\end{proof}
The following characterization of $\varepsilon$-fibrant simplicial sets
will be useful in the proof of the stability result.
\begin{proposition}\label{prop:counitIsTriv}
Let $X$ be an $\varepsilon$-simplicial set. Then
\begin{enumerate}
    \item[(i)] $X$ admits fillers of $\varepsilon$-horns iff the counit $\epsilon'_X\colon P^*P_*(X)\rightarrow X$ is a trivial Kan fibration.
    \item[(ii)] Assume that $X$ admits fillers of $\varepsilon$-horns.
Then $X$ satisfies the $\varepsilon$-coherence conditions if and only if
\(
(\Delta^1\hookrightarrow\Sim^1)\boxslash
(\widehat X\rightarrow X).
\)
\end{enumerate}
Consequently, $X$ is $\varepsilon$-fibrant iff $\widehat{\epsilon}_X\colon \widehat{X}\rightarrow X$ has the right lifting property with respect to all monomorphisms in $\eSSet$.
\end{proposition}
\begin{proof}
    First focus on statement (i).
    $\Rightarrow$: Assume first that $X$ admits fillers of all $\varepsilon$-horns. We first show $\epsilon'_X$ is surjective on vertices.
    Let $x$ be a vertex of $X$. According to formula~\eqref{eq:counitDesc},
     vertices in the preimage of $x$ are in bijection with infinite chains of
     $\varepsilon$-edges
        $\hat{x}\colon C_{\infty}[0]\rightarrow X$, such that $\hat{x}$
        sends $\{0\}$ to $x$.
The existence of $\hat{x}$ follows from the one-dimensional
$\varepsilon$-horn filling conditions, which ensure that
each vertex is both the source and the target of an
$\varepsilon$-edge.

Next assume for $n\geq 1$ a lifting problem
    \begin{equation}\label{eq:eqLifProb1}
        \begin{tikzcd}
            \partial\Delta^n\arrow[d,hook]\arrow[r,"f"]&P^*P_*(X)\arrow[d,"\epsilon'_X"]\\
            \Delta^n\arrow[r,"g"]\arrow[ru,dotted,"l"'description] & X
        \end{tikzcd}
    \end{equation}
Applying Lemma~\ref{lem:eqLifProb} with $Y=\Sim^0$ (the terminal object), we obtain an equivalent lifting problem:
\begin{equation}
\begin{tikzcd}
P^*P_!(\partial\Delta^n) \coprod_{\partial\Delta^n} \Delta^n
  \arrow[rr,"{(\bar{f},g)}"]
  \arrow[d, hook]
&&
X
\\
P^*P_!(\Delta^n)
  \arrow[rru, "\bar{l}"'description, dotted]
&&
\end{tikzcd}
\end{equation}
By Proposition~\ref{prop:trivKanFibHelp}, the left leg belongs to $\overline{\eHorn}$, and so by the assumption on $X$, there is the desired lift $\bar{l}$.
We have proved that $\epsilon'_X$ is a trivial Kan fibration.

Next consider the opposite implication. Assume $\epsilon'_X$ is a trivial Kan fibration.
We will prove by induction over $n$ that $X$ admits fillers of $\varepsilon$-horns $\ELambda^{n}_i\hookrightarrow \Sim^n$, $i=0,n$.
Given any vertex $x\in X$, surjectivity of $\epsilon'_X$ ensures that $x$ belongs to an infinite path of $\varepsilon$-edges $C_\infty[0]\rightarrow X$. In particular, $x$ is the source of an $\varepsilon$-edge as well as a target of an $\varepsilon$-edge. This establishes the induction base $n=1$.
Assume $X$ admits fillers of all $\varepsilon$-horns $\ELambda^k_{i}$, $i=0,k$, with $k\leq n$. The assumption $(\partial\Delta^n\hookrightarrow\Delta^n)\boxslash\epsilon'_X$ guarantees that $X$ has the right lifting property with respect to the lower horizontal morphism in the diagram:
\begin{equation}\label{diag:leftLeg}
    \begin{tikzcd}
        \ELambda^{n+1}_0\arrow[d,hook]\arrow[r,hook]&\Sim^{n+1}\arrow[d,hook]\\
        \partial C_\infty[n]\cup \Delta^{\{0,\ldots,n\}}\arrow[r,hook]& C_\infty[n]
    \end{tikzcd}
\end{equation}
Hence, it is enough to show that $X$ has the right lifting property with respect to the left leg. Then, using Lemma~\ref{lem:fibHelp}, we can conclude that $X$ has the RLP with respect to $\ELambda^{n+1}_0\hookrightarrow \Sim^{n+1}$.

Recall the notation for rotations from Section~\ref{sec:3}, in
particular~\eqref{eq:rotHat} and~\eqref{eq:rotBeta}. Observe that
the $\varepsilon$-horn can be written as
\[
\ELambda^{n+1}_0
=
\rot{(\partial\Delta^n)}\cup\Delta^n.
\]
Consequently, the left vertical map in~\eqref{diag:leftLeg} is a
pushout of
\[
\rot{(\partial\Delta^n)}
\hookrightarrow
\partial C_\infty[n].
\]

By Proposition~\ref{prop:cohRotation}, the inclusion
\[
\partial\Delta^n\hookrightarrow\partial C_\infty[n]
\]
admits a cellular decomposition obtained by successively attaching
copies of $\rot{(\partial\Delta^n)}$ along
$\partial\Delta^n$ or $\beta(\partial\Delta^n)\subset \rot{(\partial\Delta^n)}$, and every attaching
map belongs to
\begin{equation}\label{eq:satOfEhornsLeqN}
    \overline{\bigl\{
    \ELambda^r_i\hookrightarrow\Sim^r
    \mid 1\leq r\leq n,\ i\in\{0,r\}
\bigr\}}.
\end{equation}
Indeed, in the proof of Proposition~\ref{prop:cohRotation}(i),
each attached simplex is obtained by adjoining one vertex
to a face of $\partial\Delta^n$, and therefore has dimension
at most $n$.
Moreover, the first stage of this
decomposition is
\[
\partial\Delta^n
\hookrightarrow
\rot{(\partial\Delta^n)}.
\]
Therefore the remaining inclusion
\[
\rot{(\partial\Delta^n)}
\hookrightarrow
\partial C_\infty[n]
\]
belongs to~\eqref{eq:satOfEhornsLeqN}.
Thus the left vertical map in~\eqref{diag:leftLeg} also belongs to~\eqref{eq:satOfEhornsLeqN} and consequently has the left lifting property with respect to $X$ by the induction hypothesis.
The same construction with the reversed arrows gives fillers
for $\ELambda^{n+1}_{n+1}\hookrightarrow\Sim^{n+1}$.
This completes the induction.

Next focus on statement (ii). By the definition of $\widehat X$, a morphism
$\Delta^1\to\widehat X$ corresponds to a morphism
$C_\infty[1]\to X$, while a morphism
$\Sim^1\to\widehat X$ corresponds to a morphism
$C_\infty[\varepsilon]\to X$. Hence the following two lifting
problems are equivalent:
    \begin{equation}\label{eq:eqLifProb1x}
        \begin{tikzcd}
            \Delta^1\arrow[d,hook]\arrow[r,"f"]&\widehat{X}\arrow[d,"\hat{\epsilon}_X"]\\
            \Sim^1\arrow[r,"g"]\arrow[ru,dotted,"l"'description] & X
        \end{tikzcd}
        \quad\quad
        \begin{tikzcd}
C_\infty[1] \coprod_{\Delta^1} \Sim^1
  \arrow[rr,"{(\bar{f},g)}"]
  \arrow[d, hook]
&&
X
\\
C_\infty[\varepsilon]
  \arrow[rru, "\bar{l}"'description, dotted]
&&
{}
\end{tikzcd}
    \end{equation}
The left leg is easy to picture:
\begin{equation}
\begin{gathered}
\cdots
\begin{tikzcd}[row sep=2.2em, column sep=1.6em]
&-2\arrow[rr,LA=1.5pt]\arrow[rd]
&&0\arrow[rr,LA=1.5pt]\arrow[rd,LA=1.5pt]
&&2\arrow[rr,LA=1.5pt]\arrow[rd]
&&4\arrow[rr,LA=1.5pt]\arrow[rd]
&&6\\
{-3}\arrow[rr,LA=1.5pt]\arrow[ru]
&&{-1}\arrow[rr,LA=1.5pt]\arrow[ru]
&&{1}\arrow[rr,LA=1.5pt]\arrow[ru]
&&{3}\arrow[rr,LA=1.5pt]\arrow[ru]
&&{5}\arrow[ru]&
\end{tikzcd}
\cdots
\\[1em]
\downarrow
\\[1em]
\cdots
\begin{tikzcd}[row sep=2.2em, column sep=1.6em]
&-2\arrow[rr,LA=1.5pt]\arrow[rd,LA=1.5pt]
&&0\arrow[rr,LA=1.5pt]\arrow[rd,LA=1.5pt]
&&2\arrow[rr,LA=1.5pt]\arrow[rd,LA=1.5pt]
&&4\arrow[rr,LA=1.5pt]\arrow[rd,LA=1.5pt]
&&6\\
{-3}\arrow[rr,LA=1.5pt]\arrow[ru]
&&{-1}\arrow[rr,LA=1.5pt]\arrow[ru]
&&{1}\arrow[rr,LA=1.5pt]\arrow[ru]
&&{3}\arrow[rr,LA=1.5pt]\arrow[ru]
&&{5}\arrow[ru]&
\end{tikzcd}
\cdots
\end{gathered}
\end{equation}
Denote the inclusion by $i\colon C_\infty[1]\coprod_{\Delta^1}\Sim^1\hookrightarrow C_\infty[\varepsilon]$.
Obviously $i$ belongs to the weakly saturated class generated by \eqref{L1} and \eqref{L2}.
Conversely, let $j\colon A\hookrightarrow B$ be either~\eqref{L1} or~\eqref{L2}.
It is easy to see that there is an inclusion $f\colon A\hookrightarrow C_\infty[1]\coprod_{\Delta^1}\Sim^1$ belonging to $\overline{\eHorn}$ and its extension $g\colon B\hookrightarrow C_\infty[\varepsilon]$, fitting into the square:
\[
\begin{tikzcd}
A\arrow[d,hook,"j"]\arrow[r,hook,"f"]
&
C_\infty[1]\coprod_{\Delta^1}\Sim^1
\arrow[d,hook,"i"]
\\
B\arrow[r,hook,"g"]
&
C_\infty[\varepsilon]
\end{tikzcd}
\]
Hence, if $X$ admits fillers of $\varepsilon$-horns and $i\boxslash X$ then $j\boxslash X$ as well. Indeed, given a morphism $A\to X$, first extend it along
$f$ using the $\varepsilon$-horn filling property. Then extend the
resulting map along $i$, and finally restrict the extension along
$g\colon B\hookrightarrow C_\infty[\varepsilon]$.
\end{proof}
Let $X$ be an $\varepsilon$-simplicial set.
For an $n$-simplex $\sigma\colon\Delta^n\to X$, define a simplicial set
\[
    \mathrm{Rot}_X(\sigma)
    =
    \{\sigma\}
    \times_{\mathrm{Fun}(\Delta^n,X)}
    \mathrm{Fun}(\Delta^n,P^*P_*X),
\]
where the map to $\mathrm{Fun}(\Delta^n,X)$ is induced
by the counit $\epsilon'_X$.
Its $k$-simplices are maps
$h\colon\Delta^k\times\Delta^n\to P^*P_*X$
such that $\epsilon'_X h=\sigma\operatorname{pr}_2$.
By \eqref{eq:counitDesc}, vertices of $\mathrm{Rot}_X(\sigma)$ correspond to
extensions of $\sigma$ to $C_\infty[n]$.
Thus $\mathrm{Rot}_X(\sigma)$ is a simplicial set
of coherent choices of rotations of $\sigma$.

\begin{corollary}\label{cor:contractibleRotations}
For every simplex $\sigma$ of an $\varepsilon$-simplicial set
$X$ admitting fillers for $\varepsilon$-horns, the simplicial set $\mathrm{Rot}_X(\sigma)$
is a contractible Kan complex.
\end{corollary}

\begin{proof}
The counit $\epsilon'_X$ is a trivial Kan fibration,
so it has the right lifting property with respect to
every monomorphism, in particular
\[
    \partial\Delta^k\times\Delta^n
    \hookrightarrow
    \Delta^k\times\Delta^n
    \qquad(k\geq 0).
\]
By the internal-hom adjunction, the induced map
\[
    \mathrm{Fun}(\Delta^n,P^*P_*X)
    \longrightarrow
    \mathrm{Fun}(\Delta^n,X)
\]
is a trivial Kan fibration.
Its fiber over $\sigma$ is therefore a contractible
Kan complex.
\end{proof}
\begin{corollary}\label{cor:retract}
    Every $\varepsilon$-fibrant $X\in\eSSet$ is a retract of $\widehat{X}$.
\end{corollary}
\begin{proof}
    By Proposition~\ref{prop:counitIsTriv} there is a section $s\colon X\rightarrow \widehat{X}$
    \begin{equation*}
        \begin{tikzcd}
            & {\widehat{X}}\arrow[d,"\hat{\epsilon}_X"]\\
            X\arrow[ru,"s",dashed]\arrow[r,"{\mathrm{id}}"]&X
        \end{tikzcd}
    \end{equation*}
\end{proof}

The previous corollary gives us a bridge through which we can transfer lifting problems of an $\varepsilon$-fibrant $X$ to those of $\widehat{X}$.
The essential advantage of $\widehat{X}$ is that it admits a convenient cellular description.
\begin{lemma}\label{lem:cellularStructureOfHatX}
    Let \(X\) be an arbitrary \(\varepsilon\)-simplicial set. Then
    \(\emptyset\hookrightarrow \widehat{X}\) belongs to the weakly saturated class generated by
    \begin{align}
        \{\partial C_k[n]\hookrightarrow C_k[n]
        \mid \infty\geq k\geq 1,\ n\geq 0\}
        \quad \text{and} \label{eq:genCyc}\\
        \{C_{2r}[1]\hookrightarrow C_{2r}[\varepsilon]\mid r\geq 1\}
        \cup
        \{C_\infty[1]\hookrightarrow C_\infty[\varepsilon]\}.
        \label{eq:genEpsCyc}
    \end{align}
\end{lemma}
\begin{proof}
Set \(A=P^*P_*X\). The inclusion
\(\emptyset\hookrightarrow\widehat X\) factors as
\[
    \emptyset\hookrightarrow A\hookrightarrow\widehat X.
\]
We first show that
\(\emptyset\hookrightarrow A\) belongs to the weakly saturated class
generated by~\eqref{eq:genCyc}.
This follows from Lemma~\ref{lem:genOfMonoInPar},
the isomorphisms~\eqref{eq:repOfCycles}, and the fact that
\(P^*\) preserves colimits.

It remains to show that
\[
    A\hookrightarrow\widehat X
\]
belongs to the weakly saturated class generated by~\eqref{eq:genEpsCyc}.
By Lemma~\ref{lem:canonicalOrbitsEps},
the set
\(
    \widehat X_\varepsilon\setminus A_\varepsilon
\)
is a disjoint union of \(t_\varepsilon^2\)-orbits.
Choose a well-ordering of these orbits and attach them one by one.

Let
\(
    e\in\widehat X_\varepsilon\setminus A_\varepsilon
\)
represent an orbit, and let \(a\in A_1\) be its underlying edge.
By Lemma~\ref{lem:canonicalOrbitsEps}, there is a canonical orbit morphism
\(
    \rho_e\colon C_\infty[\varepsilon]\rightarrow\widehat X
\)
whose restriction to \(C_\infty[1]\) is the canonical orbit morphism
\(
    \rho_a\colon C_\infty[1]\rightarrow A.
\)

Suppose first that the \(t_\varepsilon^2\)-orbit of \(e\) is infinite.
Then all its elements are distinct, and the whole orbit is attached
by pushing out
\(
    \rho_a\colon C_\infty[1]\rightarrow A
\)
along
\[
    C_\infty[1]\hookrightarrow C_\infty[\varepsilon].
\]

Now suppose that the orbit of \(e\) has exactly \(r\) elements. Then
\(
    (t_\varepsilon^2)^r(e)=e.
\)
Consequently,
\(
    t_1^{2r}(a)=a,
\)
so the morphism \(\rho_a\) factors through
\[
    C_\infty[1]\longrightarrow C_{2r}[1].
\]
Moreover, the periodicity
\((t_\varepsilon^2)^r(e)=e\) implies that \(\rho_e\) is invariant
under the shift by \(2r\), and hence factors through
\(
    C_\infty[\varepsilon]\longrightarrow C_{2r}[\varepsilon].
\)
Thus the whole orbit of \(e\) is attached by a pushout along
\[
    C_{2r}[1]\hookrightarrow C_{2r}[\varepsilon].
\]

The missing \(\varepsilon\)-witnesses form disjoint orbits, which may be attached one by one.
Hence \(A\hookrightarrow\widehat X\) is a transfinite composition of pushouts of morphisms in~\eqref{eq:genEpsCyc}.
Thus \(\emptyset\hookrightarrow\widehat X\) belongs to the weakly
saturated class generated by~\eqref{eq:genCyc} and~\eqref{eq:genEpsCyc}.
\end{proof}

\section{Stability theorems}\label{sec:stability}

In this section, we show that if $X$ and $Y$ are
$\varepsilon$-fibrant, then the $\varepsilon$-simplicial set
$\mathrm{Fun}(X,Y)$ is almost $\varepsilon$-fibrant. The main
difficulty is to construct fillers for $\varepsilon$-horns. This
horn-filling condition is expressed by the following two equivalent
lifting problems (for $n\geq 1$, $l=0,n$):
\begin{equation}
\begin{tikzcd}
{\ELambda^{n}_l}\arrow[r,"g"]\arrow[d,hook,"i"]
    &{\mathrm{Fun}(X,Y)}\\
{\Sim^n}\arrow[ru,dotted,"f"']
\end{tikzcd}
\qquad\Longleftrightarrow\qquad
\begin{tikzcd}
{\ELambda^{n}_l\times X}\arrow[d,hook]\arrow[r,"\overline{g}"]
    &Y\\
{\Sim^n\times X}\arrow[ru,dotted,"\overline{f}"']
\end{tikzcd}
\end{equation}
By the results of the previous section, $X$ is a retract of an
$\varepsilon$-simplicial set $\widehat{X}$ admitting a convenient
cellular decomposition
(Lemma~\ref{lem:cellularStructureOfHatX}). The construction of the
required fillers proceeds in several steps. The following
proposition provides the main technical input for handling the
subobject $P^*P_*X\subseteq\widehat{X}$. Its proof is the most
involved simplicial-combinatorial argument in the paper. In the statement we use the notation $\overline{\mathbb N}:=\mathbb{N}\cup\{\infty\}$.
\begin{proposition}\label{prop:key}
    For each $n\geq 1$, $k\geq 0$, $l=0,n$, and $r\in\overline{\mathbb{N}}$, we have
    \begin{equation}\label{eq:key}
        (\ELambda^n_l\hookrightarrow\Sim^n)\Box (\partial C_r[k]\hookrightarrow C_r[k])\in\overline{\eHorn}.
    \end{equation}
\end{proposition}
Let us firstly explore the easiest case $n=1$, $l=0$, $k=0$, and $r=\infty$:
$$(\ELambda^1_0\hookrightarrow\Sim^1)\Box(\partial C_\infty[0]\hookrightarrow C_\infty[0])$$ Here the domain is
$A_0:=\{0\}\times C_\infty[0]\coprod_{\{0\}\times\emptyset}\Sim^1\times\emptyset=\{0\}\times C_\infty[0]$
and it corresponds to the bottom horizontal line in~\eqref{diag:easiesEx}.
The codomain is shown below:
\begin{equation}\label{diag:easiesEx}
\begin{tikzcd}[
  column sep=1.4em,
  row sep=1.4em
]
\cdots \arrow[r] &
(1,-2) \arrow[r] &
(1,-1) \arrow[r] &
(1,0) \arrow[r] &
(1,1) \arrow[r] &
(1,2) \arrow[r] &
\cdots
\\
\cdots \arrow[r] \arrow[ur,LA=1.5pt]&
(0,-2) \arrow[u] \arrow[ur,LA=1.5pt] \arrow[r] &
(0,-1) \arrow[u] \arrow[ur,LA=1.5pt] \arrow[r] &
(0,0) \arrow[u] \arrow[ur,LA=1.5pt] \arrow[r] &
(0,1) \arrow[u] \arrow[ur,LA=1.5pt] \arrow[r] &
(0,2) \arrow[u] \arrow[r]\arrow[ur,LA=1.5pt] &
\cdots
\end{tikzcd}
\end{equation}
Hence the codomain is $2$-dimensional, since $n+k+1=2$, and
we can divide its nondegenerate $2$-simplices into two sets
$U_1$ and $U_2$ as follows. The simplices in $U_1$ have exactly
one vertex on the upper horizontal line, whereas those in $U_2$
have two such vertices.

Using $\varepsilon$-horn attachments, we can first attach all
the simplices in $U_1$ and then all the simplices in $U_2$.
The key observation is that, for every simplex
$\sigma\in U_1$, the face $d_2\sigma$ lies in $A_0$, whereas,
for every $\sigma\in U_2$, the face $d_2\sigma$ lies in
\(
A_0\cup\bigcup_{\tau\in U_1}\operatorname{Im}(\tau).
\)
The proof of the general case follows the same strategy.

\begin{proof}

By duality, it suffices to consider $l=0$. Denote the domain
of~\eqref{eq:key} by
\[
A_0:=\ELambda^n_0\times C_r[k]
\coprod_{\ELambda^n_0\times\partial C_r[k]}
\Sim^n\times\partial C_r[k],
\]
and put $B:=\Sim^n\times C_r[k]$.

Recall that a $d$-simplex $\sigma$ of $B$ is a pair
$(\sigma_1,\sigma_2)$, where
$\sigma_1=\langle i_0,\ldots,i_d\rangle$ with
$0\leq i_0\leq\cdots\leq i_d\leq n$, and
$\sigma_2=[\langle j_0,\ldots,j_d\rangle]$ with
$j_0\leq\cdots\leq j_d\leq j_0+k+1$. We denote it by
\begin{equation}\label{eq:codeSigma}
\sigma=
\left[
\left\langle
\begin{matrix}
i_0 & i_1 & \cdots & i_d \\
j_0 & j_1 & \cdots & j_d
\end{matrix}
\right\rangle
\right].
\end{equation}

Set $m:=n+k+1$. Then $B$ is the union of the images of its
nondegenerate $m$-simplices $\sigma\colon\Sim^m\to B$. A
simplex of the form~\eqref{eq:codeSigma} is a nondegenerate
$m$-simplex precisely when $d=m$ and
\[
(i_s-i_{s-1},j_s-j_{s-1})\in\{(1,0),(0,1)\}
\qquad (1\leq s\leq m).
\]
In particular, $i_0=0$, $i_m=n$, and $j_m=j_0+k+1$.

Partition these nondegenerate $m$-simplices into $k+2$
disjoint sets $U_1,\ldots,U_{k+2}$ by declaring that
$\sigma\in U_p$ if exactly its last $p$ indices
$i_{m-p+1},\ldots,i_m$ are equal to $n$. Their number can
therefore vary from $1$ to $k+2$. Recursively define
\[
A_p:=A_{p-1}\cup\bigcup_{\sigma\in U_p}\operatorname{Im}(\sigma),
\qquad 1\leq p\leq k+2.
\]
Thus $A_{k+2}=B$.

A simplex $\sigma\colon\Sim^d\to B$ of the
form~\eqref{eq:codeSigma} lies in $A_0$ if and only if either
some $i\in\{1,\ldots,n\}$ does not occur among
$i_0,\ldots,i_d$, or some integer $j$ with
$j_0<j<j_0+k+1$ does not occur among $j_0,\ldots,j_d$.

\emph{Claim 1.}
Let $\sigma\colon\Sim^m\to B$ be a nondegenerate simplex with
representative~\eqref{eq:codeSigma}. Let
$\pi_2\colon\Sim^{k+1}\to C_r[k]$ be represented by
$\langle j_0,j_0+1,\ldots,j_0+k+1\rangle$ and put
$\pi:=(\mathrm{id},\pi_2)$. Then $\sigma$ admits a unique
factorization
\[
\begin{tikzcd}[ampersand replacement=\&]
\Sim^m\arrow[r,hook,"\hat{\sigma}"]
\&\Sim^n\times\Sim^{k+1}\arrow[r,"\pi"]
\&\Sim^n\times C_r[k].
\end{tikzcd}
\]

\emph{Proof of Claim 1.}
The required lift is represented by
\[
\hat{\sigma}=
\left\langle
\begin{matrix}
i_0&i_1&\cdots&i_m\\
j_0-j_0&j_1-j_0&\cdots&j_m-j_0
\end{matrix}
\right\rangle.
\]
Indeed, $\pi\hat{\sigma}=\sigma$, and the description of
nondegenerate $m$-simplices above shows that $\hat{\sigma}$ is a
monomorphism. Any other lift of $\sigma$ is again a nondegenerate
$m$-simplex and hence lies outside
$\Sim^n\times\ELambda^{k+1}_0$, whose dimension is $m-1$.
Uniqueness therefore follows from
Lemma~\ref{lem:PulbackAdPushout}(iii). \qed

For $A\subseteq B$ and a nondegenerate
$\sigma\colon\Sim^m\to B$, consider the pullback diagram
\begin{equation}\label{diag:boss}
\begin{tikzcd}
\sigma^*A \arrow[r,hook]\arrow[d,hook]
    & D_A^\sigma \arrow[d,hook]\arrow[r]
    & A\arrow[d,hook]\\
\Sim^m\arrow[r,hook,"\hat{\sigma}"]
    & \Sim^n\times\Sim^{k+1}\arrow[r,"\pi"]
    & \Sim^n\times C_r[k].
\end{tikzcd}
\end{equation}

We call $\sigma^*A\subseteq\Sim^m$ \emph{admissible} if
it is a union of simplices, each either contained in
$\partial_m\Sim^m$ or containing the unique
$\varepsilon$-edge of $\Sim^m$. Part~(ii) of
Proposition~\ref{prop:cohRotation} implies that
\begin{equation}\label{eq:admissible}
\sigma^*A\text{ admissible}
\quad\Longrightarrow\quad
(\sigma^*A\cup\partial_m\Sim^m\hookrightarrow\Sim^m)
\in\overline{\eHorn}.
\end{equation}
Admissibility is preserved under unions, and $\sigma^*$ preserves
unions. Hence, for every family $(C_i)_{i\in I}$ of subobjects of
$B$,
\begin{equation}\label{eq:unionIsAdm}
\bigl(\forall i\in I,\ \sigma^*C_i\text{ admissible}\bigr)
\Longrightarrow
\sigma^*\Bigl(\bigcup_{i\in I}C_i\Bigr)
=\bigcup_{i\in I}\sigma^*C_i
\text{ admissible}.
\end{equation}

\emph{Claim 2.} Fix $1\leq p\leq k+2$,
$\sigma\in U_p$, and $U\subseteq U_p\setminus\{\sigma\}$, and set
\[
A':=A_{p-1}\cup\bigcup_{\tau\in U}\operatorname{Im}(\tau).
\]
Then:
\begin{enumerate}
\item[(i)] $\partial_m\Sim^m\subseteq\sigma^*A_{p-1}$.

\item[(ii)] The following two squares are pushouts:
\begin{equation}\label{diag:pushOut}
\begin{tikzcd}
\sigma^*A' \arrow[r]\arrow[d,hook]
    & A'\cap\operatorname{Im}(\sigma)
      \arrow[d,hook]\arrow[r,hook]
    & A'\arrow[d,hook]\\
\Sim^m\arrow[r,"\sigma"]
    & \operatorname{Im}(\sigma)\arrow[r,hook]
    & A'\cup\operatorname{Im}(\sigma).
\end{tikzcd}
\end{equation}

\item[(iii)] $\sigma^*A'\subseteq\Sim^m$ is admissible.
\end{enumerate}

After proving Claim~2, we are done. Indeed, well-order each $U_p$ and attach its simplices
successively. By Claim~2(ii), the inclusion
$A'\hookrightarrow A'\cup\operatorname{Im}(\sigma)$ is a pushout
of $\sigma^*A'\hookrightarrow\Sim^m$. By Claim~2(i),
$\partial_m\Sim^m\subseteq\sigma^*A'$, so Claim~2(iii)
and~\eqref{eq:admissible} imply that
$\sigma^*A'\hookrightarrow\Sim^m$ belongs to
$\overline{\eHorn}$.
Thus every inclusion $A_{p-1}\hookrightarrow A_p$ belongs to
$\overline{\eHorn}$, and hence so does the composite
$A_0\hookrightarrow B$. It remains to prove Claim~2(i)--(iii).

Fix $\sigma\in U_p$. After translating the $C_r[k]$-coordinate, we may assume
the following representatives, with $0=j_0\leq\cdots\leq j_m=k+1$:
\begin{equation}\label{eq:matRep}
    \hat{\sigma}=
\left\langle
\begin{matrix}
i_0 & i_1 & \cdots & i_{m} \\
j_0 & j_1 & \cdots & j_{m}
\end{matrix}
\right\rangle
\quad\text{and}\quad
 \sigma=\left[
\left\langle
\begin{matrix}
i_0 & i_1 & \cdots & i_{m} \\
j_0 & j_1 & \cdots & j_{m}
\end{matrix}
\right\rangle
\right]
\end{equation}

\emph{Proof of Claim 2(i).}
If $p=1$, then $i_{m-1}<n$ by the definition of $U_1$.
Hence
\begin{equation}
    d_m\sigma=\left[
\left\langle
\begin{matrix}
i_0 & i_1 & \cdots & i_{m-1} \\
j_0 & j_1 & \cdots & j_{m-1}
\end{matrix}
\right\rangle
\right]\in\partial_n\Sim^n\times C_r[k]
\subseteq A_0=A_{p-1}
\end{equation}
and therefore
$\partial_m\Sim^m\subseteq\sigma^*A_{p-1}$.

Now suppose that $p\geq2$. Then $i_m=i_{m-1}$, so
$j_m=j_{m-1}+1$ and hence $j_{m-1}=k$. The simplex
\[
\tau:=
\left[
\left\langle
\begin{matrix}
i_0 & i_0 & i_1 & \cdots & i_{m-1}\\
-1  & 0   & j_1 & \cdots & j_{m-1}
\end{matrix}
\right\rangle
\right]
\]
is a nondegenerate $m$-simplex in $U_{p-1}$ satisfying
$d_0\tau=d_m\sigma$. Consequently,
\[
\sigma(\partial_m\Sim^m)
=\operatorname{Im}(d_m\sigma)
\subseteq\operatorname{Im}(\tau)
\subseteq A_{p-1},
\]
which proves
$\partial_m\Sim^m\subseteq\sigma^*A_{p-1}$. \qed

\emph{Proof of Claim 2(ii).}
Set
\[
E:=\bigcup_{1\leq j\leq k}\partial_j\Sim^{k+1}
\cup\Delta^{\{1,\ldots,k\}}.
\]
By Lemma~\ref{lem:PulbackAdPushout}(i) and
$A_0\subseteq A'$, we have
\begin{equation}\label{eq:Dsigma}
D^\sigma_{A_0}
=
\bigl(\ELambda^n_0\times\Sim^{k+1}\bigr)
\cup
\bigl(\Sim^n\times E\bigr)
\subseteq D^\sigma_{A'}.
\end{equation}
Moreover, Claim~2(i) gives
\begin{equation}\label{eq:intersectWithBottom}
\bigl(\Sim^n\times\partial_{k+1}\Sim^{k+1}\bigr)
\cap\operatorname{Im}(\hat{\sigma})
\subseteq
\partial_m\operatorname{Im}(\hat{\sigma})
\subseteq D^\sigma_{A'}.
\end{equation}
Since
$E\cup\partial_{k+1}\Sim^{k+1}=\ELambda^{k+1}_0$,
\eqref{eq:Dsigma} and~\eqref{eq:intersectWithBottom} imply
\[
\operatorname{Im}(\hat{\sigma})
\cap\bigl(\Sim^n\times\ELambda^{k+1}_0\bigr)
\subseteq D^\sigma_{A'}.
\]
Lemma~\ref{lem:PulbackAdPushout}(iii) states that $\pi$ is
levelwise injective outside
$\Sim^n\times\ELambda^{k+1}_0$. Hence $\sigma$ is levelwise
injective on $\Sim^m\setminus\sigma^*A'$. Therefore the
left-hand square in~\eqref{diag:pushOut} is a pushout, while
the right-hand square is the canonical pushout square for a
union of subobjects. This proves Claim~2(ii). \qed

\emph{Proof of Claim 2(iii).}
Under our identification,
\[
\sigma^*A'
=\operatorname{Im}(\hat{\sigma})\cap D^\sigma_{A'}.
\]
We first show that
$\sigma^*A_0=\operatorname{Im}(\hat{\sigma})\cap
D^\sigma_{A_0}$ is admissible. Using~\eqref{eq:Dsigma},
we consider the following intersections:
\begin{itemize}
\item For $1\leq q\leq n$, the intersection with
$\partial_q\Sim^n\times\Sim^{k+1}$ is the face obtained
from~\eqref{eq:matRep} by deleting all columns with $i_*=q$.
For $q=n$ it is contained in
$\partial_m\operatorname{Im}(\hat{\sigma})$; for $q<n$ it
contains the unique $\varepsilon$-edge.

\item For $1\leq s\leq k$, the intersection with
$\Sim^n\times\partial_s\Sim^{k+1}$ is obtained by deleting
all columns with $j_*=s$. The resulting simplex contains the
$\varepsilon$-edge.

\item The intersection with
$\Sim^n\times\Delta^{\{1,\ldots,k\}}$ is contained in
$\partial_m\operatorname{Im}(\hat{\sigma})$.
\end{itemize}
Thus $\sigma^*A_0$ is admissible. Put
\[
S:=\sigma^*A_0\cup\partial_m\Sim^m;
\]
then $S$ is admissible as well.

Set $V:=U_1\cup\cdots\cup U_{p-1}\cup U$, and let
$U'\subseteq V$ consist of the simplices factoring through the
fixed map $\pi$. Thus, for $\tau\in U'$,
\begin{equation}\label{eq:repOfTau}
\tau=\pi\hat{\tau},
\qquad
\hat{\tau}=
\left\langle
\begin{matrix}
i'_0&i'_1&\cdots&i'_m\\
j'_0&j'_1&\cdots&j'_m
\end{matrix}
\right\rangle,
\qquad
j'_0=0,\quad j'_m=k+1.
\end{equation}
By Lemma~\ref{lem:PulbackAdPushout}(iii),
\[
\pi^*\operatorname{Im}(\tau)
\subseteq
\Sim^n\times\ELambda^{k+1}_0
\cup\operatorname{Im}(\hat{\tau}).
\]
Since
$\ELambda^{k+1}_0
=E\cup\partial_{k+1}\Sim^{k+1}$,
equations~\eqref{eq:Dsigma}
and~\eqref{eq:intersectWithBottom} give
\begin{equation}\label{eq:pullbackOfTau}
\sigma^*\operatorname{Im}(\tau)
\subseteq
S\cup
\bigl(\operatorname{Im}(\hat{\sigma})
\cap\operatorname{Im}(\hat{\tau})\bigr).
\end{equation}
The last intersection is a common face containing the
$\varepsilon$-edge and is therefore admissible.

\emph{Claim 3.}
For every $\tau\in V$, either
$\sigma^*\operatorname{Im}(\tau)\subseteq S$, or there exists
$\rho\in U'$ such that
\[
\sigma^*\operatorname{Im}(\tau)
\subseteq S\cup\sigma^*\operatorname{Im}(\rho).
\]

Assuming Claim~3, preservation of unions under pullback
and~\eqref{eq:pullbackOfTau} yield
\begin{align}
\sigma^*A'
&\subseteq S\cup
  \bigcup_{\tau\in V}\sigma^*\operatorname{Im}(\tau)\notag \\
&\subseteq S\cup
  \bigcup_{\tau\in U'}\sigma^*\operatorname{Im}(\tau)\notag \\
&\subseteq S\cup
  \bigcup_{\tau\in U'}
  \bigl(\operatorname{Im}(\hat{\sigma})
  \cap\operatorname{Im}(\hat{\tau})\bigr).
\label{eq:pullbackOfAprime}
\end{align}
Conversely, the last term is contained in $\sigma^*A'$:
Indeed, $A_0\subseteq A'$ and Claim~2(i) give
$S\subseteq\sigma^*A'$, while
\[
\operatorname{Im}(\hat{\sigma})
\cap\operatorname{Im}(\hat{\tau})
\subseteq\sigma^*\operatorname{Im}(\tau)
\subseteq\sigma^*A'
\]
for every $\tau\in U'$. Hence all the inclusions
in~\eqref{eq:pullbackOfAprime} are equalities. Its last term is
admissible, and therefore so is $\sigma^*A'$. %

It remains to prove Claim~3.

\begin{proof}[Proof of Claim~3]
If $r=1$, then every $\tau\in V$ belongs to $U'$, and we may
take $\rho=\tau$. Hence assume $r\geq2$ and fix
$\tau=(\tau_1,\tau_2)\in V$. Again, there is nothing to prove
if $\tau\in U'$, so suppose that $\tau\notin U'$.

By Lemma~\ref{lem:PulbackAdPushout}(ii),
\[
\pi^*\operatorname{Im}(\tau)
\subseteq\Sim^n\times\partial\Sim^{k+1}.
\]
After intersecting with $\operatorname{Im}(\hat{\sigma})$, the
part lying over $E$ is contained in $\sigma^*A_0$, while the part
lying over $\partial_{k+1}\Sim^{k+1}$ is contained in
$\partial_m\operatorname{Im}(\hat{\sigma})$. Thus only the part
lying over $\partial_0\Sim^{k+1}$ but outside $D^\sigma_{A_0}$
remains. If this part is empty, the first alternative in Claim~3
holds.
If this part is nonempty, then $\tau_2$ has a representative
whose vertices contain all the integers $1,\ldots,k+1$.
Since the endpoints of such a representative differ by $k+1$,
it must have one of the following two forms:
\[
\tau=
\left[
\left\langle
\begin{matrix}
0&*&\cdots&*&n\\
0&*&\cdots&*&k+1
\end{matrix}
\right\rangle
\right]
\quad\text{or}\quad
\tau=
\left[
\left\langle
\begin{matrix}
0&*&\cdots&*&n\\
1&*&\cdots&*&k+2
\end{matrix}
\right\rangle
\right].
\]
The first possibility would imply $\tau\in U'$. We may
therefore choose a representative
\[
\tau=
\left[
\left\langle
\begin{matrix}
i'_0&i'_1&\cdots&i'_m\\
j'_0&j'_1&\cdots&j'_m
\end{matrix}
\right\rangle
\right],
\qquad j'_0=1,\quad j'_m=k+2.
\]

Let us examine simplices in $\pi^*\mathrm{Im}(\tau)$. Consider a $d$-simplex
\[
\upsilon=
\left\langle
\begin{matrix}
u_0&\cdots&u_d\\
v_0&\cdots&v_d
\end{matrix}
\right\rangle
\colon\Delta^d\longrightarrow
\Sim^n\times\partial_0\Sim^{k+1}
\]
in $\pi^*\operatorname{Im}(\tau)$. By definition, there are
indices $0\leq a_0\leq\cdots\leq a_d\leq m$ such that
\begin{equation}\label{eq:eqAndCong}
\langle u_0,\ldots,u_d\rangle
=\langle i'_{a_0},\ldots,i'_{a_d}\rangle,
\qquad
\langle v_0,\ldots,v_d\rangle
\equiv
\langle j'_{a_0},\ldots,j'_{a_d}\rangle
\pmod r.
\end{equation}
Suppose that $\upsilon$ is not contained in $D^\sigma_{A_0}$.
Then $v_0,\ldots,v_d$ contain all the integers
$1,\ldots,k+1$.

For $r=\infty$, the congruence in~\eqref{eq:eqAndCong} is
equality. If $2\leq r<\infty$, it is a simultaneous translation
by a multiple of $r$. Since $r\geq2$
and $1\leq j'_s\leq k+2$, this translation must be zero.
Hence
\[
v_s=j'_{a_s}\qquad(0\leq s\leq d).
\]

Let $t$ be the largest index such that $j'_t\leq k+1$, and set
\[
\theta:=
\left\langle
\begin{matrix}
i'_0&\cdots&i'_t\\
1=j'_0&\cdots&j'_t=k+1
\end{matrix}
\right\rangle.
\]
The preceding equality shows that
$\operatorname{Im}(\upsilon)\subseteq\operatorname{Im}(\theta)$.

If $\tau\in U_1$, then $i'_t<n$. Consequently,
\(
\operatorname{Im}(\theta)\cap
\operatorname{Im}(\hat{\sigma})
\subseteq\partial_m\operatorname{Im}(\hat{\sigma}),
\)
and therefore
$\sigma^*\operatorname{Im}(\tau)\subseteq S$.

Now suppose that $\tau\in U_q$ for some $2\leq q\leq p$.
Then $i'_m=i'_{m-1}$, so $j'_m=j'_{m-1}+1$ and consequently
$t=m-1$. Define
\[
\rho:=
\left[
\left\langle
\begin{matrix}
0&i'_0&\cdots&i'_{m-1}\\
0=j'_0-1&j'_0&\cdots&j'_{m-1}
\end{matrix}
\right\rangle
\right].
\]
This is a nondegenerate $m$-simplex in $U_{q-1}$, and its lift
along $\pi$ satisfies $\theta=d_0\hat{\rho}$. Hence
\(
\operatorname{Im}(\theta)
\subseteq\pi^*\operatorname{Im}(\rho).
\)
Moreover, $\rho\in U_{q-1}\subseteq V$, and its displayed
representative shows that $\rho\in U'$. We therefore obtain
\[
\sigma^*\operatorname{Im}(\tau)
\subseteq S\cup\sigma^*\operatorname{Im}(\rho),
\]
which proves Claim~3.
\end{proof}
This completes the proof of Claim~2(iii), and hence the proof of
the proposition.
\end{proof}

We now prove the stability result for $\varepsilon$-fibrant simplicial sets.
It falls just short of showing that $\mathrm{Fun}(X,Y)$ is $\varepsilon$-fibrant,
since the argument does not cover the $\varepsilon$-horns in dimension
$n=1$; the obstruction appears in case~(2) of the proof.
This deficiency is addressed and resolved in the next section (Corollary~\ref{cor:lastCor}).

\begin{theorem}\label{thm:stability}
    Let $X,Y$ be $\varepsilon$-fibrant. Then $\mathrm{Fun}(X,Y)$
    \begin{itemize}
        \item[(i)] admits fillers of $\varepsilon$-horns $\ELambda^{n}_0,\ELambda^{n}_n$, for $n\geq 2$;
        \item[(ii)] satisfies the $\varepsilon$-coherence conditions.
    \end{itemize}
\end{theorem}
\begin{proof}
By Corollary~\ref{cor:retract}, \(X\) is a retract of \(\widehat X\).
Hence \(\mathrm{Fun}(X,Y)\) is a retract of
\(\mathrm{Fun}(\widehat X,Y)\).
Thus it is enough to prove the assertions for
\(\mathrm{Fun}(\widehat X,Y)\).

Recall that \(\overline{\mathbb N}:=\mathbb N \cup\{\infty\}\).
By Lemma~\ref{lem:cellularStructureOfHatX}, it is enough to show that
the following four classes of morphisms belong to the weakly saturated
class generated by the \(\varepsilon\)-horns and the
\(\varepsilon\)-coherence morphisms:
\begin{enumerate}[(1)]
\item
\(\{\ELambda^n_k\hookrightarrow\Sim^n\mid n\geq2,\;k=0,n\}
\Box
\{\partial C_r[m]\hookrightarrow C_r[m]\mid
r\in\overline{\mathbb N},\;m\geq0\};\)

\item
\(\{\ELambda^n_k\hookrightarrow\Sim^n\mid n\geq2,\;k=0,n\}
\Box
\{C_{2r}[1]\hookrightarrow C_{2r}[\varepsilon]\mid
r\in\overline{\mathbb N}\};\)

\item
\(\{\eqref{L1},\eqref{L2}\}
\Box
\{\partial C_r[m]\hookrightarrow C_r[m]\mid
r\in\overline{\mathbb N},\;m\geq0\};\)

\item
\(\{\eqref{L1},\eqref{L2}\}
\Box
\{C_{2r}[1]\hookrightarrow C_{2r}[\varepsilon]\mid
r\in\overline{\mathbb N}\}.\)
\end{enumerate}
Indeed, for fixed \(u\), the functor \(u\Box(-)\) preserves colimits
and retracts.

Case~(1) follows from Proposition~\ref{prop:key} and the duality
principle~\ref{lem:dual}.
The morphisms in~(2) are isomorphisms, since for \(n\geq2\)
the inclusion \(\ELambda_k^n\hookrightarrow\Sim^n\) is a bijection
on \(\varepsilon\)-witnesses.
By the duality principle~\ref{lem:dual}, in cases~(3)--(4)
it is enough to consider~\eqref{L1}.
Recall~\eqref{L1}:
\[
\tag{L1}
\begin{aligned}
    \begin{tikzcd}[row sep={1.8em,between origins}, column sep={2.4em,between origins}, ampersand replacement=\&,execute at end picture={
\foreach \Nombre in  {A,B,C}
  {\coordinate (\Nombre) at (\Nombre.center);}
\fill[pattern=vertical lines,opacity=0.3]
  (A) -- (B) -- (C) -- cycle;
  \foreach \Nombre in  {A,C,D}
  {\coordinate (\Nombre) at (\Nombre.center);}
\fill[pattern=vertical lines,opacity=0.3]
  (B) -- (C) -- (D) -- cycle;
}]
    \&
    |[alias=C]|2\arrow[rr,"e_3"{description}] \&\& |[alias=D]|3
    \\[10pt]\&
    \&\&\\
    |[alias=A]|0
    \arrow[ruu,LA=1.5pt, "e_1"{description}]
    \arrow[rr,LA=1.5pt,"e_4"{description}]
    \&\&
    |[alias=B]|1\arrow[uul, "a"{description}]
    \arrow[ruu,LA=1.5pt,"e_2"{description}]\&
\end{tikzcd}
\quad &\hookrightarrow\quad
\begin{tikzcd}[row sep={1.8em,between origins}, column sep={2.4em,between origins}, ampersand replacement=\&,execute at end picture={
\foreach \Nombre in  {A,B,C}
  {\coordinate (\Nombre) at (\Nombre.center);}
\fill[pattern=vertical lines,opacity=0.3]
  (A) -- (B) -- (C) -- cycle;
  \foreach \Nombre in  {A,C,D}
  {\coordinate (\Nombre) at (\Nombre.center);}
\fill[pattern=vertical lines,opacity=0.3]
  (B) -- (C) -- (D) -- cycle;
}]
    \&
    |[alias=C]|2\arrow[rr,LA=1.5pt,"e_3"{description}] \&\& |[alias=D]|3
    \\[10pt]\&
    \&\&\\
    |[alias=A]|0
    \arrow[ruu,LA=1.5pt, "e_1"{description}]
    \arrow[rr,LA=1.5pt,"e_4"{description}]
    \&\&
    |[alias=B]|1\arrow[uul, "a"{description}]
    \arrow[ruu,LA=1.5pt,"e_2"{description}]\&
\end{tikzcd}
\end{aligned}
\]
and denote this morphism by \(l\colon A\hookrightarrow B\).

\smallskip
\noindent
\emph{Case (3).}
For \(m\geq1\), the morphism
\(l\Box(\partial C_r[m]\hookrightarrow C_r[m])\)
is an isomorphism.
For \(m=0\), it has the form
\[
l\times C_r[0]\colon
A\times C_r[0]\hookrightarrow B\times C_r[0].
\]
This morphism is a bijection on simplices. The only additional
\(\varepsilon\)-witnesses in the codomain are those on the edges
\((e_3,c_i)\), where
\(c_i=[\langle i,i+1\rangle]\).
For finite \(r\), the indices are taken modulo \(r\); for
\(r=\infty\), they range over \(\mathbb Z\).
Put \(b_i=[\langle i-1,i\rangle]\).
In \(A\times C_r[0]\) we have the rhombus
\[
\begin{tikzcd}[row sep={1.8em,between origins}, column sep={5.2em,between origins}, ampersand replacement=\&,execute at end picture={
\foreach \Nombre in  {A,B,C}
  {\coordinate (\Nombre) at (\Nombre.center);}
\fill[pattern=vertical lines,opacity=0.3]
  (A) -- (B) -- (C) -- cycle;
  \foreach \Nombre in  {A,C,D}
  {\coordinate (\Nombre) at (\Nombre.center);}
\fill[pattern=vertical lines,opacity=0.3]
  (B) -- (C) -- (D) -- cycle;
}]
    \&
    |[alias=C]|(2,i)\arrow[rr,"{(e_3,c_i)}"{description}]
    \&\& |[alias=D]|(3,i+1)
    \\[10pt]\&
    \&\&\\
    |[alias=A]|(0,i-1)
    \arrow[ruu,LA=1.5pt,"{(e_1,b_i)}"{description}]
    \arrow[rr,LA=1.5pt,"{(e_4,b_i)}"{description}]
    \&\&
    |[alias=B]|(1,i)
    \arrow[uul,"{(a,i)}"{description}]
    \arrow[ruu,LA=1.5pt,"{(e_2,c_i)}"{description}]
    \&
\end{tikzcd}
\]
Hence the missing \(\varepsilon\)-witness on \((e_3,c_i)\)
can be attached by a pushout along~\eqref{L1}.
Attaching these witnesses one by one proves~(3).

\smallskip
\noindent
\emph{Case (4).}
First consider
\(C_{2r}[1]\hookrightarrow C_{2r}[\varepsilon]\), for \(r\geq1\).
The corresponding pushout-product is
\[
v\colon
A\times C_{2r}[\varepsilon]
\coprod_{A\times C_{2r}[1]}
B\times C_{2r}[1]
\longrightarrow
B\times C_{2r}[\varepsilon].
\]
The only additional \(\varepsilon\)-witnesses in the codomain are
those on the edges
\((e_3,b_i)\), where
\(b_i=[\langle2i,2i+1\rangle]\)
and \(i\) is taken modulo \(r\).
For each such \(i\), the domain of \(v\) contains the rhombus
\[
\begin{tikzcd}[row sep={1.8em,between origins}, column sep={5.2em,between origins}, ampersand replacement=\&,execute at end picture={
\foreach \Nombre in  {A,B,C}
  {\coordinate (\Nombre) at (\Nombre.center);}
\fill[pattern=vertical lines,opacity=0.3]
  (A) -- (B) -- (C) -- cycle;
  \foreach \Nombre in  {A,C,D}
  {\coordinate (\Nombre) at (\Nombre.center);}
\fill[pattern=vertical lines,opacity=0.3]
  (B) -- (C) -- (D) -- cycle;
}]
    \&
    |[alias=C]|(2,2i)
    \arrow[rr,"{(e_3,b_i)}"{description}]
    \&\& |[alias=D]|(3,2i+1)
    \\[10pt]\&
    \&\&\\
    |[alias=A]|(0,2i-2)
    \arrow[ruu,LA=1.5pt]
    \arrow[rr,LA=1.5pt,"{(e_4,b_{i-1})}"{description}]
    \&\&
    |[alias=B]|(1,2i-1)
    \arrow[uul]
    \arrow[ruu,LA=1.5pt]
    \&
\end{tikzcd}
\]
and therefore the missing \(\varepsilon\)-witness on
\((e_3,b_i)\) can be attached by a pushout along~\eqref{L1}.
Attaching these witnesses one by one proves the finite case.
For \(C_\infty[1]\hookrightarrow C_\infty[\varepsilon]\),
the same argument applies with \(i\in\mathbb Z\).

This proves~(4), and hence both assertions of the theorem.
\end{proof}

\section{$\eta\varepsilon$-simplicial sets}

Our goal is to associate to $\varepsilon$-fibrant simplicial sets $X$ and $Y$
an $\varepsilon$-fibrant
$\mathrm{Fun}^\eta(X,Y)\subseteq\mathrm{Fun}(X,Y)$.
The previous section indicates that the only issue
may occur for $\varepsilon$-horns in
dimension $1$ and can be described as follows.
Suppose that $X$ contains two $\varepsilon$-edges $e_1,e_2$ with a common
source $x$. Then there is a $2$-simplex
$\sigma\colon\Delta^2\to X$ such that
$d_1\sigma=e_1$, $d_2\sigma=e_2$, and
$d_0\sigma$ is an $\eta$-edge, denoted by $a$.
The prism $\Sim^1\times\sigma$ has the following form:
\[
\begin{tikzcd}[row sep=small]
    && {(1,2)} &\\
    {(1,0)}\arrow[rru]\arrow[rrr]
    &&&{(1,1)}
        \arrow[lu,"{(\langle1,1\rangle,a)}"',pos=.45]\\
    &&{(0,2)}\arrow[uu]&\\
    {(0,0)}
        \arrow[uu]
        \arrow[uuurr,LA=1.5pt,
            "{(\langle0,1\rangle,e_1)}"description,pos=.4]
        \arrow[rrruu,LA=1.5pt,
            "{(\langle0,1\rangle,e_2)}"',pos=.75]
        \arrow[rru]
        \arrow[rrr]
    &&&{(0,1)}\arrow[lu]\arrow[uu]
\end{tikzcd}
\]
Any morphism $\Sim^1\times X\to Y$ sends
$(\langle1,1\rangle,a)$ to an $\eta$-edge.
On the other hand,
$\{1\}\times X$ has no $\varepsilon$-edges, since products in
$\eSSet$ are computed componentwise and $\{1\}$ has no
$\varepsilon$-witnesses.
Thus a morphism $\{1\}\times X\to Y$ is not required
a priori to preserve $\eta$-edges.
If it fails to preserve the edge $a$, it cannot
extend to $\Sim^1\times X\to Y$.

Fortunately, this is essentially the only obstruction. It can be removed
by restricting $\mathrm{Fun}(X,Y)$ to simplices whose vertices correspond to
$\eta$-edge-preserving morphisms. To implement this restriction, we
introduce $\eta\varepsilon$-simplicial sets, which retain the information
about $\eta$-edges in addition to the $\varepsilon$-structure.

Most results and proofs in this section are direct analogues of those in
the previous section, so we will proceed more quickly whenever the
arguments are unchanged.

\begin{defn}
    Let $\Delta^{\eta\varepsilon}$ denote the category obtained from $\Sim$ by adding an extra object $[\eta]$ and a splitting of $\sigma^0_0$ as
\begin{equation}
\begin{tikzcd}[ampersand replacement=\&]
    {[0]} \&\& {[1]} \arrow[ll,"\sigma_0^0"] \arrow[ld,"\eta"]\\
    \& {[\eta]} \arrow[lu] \&
\end{tikzcd}
\end{equation}
   By an $\eta\varepsilon$-simplicial set we mean a presheaf on $\Delta^{\eta\varepsilon}$.
\end{defn}
Note that the definition is similar to that of marked simplicial sets.
The idea is that $\varepsilon$-edges are those edges which ought to be degenerate as images under $s_1^0$ of vertices, while $\eta$-edges are those edges which ought to be degenerate as images under $s^0_0$ of vertices.
Note that each degenerate edge of an $\eta\varepsilon$-simplicial set has a unique \emph{degenerate \(\eta\)-witness}.
For an $\eta\varepsilon$-simplicial set, by $\eta$-edge we mean an edge that has an $\eta$-witness, not an edge satisfying Definition~\ref{def:etaEdge}.

There is a commutative triangle of functors
\[
\begin{tikzcd}
    {\Lambda_\infty}&&{\Sim}\arrow[ll,"P"]\arrow[d,"I"]\\
    &&{\Delta^{\eta\varepsilon}}\arrow[llu,"Q"]
\end{tikzcd}
\]
where $I$ is the obvious inclusion and $Q$ sends $Q([1]\rightarrow [\eta])=\sigma_0^0$ and  $Q([\eta]\rightarrow [0])=\mathrm{id}_{[0]}$.
The commutative triangle induces $P_!\cong Q_!I_!$, $P^*=I^*Q^*$, and  $P_*\cong Q_*I_*$.
Moreover, for each $F=P,Q,I$, we denote by $\eta^F$ and $\epsilon^F$ the natural transformations
\[\eta^F\colon \mathrm{Id}\rightarrow F^*F_!\quad\text{resp.}\quad \epsilon^F\colon F^*F_*\rightarrow \mathrm{Id}.\]
It is easy to see that the functor $I_!\colon \eSSet\rightarrow \neSSet$ is a fully faithful embedding which adds degenerate $\eta$-witnesses to degenerate edges, while $I^*$ forgets all $\eta$-data.

Via the fully faithful embedding $I_!$, we will regard
$\varepsilon$-simplicial sets such as $\Delta^n$ and $\Sim^n$ as
$\eta\varepsilon$-simplicial sets whenever no confusion can arise.
Moreover, we will denote by $\etaSim^1$ the $\eta\varepsilon$-simplicial set having a unique nondegenerate $\eta$-edge, that is $\etaSim^1=y([\eta])$.

Under the embedding $I_!$ we can apply most of the results about $\eSSet$ from previous sections to $\neSSet$. The key observation is the following.
\begin{lemma}\label{lem:recycling}
There is a natural isomorphism
\[
\gamma_X:=\epsilon_{Q^*Q_!I_!X}\colon
I_!P^*P_!X
\cong
I_!I^*Q^*Q_!I_!X
\longrightarrow
Q^*Q_!I_!X,
\]
where $\epsilon\colon I_!I^*\to\mathrm{Id}$ is the counit of
$I_!\dashv I^*$.
Moreover,
\[
\eta^Q_{I_!X}=\gamma_X\circ I_!(\eta^P_X).
\]
\end{lemma}

\begin{proof}
For any $Y\in\neSSet$, the object $I_!I^*Y$ is obtained from $Y$ by
forgetting all nondegenerate $\eta$-witnesses. Hence
$\epsilon_Y\colon I_!I^*Y\to Y$ is an isomorphism whenever $Y$
contains only degenerate $\eta$-witnesses. In particular, this holds for
objects in the image of $Q^*$, and therefore
$\gamma_X=\epsilon_{Q^*Q_!I_!X}$ is an isomorphism.

For the second statement, using the formula for the unit of a composite
adjunction, naturality of $\epsilon$, and the triangle identity, we get
\[
\begin{aligned}
\gamma_X\circ I_!(\eta^P_X)
&=\epsilon_{Q^*Q_!I_!X}\circ I_!(\eta^P_X)\\
&=\epsilon_{Q^*Q_!I_!X}
  \circ I_!I^*(\eta^Q_{I_!X})
  \circ I_!(\eta^I_X)\\
&=\eta^Q_{I_!X}\circ
  \epsilon_{I_!X}\circ I_!(\eta^I_X)\\
&=\eta^Q_{I_!X}.
\end{aligned}
\]
\end{proof}

Working with $\eta\varepsilon$-simplicial sets, $\eta$-edges will be pictured as dotted arrows. Consider the following $\eta\varepsilon$-simplicial set morphisms:
\begin{equation}\label{I1}
\tag{I1}
\begin{aligned}
       \begin{tikzcd}[row sep={2.4em,between origins}, column sep={3.3em,between origins}, ampersand replacement=\&,execute at end picture={
\foreach \Nombre in  {A,B,C}
  {\coordinate (\Nombre) at (\Nombre.center);}
\fill[pattern=vertical lines,opacity=0.3]
  (A) -- (B) -- (C) -- cycle;
}]
    \&
    |[alias=C]|2 \&
    \\[10pt]\&
    \&\\
    |[alias=A]|0
    \arrow[ruu,LA=1.5pt, "e_1"{description}]
    \arrow[rr, LA=1.5pt,"e_2"{description}]
    \&\&
    |[alias=B]|1\arrow[uul, "a"{description}]
\end{tikzcd}  &\quad \hookrightarrow \quad
        \begin{tikzcd}[row sep={2.4em,between origins}, column sep={3.3em,between origins}, ampersand replacement=\&,execute at end picture={
\foreach \Nombre in  {A,B,C}
  {\coordinate (\Nombre) at (\Nombre.center);}
\fill[pattern=vertical lines,opacity=0.3]
  (A) -- (B) -- (C) -- cycle;
}]
    \&
    |[alias=C]|2
    \&\\[10pt]\&
    \&\\
    |[alias=A]|0
    \arrow[ruu,LA=1.5pt, "e_1"{description}]
    \arrow[rr, LA=1.5pt,"e_2"{description}]
    \&\&
    |[alias=B]|1\arrow[uul, dotted, "a"{description}]
\end{tikzcd}
\end{aligned}
\end{equation}

\begin{equation}\label{I2}
\tag{I2}
\begin{aligned}
       \begin{tikzcd}[row sep={2.4em,between origins}, column sep={3.3em,between origins}, ampersand replacement=\&,execute at end picture={
\foreach \Nombre in  {A,B,C}
  {\coordinate (\Nombre) at (\Nombre.center);}
\fill[pattern=vertical lines,opacity=0.3]
  (A) -- (B) -- (C) -- cycle;
}]
    \&
    |[alias=C]|2 \&
    \\[10pt]\&
    \&\\
    |[alias=A]|0
    \arrow[ruu,LA=1.5pt, "e_1"{description}]
    \arrow[rr,"a"{description}]
    \&\&
    |[alias=B]|1\arrow[uul,LA=1.5pt, "e_2"{description}]
\end{tikzcd}  &\quad \hookrightarrow \quad
         \begin{tikzcd}[row sep={2.4em,between origins}, column sep={3.3em,between origins}, ampersand replacement=\&,execute at end picture={
\foreach \Nombre in  {A,B,C}
  {\coordinate (\Nombre) at (\Nombre.center);}
\fill[pattern=vertical lines,opacity=0.3]
  (A) -- (B) -- (C) -- cycle;
}]
    \&
    |[alias=C]|2 \&
    \\[10pt]\&
    \&\\
    |[alias=A]|0
    \arrow[ruu,LA=1.5pt, "e_1"{description}]
    \arrow[rr, dotted,"a"{description}]
    \&\&
    |[alias=B]|1\arrow[uul,LA=1.5pt, "e_2"{description}]
\end{tikzcd}
\end{aligned}
\end{equation}

\begin{equation}\label{R1}
\tag{R1}
\begin{aligned}
       \begin{tikzcd}[row sep={2.4em,between origins}, column sep={3.3em,between origins}, ampersand replacement=\&,execute at end picture={
\foreach \Nombre in  {A,B,C}
  {\coordinate (\Nombre) at (\Nombre.center);}
\fill[pattern=vertical lines,opacity=0.3]
  (A) -- (B) -- (C) -- cycle;
}]
    \&
    |[alias=C]|2 \&
    \\[10pt]\&
    \&\\
    |[alias=A]|0
    \arrow[ruu,LA=1.5pt, "e_1"{description}]
    \arrow[rr, dotted,"a"{description}]
    \&\&
    |[alias=B]|1\arrow[uul, "e_2"{description}]
\end{tikzcd}  &\quad \hookrightarrow \quad
         \begin{tikzcd}[row sep={2.4em,between origins}, column sep={3.3em,between origins}, ampersand replacement=\&,execute at end picture={
\foreach \Nombre in  {A,B,C}
  {\coordinate (\Nombre) at (\Nombre.center);}
\fill[pattern=vertical lines,opacity=0.3]
  (A) -- (B) -- (C) -- cycle;
}]
    \&
    |[alias=C]|2 \&
    \\[10pt]\&
    \&\\
    |[alias=A]|0
    \arrow[ruu,LA=1.5pt, "e_1"{description}]
    \arrow[rr, dotted,"a"{description}]
    \&\&
    |[alias=B]|1\arrow[uul,LA=1.5pt, "e_2"{description}]
\end{tikzcd}
\end{aligned}
\end{equation}

\begin{equation}\label{R2}
\tag{R2}
\begin{aligned}
       \begin{tikzcd}[row sep={2.4em,between origins}, column sep={3.3em,between origins}, ampersand replacement=\&,execute at end picture={
\foreach \Nombre in  {A,B,C}
  {\coordinate (\Nombre) at (\Nombre.center);}
\fill[pattern=vertical lines,opacity=0.3]
  (A) -- (B) -- (C) -- cycle;
}]
    \&
    |[alias=C]|2 \&
    \\[10pt]\&
    \&\\
    |[alias=A]|0
    \arrow[ruu,LA=1.5pt, "e_2"{description}]
    \arrow[rr,"e_1"{description}]
    \&\&
    |[alias=B]|1\arrow[uul, dotted, "b"{description}]
\end{tikzcd}  &\quad \hookrightarrow \quad
         \begin{tikzcd}[row sep={2.4em,between origins}, column sep={3.3em,between origins}, ampersand replacement=\&,execute at end picture={
\foreach \Nombre in  {A,B,C}
  {\coordinate (\Nombre) at (\Nombre.center);}
\fill[pattern=vertical lines,opacity=0.3]
  (A) -- (B) -- (C) -- cycle;
}]
    \&
    |[alias=C]|2 \&
    \\[10pt]\&
    \&\\
    |[alias=A]|0
    \arrow[ruu,LA=1.5pt, "e_2"{description}]
    \arrow[rr, LA=1.5pt,"e_1"{description}]
    \&\&
    |[alias=B]|1\arrow[uul,dotted, "b"{description}]
\end{tikzcd}
\end{aligned}
\end{equation}
If an $\eta\varepsilon$-simplicial set $X$ has the right lifting property with respect to the morphisms~(I1--2) and (R1--2), then we say $X$ satisfies the \emph{$\eta\varepsilon$-coherence conditions}.
The four conditions~\eqref{I1}--\eqref{R2} guarantee that rotations of $\varepsilon$-edges are $\eta$-edges and vice versa. An $\eta\varepsilon$-simplicial set which admits fillers for
$\varepsilon$-horns and satisfies the coherence conditions (I1-2) and (R1-2) is called \emph{$\eta\varepsilon$-fibrant}.

For $X\in\eSSet$, define $X^\eta$ to be the $\eta\varepsilon$-simplicial set, where each degenerate edge and each $\eta$-edge (in the sense of Definition~\ref{def:etaEdge}) has a unique $\eta$-witness.
Note that $(-)^\eta$ extends to a functor $\eSSet\rightarrow\neSSet$.

The following lemma is easy.
\begin{lemma}\label{lem:bridge}
    Let $X\in\eSSet$ and $Y\in\neSSet$.
    Then $X$ is $\varepsilon$-fibrant if and only if $X^\eta$ is
$\eta\varepsilon$-fibrant. Moreover, if $Y$ is
$\eta\varepsilon$-fibrant, then $I^*Y$ is $\varepsilon$-fibrant.
\end{lemma}
\begin{proof}
    Since the construction $X \mapsto X^\eta$ leaves the
    underlying simplicial set and its $\varepsilon$-edges
    unchanged, $X$ admits fillers for $\varepsilon$-horns
    if and only if $X^\eta$ does.

    Suppose first that $X$ is $\varepsilon$-fibrant.
    By definition, $X^\eta$ satisfies the right lifting
    property with respect to (I1) and (I2).
    To verify (R1), let $A \hookrightarrow B$ be the
    corresponding inclusion and consider a map
    $A \to X^\eta$.
    By Lemma~\ref{lem:defOfEtaEdge}, this map extends
    to a rhombus with three $\varepsilon$-edges.
    Theorem~\ref{thm:4for3} implies that the remaining
    outer edge is also an $\varepsilon$-edge, giving
    the required extension to $B$.
    The verification of (R2) is analogous.
    Thus $X^\eta$ is $\eta\varepsilon$-fibrant.

    Now let $Y$ be any $\eta\varepsilon$-fibrant object.
    Then $I^*Y$ admits fillers for $\varepsilon$-horns.
    To verify its $\varepsilon$-coherence conditions,
    consider a map $A \to I^*Y$, where $A$ is the
    domain of (L1) or (L2).
    Applying (I1) or (I2) in $Y$ shows that the
    diagonal of the corresponding rhombus is
    an $\eta$-edge.
    Applying (R1) or (R2) then shows that the remaining
    outer edge is an $\varepsilon$-edge.
    Hence $I^*Y$ satisfies (L1) and (L2), and is
    therefore $\varepsilon$-fibrant.

    In particular, if $X^\eta$ is
    $\eta\varepsilon$-fibrant, then
    $X = I^*(X^\eta)$ is $\varepsilon$-fibrant.
\end{proof}
For an $\varepsilon$-simplicial set $X$ without nontrivial $\eta$-edges we have $X^\eta\cong I_!X$, hence we can safely identify
\[C^\eta_r[n]:=(C_r[n])^\eta\cong I_!C_r[n].\]
Note that by Lemma~\ref{lem:recycling} we can represent $\eta^Q_{\Delta^n}\colon \Delta^n\rightarrow Q^*Q_!\Delta^n$ as $\Delta^n\hookrightarrow C^\eta_\infty[n]$. It follows that we have a natural identification
\[Q^*Q_*(X)_n\cong\Hom(C^\eta_\infty[n],X)\]
Moreover, we have the following identifications:
\begin{lemma}
Under the identifications $Q^*y([n]) = C^\eta_\infty[n]$,
we have
\begin{align}
    Q^*Q_!(\Delta^1 \to \Sim^1)
    &\cong \bigl(q\colon C_\infty^\eta[1]
    \to C_\infty^\eta[0]\bigr),
    && q(2i)=q(2i-1)=i,
    \label{eq:QQeps}\\
    Q^*Q_!(\Delta^1 \to \etaSim^1)
    &\cong \bigl(p\colon C_\infty^\eta[1]
    \to C_\infty^\eta[0]\bigr),
    && p(2i)=p(2i+1)=i,
    \label{eq:QQeta}\\
    Q^*Q_!(\etaSim^1 \to \Delta^0)
    &\cong \bigl(\mathrm{id}\colon C_\infty^\eta[0]
    \to C_\infty^\eta[0]\bigr),
    &&
    \label{eq:QQid}
\end{align}
where $i\in\mathbb{Z}$.
\end{lemma}
\begin{proof}
    Formula~\eqref{eq:QQeps} follows by the same argument
    as in the proof of Lemma~\ref{lem:defOfq}.
    Formula~\eqref{eq:QQeta} is proved analogously.
    Finally, formula~\eqref{eq:QQid} follows from
    $Q([\eta]\to[0])=\mathrm{id}_{[0]}$
    and the natural identifications
    $Q_!y(a)\cong y(Qa)$.
\end{proof}
In analogy with $\varepsilon$-simplicial sets, we now enlarge \(Q^*Q_*X\) by adding the missing \(\eta\)- and \(\varepsilon\)-witnesses.
For each $1\leq l\leq\infty$, we define
\[
C^\eta_l[\varepsilon]=\big(C_l[\varepsilon]\big)^\eta
\qquad\text{and}\qquad
C^\eta_l[\eta]=\bigg(C_l[1]\cup\bigcup_i\Sim^{\{2i-1,2i\}}\bigg)^\eta.
\]
The following diagram illustrates $C^\eta_\infty[\varepsilon]$:
\begin{equation*}
  \cdots
\begin{tikzcd}[row sep=2.2em, column sep=1.6em]
&-2\arrow[rr,LA=1.5pt]\arrow[rd,LA=1.5pt]
&&0\arrow[rr,LA=1.5pt]\arrow[rd,LA=1.5pt]
&&2\arrow[rr,LA=1.5pt]\arrow[rd,LA=1.5pt]
&&4\arrow[rr,LA=1.5pt]\arrow[rd,LA=1.5pt]
&&6\\
{-3}\arrow[rr,LA=1.5pt]\arrow[ru,dotted]
&&{-1}\arrow[rr,LA=1.5pt]\arrow[ru,dotted]
&&{1}\arrow[rr,LA=1.5pt]\arrow[ru,dotted]
&&{3}\arrow[rr,LA=1.5pt]\arrow[ru,dotted]
&&{5}\arrow[ru,dotted]&
\end{tikzcd}
\cdots
\end{equation*}
Note that shifting the indices by one induces isomorphisms
$C^\eta_l[\eta]\cong C^\eta_l[\varepsilon]$.
The maps $p$ and $q$ from~(\ref{eq:QQeps}--\ref{eq:QQeta}) admit unique extensions
\[q_\varepsilon\colon C^\eta_\infty[\varepsilon]\to C^\eta_\infty[0]\qquad\text{and}\qquad p_\eta\colon C^\eta_\infty[\eta]\to C^\eta_\infty[0].\]

\begin{defn}
For $X\in\neSSet$, define
$\widetilde X$ to have the same underlying simplicial set as $Q^*Q_*X$, with
\begin{equation}
\widetilde X_\eta=\Hom(C^\eta_\infty[\eta],X),
\qquad
\widetilde X_\varepsilon=\Hom(C^\eta_\infty[\varepsilon],X).
\end{equation}
Its structure maps $\widetilde{X}_0\rightarrow \widetilde{X}_\eta\rightarrow \widetilde{X}_1$ and $\widetilde{X}_\varepsilon\rightarrow \widetilde{X}_1$ are induced by
\begin{align}
&\Hom(C^\eta_\infty[0],X)\xrightarrow{p_\eta^*}
\Hom(C^\eta_\infty[\eta],X)\xrightarrow{\mathrm{res}}
\Hom(C^\eta_\infty[1],X) \label{eq:structureMaps} \\
&\Hom(C^\eta_\infty[\varepsilon],X)\xrightarrow{\mathrm{res}}
\Hom(C^\eta_\infty[1],X).
\end{align}

Moreover, define $\tilde{\epsilon}_X\colon \widetilde{X}\rightarrow X$ to be an extension of $\epsilon_X^Q$ by
\begin{align*}
    (\tilde{\epsilon}_X)_\varepsilon\colon& \Hom(C^\eta_\infty[\varepsilon],X)\xrightarrow{\mathrm{res}}\Hom(\Sim^{\{0,1\}},X),\\
(\tilde{\epsilon}_X)_\eta\colon& \Hom(C^\eta_\infty[\eta],X)\xrightarrow{\mathrm{res}}\Hom(\etaSim^{\{0,1\}},X).
\end{align*}
\end{defn}

\begin{lemma}
Let $X$ be an $\eta\varepsilon$-simplicial set. Then $\widetilde{X}$ is well-defined as an $\eta\varepsilon$-simplicial set and the maps
\[
(Q^*Q_*X)_\eta\cong\Hom(C^\eta_\infty[0],X)
\xrightarrow{p_\eta^*}\widetilde X_\eta,
\quad
(Q^*Q_*X)_\varepsilon\cong\Hom(C^\eta_\infty[0],X)
\xrightarrow{q_\varepsilon^*}\widetilde X_\varepsilon
\]
extend to an embedding
$j\colon Q^*Q_*X\hookrightarrow\widetilde X$
which is the identity on the underlying simplicial set.

Moreover, $\tilde{\epsilon}_X\colon\widetilde{X}\rightarrow X$ is a well-defined morphism and $\tilde{\epsilon}_X j=\epsilon^Q_X$.
\end{lemma}
\begin{proof}
To see that $\widetilde X$ is well-defined, it is enough to check that the composition of~\eqref{eq:structureMaps} agrees with
$s^0_0\colon\widetilde X_0\to\widetilde X_1$. This follows since

$$
C^\eta_\infty[1]\longrightarrow C^\eta_\infty[\eta]
\xrightarrow{p_\eta}C^\eta_\infty[0]
$$

equals $p\colon C^\eta_\infty[1]\to C^\eta_\infty[0]$, which represents
$Q^*Q_!(\Delta^1\to\Delta^0)$.

Next we show that $j$ preserves the structure maps. Compatibility with
$[1]\to[\varepsilon]$ follows from
\begin{equation}
\begin{tikzcd}
{(Q^*Q_*X)_1}\arrow[r,"\cong"]
&{\Hom(C_\infty^\eta[1],X)}\arrow[r,equal]
&{\Hom(C^\eta_\infty[1],X)}\arrow[r,"\cong"]
&{\widetilde X_1}\\
{(Q^*Q_*X)_\varepsilon}\arrow[r,"\cong"]\arrow[u]
&{\Hom(C_\infty^\eta[0],X)}\arrow[u,"q^*"]\arrow[r,"q_\varepsilon^*"]
&{\Hom(C^\eta_\infty[\varepsilon],X)}\arrow[u,"\mathrm{res.}"]\arrow[r,"\cong"]
&{\widetilde X_\varepsilon}\arrow[u]
\end{tikzcd}
\end{equation}
and compatibility with $[1]\to[\eta]\to[0]$ from
\begin{equation}
\begin{tikzcd}
{(Q^*Q_*X)_1}\arrow[r,"\cong"]
&{\Hom(C^\eta_\infty[1],X)}\arrow[r,equal]
&{\Hom(C^\eta_\infty[1],X)}\arrow[r,"\cong"]
&{\widetilde X_1}\\
{(Q^*Q_*X)_\eta}\arrow[u]\arrow[r,"\cong"]
&{\Hom(C^\eta_\infty[0],X)}\arrow[u,"p^*"]\arrow[r,"p_\eta^*"]
&{\Hom(C^\eta_\infty[\eta],X)}\arrow[u,"\mathrm{res.}"]\arrow[r,"\cong"]
&{\widetilde X_\eta}\arrow[u]\\
{(Q^*Q_*X)_0}\arrow[u]\arrow[r,"\cong"]
&{\Hom(C^\eta_\infty[0],X)}\arrow[u,equal]\arrow[r,equal]
&{\Hom(C^\eta_\infty[0],X)}\arrow[u,"p_\eta^*"]\arrow[r,"\cong"]
&{\widetilde X_0}\arrow[u].
\end{tikzcd}
\end{equation}
Thus the prescribed maps extend to a morphism
$j\colon Q^*Q_*X\to\widetilde X$. Since $p_\eta$ and $q_\varepsilon$ are
surjective, $p_\eta^*$ and $q_\varepsilon^*$ are injective. Hence $j$ is
an embedding.

It remains to check that $\tilde\epsilon_X$ preserves the structure maps.
The essential compatibilities are obtained by applying $\Hom(-,X)$ to
\begin{equation*}
\begin{tikzcd}
{C^\eta_\infty[\varepsilon]}
&{\Sim^{\{0,1\}}}\arrow[l,hook]\\
{C^\eta_\infty[1]}\arrow[u,hook]
&{\Delta^1}\arrow[l,hook]\arrow[u,hook]
\end{tikzcd}
\qquad
\begin{tikzcd}
{C^\eta_\infty[0]}&{\Delta^0}\arrow[l]\\
{C^\eta_\infty[\eta]}\arrow[u,"p_\eta"]
&{\etaSim^{\{0,1\}}}\arrow[l,hook]\arrow[u]\\
{C^\eta_\infty[1]}\arrow[u,hook]
&{\Delta^1}\arrow[l,hook]\arrow[u,hook].
\end{tikzcd}
\end{equation*}
Hence $\tilde\epsilon_X\colon\widetilde X\to X$ is a morphism.

Finally, $\tilde\epsilon_Xj=\epsilon^Q_X$ already holds on the underlying
simplicial set, so it remains only to check the $\varepsilon$- and
$\eta$-components. These are given by the following commutative triangles, which are easy to verify:
\begin{equation*}
\begin{tikzcd}[column sep=0.8em]
{\Hom(C^\eta_\infty[0],X)}
\arrow[rd,"{(\eta^Q_{\Sim^1})^*}"']
\arrow[r,"q_\varepsilon^*"]
&{\Hom(C^\eta_\infty[\varepsilon],X)}
\arrow[d,"\mathrm{res.}"]\\
&{\Hom(\Sim^1,X)}
\end{tikzcd}
\quad
\begin{tikzcd}[column sep=0.8em]
{\Hom(C^\eta_\infty[0],X)}
\arrow[rd,"{(\eta^Q_{\etaSim^1})^*}"']
\arrow[r,"p_\eta^*"]
&{\Hom(C^\eta_\infty[\eta],X)}
\arrow[d,"\mathrm{res.}"]\\
&{\Hom(\etaSim^1,X)}.
\end{tikzcd}
\end{equation*}
Therefore $\tilde\epsilon_Xj=\epsilon^Q_X$.
\end{proof}

\begin{lemma}\label{lem:etaCellStructure}
For each $X\in\neSSet$, the morphism
$\emptyset\hookrightarrow\widetilde X$ belongs to the weakly saturated
class generated by
\[
\partial C^\eta_l[n]\hookrightarrow C^\eta_l[n]
\qquad\text{and}\qquad
C^\eta_{2r}[1]\hookrightarrow C^\eta_{2r}[\varepsilon],
\]
where $n\geq0$ and $1\leq l,r\leq\infty$ (we set $2\infty=\infty$).
\end{lemma}
\begin{proof}
    This is a direct modification of the proof of
Lemma~\ref{lem:cellularStructureOfHatX}. Set $A=Q^*Q_*X$ and identify it with its image in
$\widetilde X$. By the same argument,
$\emptyset\hookrightarrow A$ belongs to the weakly saturated class
generated by $(\partial C_l[n]\hookrightarrow C_l[n])^\eta$,
$n\geq0$, $1\leq l\leq\infty$.

It remains to attach the new witnesses.
Note that precomposition with the index shift $j\mapsto j+1$
between $C^\eta_\infty[\eta]$ and
$C^\eta_\infty[\varepsilon]$ defines a bijection
\[
    S\colon
    \widetilde X_\varepsilon\amalg\widetilde X_\eta
    \longrightarrow
    \widetilde X_\varepsilon\amalg\widetilde X_\eta
\]
interchanging the two summands.
We first show that $S$ preserves
$A_\varepsilon\amalg A_\eta$ and its complement.
By the formulas for $p$ and $q$ in
\eqref{eq:QQeps}--\eqref{eq:QQeta},
a morphism $f\colon C^\eta_\infty[\varepsilon]\to X$
factors through $q_\varepsilon$ if and only if
$S(f)$ factors through $p_\eta$.
Similarly, a morphism
$g\colon C^\eta_\infty[\eta]\to X$
factors through $p_\eta$ if and only if
$S(g)$ factors through $q_\varepsilon$.
Thus $S$ preserves the original witnesses as well as the new witnesses.
Hence, the new witnesses decompose into $S$-orbits,
each infinite or of even size $2r$.

Choose an $\varepsilon$-witness $f$ in such an orbit,
and let $a\in A_1$ be its underlying edge.
The canonical orbit morphism
$\rho_a\colon C^\eta_\infty[1]\to A$
extends to
\[
    \rho_f\colon
    C^\eta_\infty[\varepsilon]\longrightarrow\widetilde X
\]
by sending the additional witness over
$\langle j,j+1\rangle$ to $S^j(f)$:
this is an $\varepsilon$-witness for even $j$ and
an $\eta$-witness for odd $j$.
These assignments are compatible with the underlying
edges, since restriction commutes with the index shift,
as in Lemma~\ref{lem:canonicalOrbitsEps}.

If the orbit has size $2r$, with $r<\infty$, the equality $S^{2r}(f)=f$
makes $\rho_f$ invariant under the shift by $2r$,
so it factors through $C^\eta_{2r}[\varepsilon]$.
Its additional witnesses map bijectively onto the orbit,
while its restriction to $C^\eta_{2r}[1]$ maps into $A$.
Consequently, adjoining this orbit is a pushout of
\[
    C^\eta_{2r}[1]\hookrightarrow C^\eta_{2r}[\varepsilon].
\]
For an infinite orbit the same argument uses
$C^\eta_\infty[1]\hookrightarrow C^\eta_\infty[\varepsilon]$.
Well-ordering and attaching all the orbits proves the claim.
\end{proof}

\begin{proposition}\label{prop:tildeRetract}
    Let $X\in\neSSet$. Then $X$ is $\eta\varepsilon$-fibrant if and only if
    \[
        \tilde{\epsilon}_X\colon \widetilde{X}\rightarrow X
    \]
    has the RLP with respect to all monomorphisms in $\neSSet$.
\end{proposition}

\begin{proof}
    We first compare the underlying simplicial morphism of
    $\tilde{\epsilon}_X\colon\widetilde X\to X$
    with that of
    \(
        \widehat\epsilon_{I^*X}\colon
        \widehat{I^*X}\longrightarrow I^*X.
    \)
    We claim that these morphisms are naturally isomorphic.

For $n\geq 0$, by Lemma~\ref{lem:recycling}, we have the following diagram:
    \[
    \begin{tikzcd}[column sep=1.6em]
        {(\widetilde{X})_n}
        \arrow[r,"\cong"]
        \arrow[d,"(\widetilde\epsilon_X)_n"']
        &
        {\Hom(Q^*Q_!I_!\Delta^n,X)}
        \arrow[r,"(\gamma_{\Delta^n})^*"]
        \arrow[d,"(\eta^Q_{I_!\Delta^n})^*"']
        &
        {\Hom(I_!P^*P_!\Delta^n,X)}
        \arrow[d,"(I_!\eta^P_{\Delta^n})^*"']\arrow[r,"\cong"]
        &
        \widehat{I^*X}_n\arrow[d,"(\eta_{\Delta^n}^P)^*"]
        \\
        {X_n}
        \arrow[r,"\cong"]
        &
        {\Hom(I_!\Delta^n,X)}
        \arrow[r,equal]
        &
        {\Hom(I_!\Delta^n,X)}\arrow[r,"\cong"]
        &
        I^*X_n
    \end{tikzcd}
\]
Consequently, by Proposition~\ref{prop:counitIsTriv} we have
    \[
    I_!(\eHorn)\boxslash X
    \quad\Leftrightarrow\quad
    \eHorn\boxslash I^*X
    \quad\Leftrightarrow\quad
    \widehat\epsilon_{I^*X}\text{ is a trivial Kan fibration}
    \]
and the latter is equivalent to $\widetilde{\epsilon}_X$ being a trivial Kan fibration.

    It remains to consider the $\eta\varepsilon$-coherence conditions.
    This is a direct modification of the proof of
    Proposition~\ref{prop:counitIsTriv}(ii).
    One shows that
    the lifting problem
    \[
        (\Delta^{\{0,1\}}\hookrightarrow\Sim^{\{0,1\}})
        \boxslash\widetilde\epsilon_X
    \]
    is equivalent to $X$ having the RLP with respect to
    \[
        \Sim^{\{0,1\}}
        \coprod_{\Delta^{\{0,1\}}}
        C^\eta_\infty[1]
        \hookrightarrow
        C^\eta_\infty[\varepsilon].
    \]
    Similarly,
    \[
        (\Delta^{\{0,1\}}\hookrightarrow \etaSim^1)
        \boxslash\widetilde\epsilon_X
    \]
    is equivalent to $X$ having the RLP with respect to
    \[
        \etaSim^{\{0,1\}}
        \coprod_{\Delta^{\{0,1\}}}
        C^\eta_\infty[1]
        \hookrightarrow
        C^\eta_\infty[\eta].
    \]

    Assuming that $X$ admits fillers for the $\varepsilon$-horns, the
    same argument as in Proposition~\ref{prop:counitIsTriv} shows that
    the $\eta\varepsilon$-coherence conditions for $X$ are equivalent to
    the RLP with respect to these two morphisms. Hence the fibrancy
    conditions for $X$ are equivalent precisely to
    $\widetilde\epsilon_X$ having the RLP with respect to
    all monomorphisms in $\neSSet$.
\end{proof}

We have the following stability result for $\eta\varepsilon$-simplicial sets which is also our main result.
\begin{theorem}\label{thm:mainEta}
    Let $X,Y\in\neSSet$ be $\eta\varepsilon$-fibrant. Then so is $\mathrm{Fun}(X,Y)$.
\end{theorem}

\begin{proof}
The proof follows the same strategy as that of
Theorem~\ref{thm:stability}. By Proposition~\ref{prop:tildeRetract},
$X$ is a retract of $\widetilde X$; hence it is enough to prove that
$\mathrm{Fun}(\widetilde X,Y)$ is $\eta\varepsilon$-fibrant.
By Lemma~\ref{lem:etaCellStructure}, it suffices to show that the following
four classes of pushout-products belong to the weakly saturated class
generated by the $\varepsilon$-horns and the coherence morphisms
(I1--2), (R1--2).

\begin{enumerate}[(1)]
\item
\(
\{\ELambda^n_k\hookrightarrow\Sim^n\mid n\geq1,\;k=0,n\}
\Box
\{\partial C^\eta_r[m]\hookrightarrow C^\eta_r[m]\mid
r\in\overline{\mathbb N},\;m\geq0\};
\)

\item
\(
\{\ELambda^n_k\hookrightarrow\Sim^n\mid n\geq1,\;k=0,n\}
\Box
\{C^\eta_{2r}[1]\hookrightarrow C^\eta_{2r}[\varepsilon]\mid
r\in\overline{\mathbb N}\};
\)

\item
\(
\{\eqref{I1},\eqref{I2},\eqref{R1},\eqref{R2}\}
\Box
\{\partial C^\eta_r[m]\hookrightarrow C^\eta_r[m]\mid
r\in\overline{\mathbb N},\;m\geq0\};
\)

\item
\(
\{\eqref{I1},\eqref{I2},\eqref{R1},\eqref{R2}\}
\Box
\{C^\eta_{2r}[1]\hookrightarrow C^\eta_{2r}[\varepsilon]\mid
r\in\overline{\mathbb N}\}.
\)
\end{enumerate}
By the duality principle, throughout the proof it is enough to consider
$k=0$ and the coherence morphisms (I1) and (R1); the cases $k=n$, (I2),
and (R2) follow by duality.

Case (1) involves no nondegenerate $\eta$-witnesses and follows directly
from Proposition~\ref{prop:key}.

For case (2), the pushout-products are isomorphisms for $n>1$.
It remains to consider
\[
\ELambda^1_0\times C^\eta_{2r}[\varepsilon]
\coprod_{\ELambda^1_0\times C^\eta_{2r}[1]}
\Sim^1\times C^\eta_{2r}[1]
\longrightarrow
\Sim^1\times C^\eta_{2r}[\varepsilon].
\]
The only missing data are the $\varepsilon$-witnesses on the edges
$\langle(0,2i),(1,2i+1)\rangle$ and the $\eta$-witnesses on the edges $\langle(1,2i+1),(1,2i+2)\rangle$, $i\in\mathbb{Z}$; the underlying edges are already present.
For each $i$, the domain contains the following pair of $2$-simplices:
\[
\begin{tikzcd}[row sep=small]
{(1,2i+1)}\arrow[r]&{(1,2i+2)}\\[1em]
{(0,2i-1)}\arrow[u,LA=1.5pt]\arrow[r,dotted]&{(0,2i)}\arrow[u,LA=1.5pt]\arrow[lu]
\end{tikzcd}
\]
We first apply (R1) to attach the missing $\varepsilon$-witness on
$\langle(0,2i),(1,2i+1)\rangle$.
Then we apply (I1) to attach the missing $\eta$-witness on
$\langle(1,2i+1),(1,2i+2)\rangle$.

For case (3), the pushout-product is an isomorphism for $m>0$, so we only
have to consider
\[
A\times C^\eta_r[0]\hookrightarrow B\times C^\eta_r[0],
\]
where $A\hookrightarrow B$ is either (I1) or (R1).

For (I1), the domain misses the $\eta$-witnesses on edges $\langle(1,i),(2,i)\rangle$. We can attach these by (I1), since the domain contains for every $i$ the configuration:
\[
\begin{tikzcd}[row sep=small]
    && {(2,i)} &\\
    {(0,i)}\arrow[rru]\arrow[rrr]&&&{(1,i)}\arrow[lu]\\
    &&{(2,i-1)}\arrow[uu]&\\
    {(0,i-1)}\arrow[uu]\arrow[uuurr,LA=1.5pt]
    \arrow[rrruu,LA=1.5pt]\arrow[rru]\arrow[rrr]
    &&&{(1,i-1)}\arrow[lu]\arrow[uu]
\end{tikzcd}
\]

For (R1), the domain misses the $\varepsilon$-witnesses on the edges $\langle(1,i),(2,i+1)\rangle$. We can attach these by (R1), since the domain contains for every $i$
the configuration:
\[
\begin{tikzcd}[row sep=small]
    && {(2,i+1)} &\\
    {(0,i+1)}\arrow[rru]\arrow[rrr]&&&{(1,i+1)}\arrow[lu]\\
    &&{(2,i)}\arrow[uu]&\\
    {(0,i)}\arrow[uu]\arrow[uuurr,LA=1.5pt]
    \arrow[rrruu]\arrow[rru]\arrow[rrr,dotted]
    &&&{(1,i)}\arrow[luuu]\arrow[lu]\arrow[uu]
\end{tikzcd}
\]

Finally, in case (4) it remains to consider (I1) and (R1). Let $A\hookrightarrow B$ be (I1) or (R1), respectively. Then the domain of
\[
A\times C^\eta_{2r}[\varepsilon]
\coprod_{A\times C^\eta_{2r}[1]}
B\times C^\eta_{2r}[1]
\longrightarrow
B\times C^\eta_{2r}[\varepsilon],
\]
misses the $\eta$-witnesses on the edges $(\langle(1,2i-1),(2,2i)\rangle)$ or the $\varepsilon$-witnesses on the edges $(\langle(1,2i),(2,2i+1)\rangle)$, respectively.
The relevant configurations are, respectively,
\begin{equation}\label{eq:case4}
\begin{tikzcd}[row sep={2.4em,between origins},
               column sep={3.3em,between origins}]
    &{(2,2i)}&\\
    &&\\[8pt]
    {(0,2i-2)}\arrow[ruu,LA=1.5pt]\arrow[rr,LA=1.5pt]
    &&{(1,2i-1)}\arrow[luu]
\end{tikzcd}
\qquad
\begin{tikzcd}[row sep={2.4em,between origins},
               column sep={3.3em,between origins}]
    &{(2,2i+1)}&\\
    &&\\[8pt]
    {(0,2i-1)}\arrow[ruu,LA=1.5pt]\arrow[rr,dotted]
    &&{(1,2i)}\arrow[luu]
\end{tikzcd}.
\end{equation}
Thus, in the first configuration, the missing $\eta$-witness on
$\langle(1,2i-1),(2,2i)\rangle$ can be attached by (I1), while in the
second one the missing $\varepsilon$-witness on
$\langle(1,2i),(2,2i+1)\rangle$ can be attached by (R1).
Attaching these witnesses for all $i$ completes case (4), and hence the proof.
\end{proof}

\begin{corollary}\label{cor:lastCor}
Let $X,Y\in\eSSet$ be $\varepsilon$-fibrant and let
$\mathrm{Fun}^\eta(X,Y)$ be the $\varepsilon$-simplicial subset of $\mathrm{Fun}(X,Y)$
consisting of those simplices whose vertices correspond to
$\eta$-edge-preserving morphisms. Then $\mathrm{Fun}^\eta(X,Y)$ is
$\varepsilon$-fibrant.
\end{corollary}

\begin{proof}
By Lemma~\ref{lem:bridge} and Theorem~\ref{thm:mainEta},
$\mathrm{Fun}(X^\eta,Y^\eta)$ is $\eta\varepsilon$-fibrant. Hence, by
Lemma~\ref{lem:bridge}, it is enough to show
\[
\mathrm{Fun}^\eta(X,Y)\cong I^*\mathrm{Fun}(X^\eta,Y^\eta).
\]

For $Z\in\eSSet$, morphisms $Z\to\mathrm{Fun}^\eta(X,Y)$ naturally correspond
to morphisms $Z\times X\to Y$ such that, for every vertex $z\in Z_0$,
the restriction $\{z\}\times X\to Y$ preserves $\eta$-edges. On the
other hand, by the adjunction $I_!\dashv I^*$ and the tensor-hom
adjunction, morphisms
\[
Z\to I^*\mathrm{Fun}(X^\eta,Y^\eta)
\]
naturally correspond to morphisms
\[
I_!Z\times X^\eta\to Y^\eta.
\]
Since $I_!Z$ has only degenerate $\eta$-witnesses, the latter are
precisely the morphisms $Z\times X\to Y$ satisfying the condition above.
Since this bijection is natural in $Z$, the desired isomorphism follows from the Yoneda lemma.
\end{proof}
\begin{corollary}\label{cor:stabilityUniqueFillers}
    Let $X,Y\in\eSSet$ be $\varepsilon$-fibrant. Assume that $X$ admits the unique right lifting property with respect to $\ELambda^n_i\hookrightarrow \Sim^n$, for $n=1,2$, $i=0,n$. Then $\mathrm{Fun}(X,Y)$ is $\varepsilon$-fibrant.
\end{corollary}
\begin{proof}
    Under the assumptions $X$ does not contain any nontrivial $\eta$-edges; hence $\mathrm{Fun}^\eta(X,Y)=\mathrm{Fun}(X,Y)$.
\end{proof}
\printbibliography

@article{Lachman2026,
  author       = {Lachman, D.},
  title        = {Simplicial Approach to Frobenius Algebras in the Category of Relations},
  journal      = {Applied Categorical Structures},
  year         = {2026},
  volume       = {34},
  pages        = {25},
  doi          = {10.1007/s10485-026-09855-1}
  }

@article{Lachman2025Simplicial,
  author    = {Lachman, Dominik},
  title     = {Simplicial Perspectives on (Pseudo) Effect Algebras},
  year      = {2025},
  note      = {arXiv:2510.02748 [math.CT]}
}

@article{JJL2025,
  author    = {Jen{\v{c}}a, Gejza and Jen{\v{c}}ov{\'a}, Anna and Lachman, Dominik},
  title     = {Coherent Configurations and Frobenius Structures},
  year      = {2025},
  eprint    = {2507.21774},
  archivePrefix = {arXiv},
  primaryClass  = {math.CO},
  doi       = {10.48550/arXiv.2507.21774}
}

@article{MZ, 
author={R. A. Mehta and R. Zhang},
title={Frobenius objects in the category of relations},
journal={Lett Math Phys},
year={2020},
volume={110},
pages={1941–1959},
doi={10.1007/s11005-020-01281-1}
}

@article{CKM,
author = {Contreras, I. and Keller, M. and Mehta, R. A.},
title = {Frobenius objects in the category of spans},
journal = {Reviews in Mathematical Physics},
volume = {34},
number = {10},
year = {2022},
doi = {10.1142/S0129055X22500362}
}

@article{CMS,
title = {Frobenius and commutative pseudomonoids in the bicategory of spans},
journal = {Journal of Geometry and Physics},
volume = {207},
year = {2025},
doi = {10.1016/j.geomphys.2024.105309},
author = {I. Contreras and R. A. Mehta and W. H. Stern}
}

@article{OCS2026,
  author       = {Okay, Cihan and Castillo, Victor
                  and Stern, Walker H.},
  title        = {Simplicial effects and weakly associative
                  partial groups},
  journaltitle = {Applied Categorical Structures},
  volume       = {34},
  date         = {2026},
  eid          = {44},
  doi          = {10.1007/s10485-026-09878-8},
  eprint       = {2502.05958},
  eprinttype   = {arxiv}
}

@article{Stern2021,
  author      = {Stern, W. H.},
  title       = {{2-Segal} objects and algebras in spans},
  journaltitle = {Journal of Homotopy and Related Structures},
  volume      = {16},
  number      = {2},
  date        = {2021},
  pages       = {297--361},
  doi         = {10.1007/s40062-021-00282-8}
}

@unpublished{MM2025Segal,
  author       = {Marx, Sophia E. and Mehta, Rajan Amit},
  title        = {{2-Segal} sets and pseudomonoids
                  in the bicategory of spans},
  date         = {2025},
  eprint       = {2505.22832},
  eprinttype   = {arxiv},
  eprintclass  = {math.CT},
  pubstate     = {prepublished}
}

@unpublished{MM2025TQFT,
  author       = {Marx, Sophia E. and Mehta, Rajan Amit},
  title        = {Coherent span-valued {2D TQFTs}},
  date         = {2025},
  eprint       = {2512.24887},
  eprinttype   = {arxiv},
  eprintclass  = {math.AT},
  pubstate     = {prepublished}
}

@book{DK2019,
  author       = {Dyckerhoff, Tobias and Kapranov, Mikhail},
  title        = {Higher {Segal} Spaces},
  series       = {Lecture Notes in Mathematics},
  number       = {2244},
  publisher    = {Springer},
  location     = {Cham},
  date         = {2019},
  doi          = {10.1007/978-3-030-27124-4}
}

@book{Kock2003,
  author       = {Kock, Joachim},
  title        = {{Frobenius} Algebras and {2D} Topological
                  Quantum Field Theories},
  series       = {London Mathematical Society Student Texts},
  number       = {59},
  publisher    = {Cambridge University Press},
  isbn         = {9780521540315},
  doi          = {10.1017/CBO9780511615443},
  date         = {2003}
}

@thesis{Roumen2017,
  author       = {Roumen, Frank},
  title        = {Effect algebroids},
  type         = {phdthesis},
  institution  = {Radboud University},
  location     = {Nijmegen},
  date         = {2017},
  url          = {https://hdl.handle.net/2066/181438}
}

@misc{PS,
  author = {Pavlovic, D. and Seidel, P. M.},
 title = {(Modular) Effect Algebras are Equivalent to (Frobenius) Antispecial Algebras},
 publisher = {arXiv},
  year = {2016},
  doi = {10.48550/arXiv.1605.06719}
}

@article{HCC2013,
  author       = {Heunen, Chris and Contreras, Ivan
                  and Cattaneo, Alberto S.},
  title        = {Relative {Frobenius} algebras are groupoids},
  journaltitle = {Journal of Pure and Applied Algebra},
  volume       = {217},
  number       = {1},
  date         = {2013},
  pages        = {114--124},
  doi          = {10.1016/j.jpaa.2012.04.002}
}

@misc{Joyal2008,
  author       = {Joyal, Andr{\'e}},
  title        = {The Theory of Quasi-Categories and its Applications},
  year         = {2008},
  howpublished = {Lecture notes, Advanced Course on Simplicial Methods
                  in Higher Categories, Volume II, Centre de Recerca
                  Matem{\`a}tica},
  url          = {https://mat.uab.cat/~kock/crm/hocat/advanced-course/Quadern45-2.pdf}
}

@article{F,
title={Coupled Physical Systems},
author={D. J. Foulis},
journal={Foundations of Physics},
year={1989},
pages={905-922},
volume={19}
}

@article{Gudder1997EffectTestSpacesAndAlgebras,
  author  = {Gudder, Stanley},
  title   = {Effect test spaces and effect algebras},
  journal = {Foundations of Physics},
  volume  = {27},
  number  = {2},
  pages   = {287--304},
  year    = {1997},
  doi     = {10.1007/BF02550455}
}

@article{BM2011,
author={C. Berger and I. Moerdijk},
title={On an extension of the notion of Reedy category},
journal={Math. Z.},
year={2011},
volume={269},
pages={977-1004},
doi={10.1007/s00209-010-0770-x}
}

@article{FPR, 
title={Realism, Operationalism, and Quantum Mechanics},
journal={Found. Phys.},
year={1983},
volume={13},
pages={813-841},
author={D. Foulis and C. Piron and C. Randall}
}

@book{HV, 
author={C. Heunen and J. Vicary},
title={Categories for quantum theory},
publisher={Oxford University Press},
year={2020},
ISBN={978-0198739616}
}

@book{L, 
author={J. Lurie},
title={Higher Topos Theory},
publisher={Princeton University Press},
year={2009},
ISBN={9780691140490}
}

@book{MH, 
author={M. Hovey},
title={Model Categories},
publisher={American Mathematical Society},
year={1998},
ISBN={978-0-8218-4361-1}
}

@article{AHS2007,
  author       = {Ad{\'a}mek, Ji{\v r}{\'i} and H{\'e}bert, Michel and Sousa, Lurdes},
  title        = {A Logic of Injectivity},
  journal = {Journal of Homotopy and Related Structures},
  year         = {2007},
  volume       = {2},
  number       = {2},
  pages        = {13-47},
  url          = {http://eudml.org/doc/223445}
}

@book {H2003,
    AUTHOR = {Hirschhorn, Philip S.},
     TITLE = {Model Categories and Their Localizations},
    SERIES = {Mathematical Surveys and Monographs},
    VOLUME = {99},
 PUBLISHER = {American Mathematical Society},
   ADDRESS = {Providence, RI},
      YEAR = {2003}
}

@article{CaSa2025,
  author  = {Cavallo, Evan and Sattler, Christian},
  title   = {Relative elegance and {C}artesian cubes with one connection},
  journal = {Canadian Journal of Mathematics},
  year    = {2025},
  doi     = {10.4153/S0008414X25101466},
}

@book{Cisinski2019,
  author    = {Cisinski, Denis-Charles},
  title     = {Higher Categories and Homotopical Algebra},
  series    = {Cambridge Studies in Advanced Mathematics},
  number    = {180},
  publisher = {Cambridge University Press},
  location  = {Cambridge},
  date      = {2019},
  doi       = {10.1017/9781108588737},
  isbn      = {9781108473200}
}

\end{document}